\documentclass[12pt]{article}%
\usepackage{amsfonts}
\usepackage{sw20bams}
\usepackage{amsmath}
\usepackage{amssymb}
\usepackage{graphicx}%
\providecommand{\U}[1]{\protect\rule{.1in}{.1in}}
\begin{document}

\title{Errata and Addenda to \textit{Mathematical Constants}}
\author{Steven Finch}
\date{May 27, 2024}
\maketitle

\begin{abstract}
We humbly and briefly offer corrections and supplements to
\textit{Mathematical Constants }(2003) and \textit{Mathematical Constants II
}(2019), both published by Cambridge University Press. \ Comments are always welcome.

\end{abstract}

\footnotetext{Copyright \copyright \ 2024 by Steven R. Finch. All rights
reserved.}

\section{First Volume}

\ \ \ \ \textbf{1.1. Pythagoras' Constant.} A geometric irrationality proof of
$\sqrt{2}$ appears in \cite{Apstl-err}; the transcendence of the numbers
\[%
\begin{array}
[c]{ccccc}%
\sqrt{2}^{\sqrt{2}^{\sqrt{2}}}, &  & i^{i^{i}}, &  & i^{e^{\pi}}%
\end{array}
\]
would follow from a proof of Schanuel's conjecture \cite{MrqSdw-err}. \ A
curious recursion in \cite{GrhmP-err, TStoll-err} gives the $n^{\text{th}}$
digit in the binary expansion of $\sqrt{2}$. Catalan \cite{Ctln-err} proved
the Wallis-like infinite product for $1/\sqrt{2}$. More references on radical
denestings include \cite{DShnks-err, BFHT-err, Zippel-err, SuLnd-err}.

\textbf{1.2. The Golden Mean}. The cubic irrational $\psi=1.3247179572...$\ is
connected to a sequence%
\[%
\begin{array}
[c]{ccccc}%
\psi_{1}=1, &  & \psi_{n}=\sqrt[3]{1+\psi_{n-1}} &  & \text{for }n\geq2
\end{array}
\]
which experimentally gives rise to \cite{Piezas-err}%
\[
\lim_{n\rightarrow\infty}\left(  \psi-\psi_{n}\right)  \left(  3(1+\tfrac
{1}{\psi})\right)  ^{n}=1.8168834242....
\]
The cubic irrational $\chi=1.8392867552...$ is mentioned elsewhere in the
literature with regard to iterative functions \cite{CiMa-err, Lot-err,
Web-err} (the four-numbers game is a special case of what are known as Ducci
sequences), geometric constructions \cite{Sel-err, Rgr-err} and numerical
analysis \cite{Ptra-err}. Infinite radical expressions are further covered in
\cite{Her-err, Lynd-err, Tkch-err}; more generalized continued fractions
appear in \cite{GmzM-err, Murru-err}. See \cite{Fi16-err} for an interesting
optimality property of the logarithmic spiral. A mean-value analog $C$ of
Viswanath's constant $1.13198824...$ (the latter applies for almost every
random Fibonacci sequence) was discovered by Rittaud \cite{Ritt-err}:
$C=1.2055694304...$ has minimal polynomial $x^{3}+x^{2}-x-2$. The Fibonacci
factorial constant $c$ arises in \cite{Ktsn1-err} with regard to the
asymptotics
\begin{align*}
-\frac{d}{ds}%
%TCIMACRO{\dsum \limits_{n=1}^{\infty}}%
%BeginExpansion
{\displaystyle\sum\limits_{n=1}^{\infty}}
%EndExpansion
\frac{1}{f_{n}^{s}}  &  \sim\frac{1}{\ln(\varphi)s^{2}}+\frac{1}{24}\left(
6\ln(5)-2\ln(\varphi)-\frac{3\ln(5)^{2}}{\ln(\varphi)}\right)  +\ln(c)\\
&  \sim\frac{1}{\ln(\varphi)s^{2}}+\ln(0.8992126807...)
\end{align*}
as $s\rightarrow0$, which gives meaning to the \textquotedblleft regularized
product\textquotedblright\ of all Fibonacci numbers.

\textbf{1.3. The Natural Logarithmic Base}. More on the matching problem
appears in \cite{Exp0-err}. Let $N$ denote the number of independent Uniform
$[0,1]$ observations $X_{k}$ necessary until $\sum_{k\leq N}X_{k}$ first
exceeds $1$. The fact that $\operatorname*{E}(N)=e$ goes back to at least
Laplace \cite{Exp1-err}; see also \cite{Exp2-err, Exp3-err, Exp4-err,
Exp5-err, Exp6-err, Exp7-err, Exp8-err, Exp9-err, Exp10-err, Exp11-err}.
Imagine guests arriving one-by-one at an infinitely long dinner table, finding
a seat at random, and choosing a napkin (at the left or at the right) at
random. If there is only one napkin available, then the guest chooses it. The
mean fraction of guests without a napkin is $(2-\sqrt{e})^{2}=0.1233967456...$
and the associated variance is $(3-e)(2-\sqrt{e})^{2}=0.0347631055...$
\cite{Exp12-err, Exp13-err, Exp14-err, Exp15-err}. See pages 280--281 for the
discrete parking problem and \cite{Fi23-err} for related annihilation processes.

Proofs of the two infinite products for $e$ are given in \cite{Ctln-err,
Pip-err}; Hurwitzian continued fractions for $e^{1/q}$ and $e^{2/q}$ appear in
\cite{Davi-err, Walt-err, Matt-err, Poor-err}. The probability that a random
permutation on $n$ symbols is \textit{simple} is asymptotically $1/e^{2}$,
where
\[%
\begin{array}
[c]{ccc}%
(2647513)\text{ is non-simple} &  & \text{(since the interval }2..5\text{ is
mapped onto }4..7\text{),}%
\end{array}
\]%
\[%
\begin{array}
[c]{ccc}%
(2314)\text{ is non-simple} &  & \text{(since the interval }1..2\text{ is
mapped onto }2..3\text{),}%
\end{array}
\]
but $(51742683)$ and $(2413)$ are simple, for example. Only intervals of
length $\ell$, where $1<\ell<n$, are considered, since the lengths $\ell=1 $
and $\ell=n$ are trivial \cite{AAK-err, Brgnl-err}.

Define the following set of integer $k$-tuples
\[
N_{k}=\left\{  (n_{1},n_{2},\ldots,n_{k}):%
%TCIMACRO{\dsum \limits_{j=1}^{k}}%
%BeginExpansion
{\displaystyle\sum\limits_{j=1}^{k}}
%EndExpansion
\frac{1}{n_{j}}=1\text{ and }1\leq n_{1}<n_{2}<\ldots<n_{k}\right\}  .
\]
Martin \cite{Mar-err} proved that
\[
\min_{(n_{1},n_{2},\ldots,n_{k})\in N_{k}}n_{k}\sim\frac{e}{e-1}k
\]
as $k\rightarrow\infty$, but it remains open whether
\[
\max_{(n_{1},n_{2},\ldots,n_{k})\in N_{k}}n_{1}\sim\frac{1}{e-1}k.
\]
Croot \cite{Cro-err} made some progress on the latter:\ He proved that
$n_{1}\geq(1+o(1))k/(e-1)$ for infinitely many values of $k$, and this bound
is best possible.

Holcombe \cite{Hlcmb-err} evaluated the infinite products%
\[%
%TCIMACRO{\dprod \limits_{n=2}^{\infty}}%
%BeginExpansion
{\displaystyle\prod\limits_{n=2}^{\infty}}
%EndExpansion
\left(  1-\frac{1}{n^{2}}\right)  ^{n^{2}}e=\frac{\pi}{e^{3/2}},
\]%
\[%
%TCIMACRO{\dprod \limits_{n=1}^{\infty}}%
%BeginExpansion
{\displaystyle\prod\limits_{n=1}^{\infty}}
%EndExpansion
\left(  1+\frac{1}{n^{2}}\right)  ^{n^{2}}\frac{1}{e}=\frac{\exp\left[
\frac{1}{2}+\frac{2\pi}{3}-\frac{1}{2\pi^{2}}\zeta(3)+\frac{1}{2\pi^{2}%
}\operatorname*{Li}\nolimits_{3}\left(  e^{-2\pi}\right)  +\frac{1}{\pi
}\operatorname*{Li}\nolimits_{2}\left(  e^{-2\pi}\right)  \right]  }%
{2\sinh(\pi)}%
\]
and similar products appear in \cite{KchTzr1-err, KchTzr2-err}. Also, define
$f_{0}(x)=x$ and, for each $n>0$,
\[
f_{n}(x)=\left(  1+f_{n-1}(x)-f_{n-1}(0)\right)  ^{\frac{1}{x}}\text{.}%
\]
This imitates the definition of $e$, in the sense that the exponent
$\rightarrow\infty$ and the base $\rightarrow1$ as $x\rightarrow0$. We have
$f_{1}(0)=e=2.718...$,
\[%
\begin{array}
[c]{ccc}%
f_{2}(0)=\exp\left(  -\frac{e}{2}\right)  =0.257..., &  & f_{3}(0)=\exp\left(
\frac{11-3e}{24}\exp\left(  1-\frac{e}{2}\right)  \right)  =1.086...
\end{array}
\]
and $f_{4}(0)=0.921...$ (too complicated an expression to include here). Does
a pattern develop here?

\textbf{1.4. Archimedes' Constant.} Vi\`ete's product
\[
\frac2\pi=\sqrt{\frac12}\cdot\sqrt{\frac12+\frac12\sqrt{\frac12}}\cdot
\sqrt{\frac12+\frac12\sqrt{\frac12+\frac12\sqrt{\frac12}}}\cdots
\]
has the following close cousin:
\[
\frac2L=\sqrt{\frac12}\cdot\sqrt{\frac12+\frac{\frac12}{\sqrt{\frac12}}}%
\cdot\sqrt{\frac12+\frac{\frac12}{\sqrt{\frac12+\frac{\frac12}{\sqrt{\frac12}%
}}}}\cdots
\]
where $L$ is the lemniscate constant (pages 420--423). Levin \cite{ALvn1-err,
ALvn2-err} developed analogs of sine and cosine for the curve $x^{4}+y^{4}=1$
to prove the latter formula; he also noted that the area enclosed by
$x^{4}+y^{4}=1$ is $\sqrt{2}L$ and that
\[
\frac{2\sqrt{3}}\pi=\left(  \frac12+\sqrt{\frac12}\right)  \cdot\left(
\frac12+\sqrt{\frac12-\frac12\sqrt{\frac12}}\right)  \cdot\left(
\frac12+\sqrt{\frac12-\frac12\sqrt{\frac12-\frac12\sqrt{\frac12}}}\right)
\cdots.
\]
Can the half-circumference of $x^{4}+y^{4}=1$ be written in terms of $L$ as
well? This question makes sense both in the usual $2$-norm and in the
$4$-norm; call the half-circumference $\pi_{4}$ for the latter. More
generally, define $\pi_{p}$ to be the half-circumference of the unit
$p$-circle $|x|^{p}+|y|^{p}=1$, where lengths are measured via the $p$-norm
and $1\leq p<\infty$. It turns out \cite{AdlTnt-err} that $\pi=\pi_{2}$ is the
minimum value of $\pi_{p}$. Additional infinite radical expressions for $\pi$
appear in \cite{Servi-err, Nyblm-err}; more on the Matiyasevich-Guy formula is
covered in \cite{Xpi1-err, Xpi2-err, Xpi3-err, Xpi4-err, Xpi5-err}; see
\cite{JGbb-err} for a revised spigot algorithm for computing decimal digits of
$\pi$ and \cite{Ypi1-err, Ypi2-err} for more on BBP-type formulas.

\textbf{1.5. Euler-Mascheroni Constant.} Impressive surveys appear in
\cite{LgrsEC-err, HavilEC-err, DencEC-err}. \ Subtracting%
\[
1+\frac{1}{2}+\frac{1}{3}+\cdots+\frac{1}{n^{2}}-\ln\left(  n^{2}\right)
\rightarrow\gamma
\]
from%
\[
2\left(  1+\frac{1}{2}+\frac{1}{3}+\cdots+\frac{1}{n}\right)  -2\ln\left(
n\right)  \rightarrow2\gamma
\]
gives a nice result \cite{MacysE-err}%
\[
\lim_{n\rightarrow\infty}\left(
%TCIMACRO{\dsum \limits_{k=1}^{n}}%
%BeginExpansion
{\displaystyle\sum\limits_{k=1}^{n}}
%EndExpansion
\frac{1}{k}-\frac{1}{n+1}-\frac{1}{n+2}-\frac{1}{n+3}-\cdots-\frac{1}{n^{2}%
}\right)  =\gamma
\]
avoiding the usual logarithm. \ De la Vall\'{e}e Poussin's theorem was, in
fact, anticipated by Dirichlet \cite{Drch-err, Dck-err}; it is a corollary of
the formula for the limiting mean value of $d(n)$ \cite{Krmr-err}. Vacca's
series was anticipated by Franklin \cite{Frnk-err}, Nielsen \cite{Niel-err}
and Jacobsthal \cite{Jcb1-err, Jcb2-err}. An extension was found by Koecher
\cite{Koec-err}:
\[
\gamma=\delta-\frac{1}{2}%
%TCIMACRO{\dsum \limits_{k=2}^{\infty}}%
%BeginExpansion
{\displaystyle\sum\limits_{k=2}^{\infty}}
%EndExpansion
\frac{(-1)^{k}}{(k-1)k(k+1)}\left\lfloor \frac{\ln(k)}{\ln(2)}\right\rfloor
\]
where $\delta=(1+\alpha)/4=0.6516737881...$ and $\alpha=%
%TCIMACRO{\tsum \nolimits_{n=1}^{\infty}}%
%BeginExpansion
{\textstyle\sum\nolimits_{n=1}^{\infty}}
%EndExpansion
1/(2^{n}-1)=1.6066951524...$ is one of the digital search tree constants.
Glaisher \cite{Glaish-err} discovered a similar formula:
\[
\gamma=%
%TCIMACRO{\dsum \limits_{n=1}^{\infty}}%
%BeginExpansion
{\displaystyle\sum\limits_{n=1}^{\infty}}
%EndExpansion
\frac{1}{3^{n}-1}-2%
%TCIMACRO{\dsum \limits_{k=1}^{\infty}}%
%BeginExpansion
{\displaystyle\sum\limits_{k=1}^{\infty}}
%EndExpansion
\frac{1}{(3k-1)(3k)(3k+1)}\left\lfloor \frac{\ln(3k)}{\ln(3)}\right\rfloor
\]
nearly eighty years earlier. The following series \cite{So1-err, So2-err,
So0-err} suggest that $\ln(4/\pi)=0.2415644752...$ is an \textquotedblleft
alternating Euler constant\textquotedblright:
\[
\gamma=%
%TCIMACRO{\dsum \limits_{k=1}^{\infty}}%
%BeginExpansion
{\displaystyle\sum\limits_{k=1}^{\infty}}
%EndExpansion
\left(  \frac{1}{k}-\ln\left(  1+\frac{1}{k}\right)  \right)  =-%
%TCIMACRO{\dint \limits_{0}^{1}}%
%BeginExpansion
{\displaystyle\int\limits_{0}^{1}}
%EndExpansion%
%TCIMACRO{\dint \limits_{0}^{1}}%
%BeginExpansion
{\displaystyle\int\limits_{0}^{1}}
%EndExpansion
\frac{1-x}{(1-xy)\ln(xy)}dx\,dy,
\]%
\[
\ln\left(  \frac{4}{\pi}\right)  =%
%TCIMACRO{\dsum \limits_{k=1}^{\infty}}%
%BeginExpansion
{\displaystyle\sum\limits_{k=1}^{\infty}}
%EndExpansion
(-1)^{k-1}\left(  \frac{1}{k}-\ln\left(  1+\frac{1}{k}\right)  \right)  =-%
%TCIMACRO{\dint \limits_{0}^{1}}%
%BeginExpansion
{\displaystyle\int\limits_{0}^{1}}
%EndExpansion%
%TCIMACRO{\dint \limits_{0}^{1}}%
%BeginExpansion
{\displaystyle\int\limits_{0}^{1}}
%EndExpansion
\frac{1-x}{(1+xy)\ln(xy)}dx\,dy
\]
(see section 1.7 later for more). \ An equally compelling analog is
$\operatorname{Ei}(\ln(2))-\gamma=0.4679481152...$, as described in
\cite{Blagou-err}. \ Evaluation of the definite integral involving $%
%TCIMACRO{\tsum \nolimits_{k=1}^{\infty}}%
%BeginExpansion
{\textstyle\sum\nolimits_{k=1}^{\infty}}
%EndExpansion
x^{2^{k}}$ was first done by Catalan \cite{Ctln-err}; related formulas%
\[%
%TCIMACRO{\dint \limits_{0}^{1}}%
%BeginExpansion
{\displaystyle\int\limits_{0}^{1}}
%EndExpansion
\dfrac{1+2x}{1+x+x^{2}}\left(
%TCIMACRO{\dsum \limits_{k=1}^{\infty}}%
%BeginExpansion
{\displaystyle\sum\limits_{k=1}^{\infty}}
%EndExpansion
x^{3^{k}}\right)  dx=1-\gamma=%
%TCIMACRO{\dint \limits_{0}^{1}}%
%BeginExpansion
{\displaystyle\int\limits_{0}^{1}}
%EndExpansion
\dfrac{1+\frac{1}{2}\sqrt{x}}{\left(  1+\sqrt{x}\right)  \left(  1+\sqrt
{x}+x\right)  }\left(
%TCIMACRO{\dsum \limits_{k=1}^{\infty}}%
%BeginExpansion
{\displaystyle\sum\limits_{k=1}^{\infty}}
%EndExpansion
x^{(3/2)^{k}}\right)  dx
\]
are due to Ramanujan \cite{BrnBwm-err}.

Sample criteria for the irrationality of $\gamma$ appear in Sondow
\cite{So3-err, So4-err, So5-err, HPHP-err, Prvst-err}. Long ago, Mahler
attempted to prove that $\gamma$ is transcendental; the closest he came to
this was to prove the transcendentality of the constant \cite{Mah1-err,
Mah2-err}
\[
\frac{\pi Y_{0}(2)}{2J_{0}(2)}-\gamma
\]
where $J_{0}(x)$ and $Y_{0}(x)$ are the zeroth Bessel functions of the first
and second kinds. (Unfortunately the conclusion cannot be applied to the terms
separately!) From Nesterenko's work, $\Gamma(1/6)$ is transcendental; from
Grinspan's work \cite{Murty0-err}, at least two of the three numbers $\pi$,
$\Gamma(1/5)$, $\Gamma(2/5)$ are algebraically independent. See
\cite{Murty1-err, Murty2-err, Murty3-err} for more such results.

Diamond \cite{Dia-err, Zha-err} proved that, if
\[
F_{k}(n)=%
%TCIMACRO{\dsum }%
%BeginExpansion
{\displaystyle\sum}
%EndExpansion
\frac{1}{\ln(\nu_{1})\ln(\nu_{2})\cdots\ln(\nu_{k})}%
\]
where the (finite) sum is over all integer multiplicative compositions
$n=\nu_{1}\nu_{2}\cdots\nu_{k}$ and each $\nu_{j}\geq2$, then
\[
\lim_{N\rightarrow\infty}\frac{1}{N}\left(  1+%
%TCIMACRO{\dsum \limits_{n=2}^{N}}%
%BeginExpansion
{\displaystyle\sum\limits_{n=2}^{N}}
%EndExpansion%
%TCIMACRO{\dsum \limits_{k=1}^{\infty}}%
%BeginExpansion
{\displaystyle\sum\limits_{k=1}^{\infty}}
%EndExpansion
\frac{F_{k}(n)}{k!}\right)  =\exp(\gamma^{\prime}-\gamma-\ln(\ln
(2))=1.2429194164...
\]
where $\gamma^{\prime}=0.4281657248...$ is the analog of Euler's constant when
$1/x$ is replaced by $1/(x\ln(x))$ (see Table 1.1). The analog when $1/x$ is
replaced by $1/\sqrt{x}$ is called Ioachimescu's constant \cite{GvrIvn-err}.
See \cite{BSRK-err} for a different generalization of $\gamma$. \ Also,
related limiting formulas include \cite{Ivan-err}%
\[
\lim_{n\rightarrow\infty}\left(  \sum_{k=1}^{n}\arctan\left(  \frac{1}%
{k}\right)  -\ln(n)\right)  =-\arg\left(  \Gamma(1+i)\right)  ,
\]%
\[
\lim_{n\rightarrow\infty}\left(  \sum_{k=2}^{n}\operatorname{arctanh}\left(
\frac{1}{k}\right)  -\ln(n)\right)  =-\frac{1}{2}\ln(2).
\]

\textbf{1.6. Ap\'{e}ry's Constant.} The famous alternating central binomial
series for $\zeta(3)$ dates back at least as far as 1890, appearing as a
special case of a formula due to Markov \cite{WZ1-err, WZ2-err, WZ3-err}:
\[%
%TCIMACRO{\dsum \limits_{n=0}^{\infty}}%
%BeginExpansion
{\displaystyle\sum\limits_{n=0}^{\infty}}
%EndExpansion
\frac{1}{(x+n)^{3}}=\frac{1}{4}%
%TCIMACRO{\dsum \limits_{n=0}^{\infty}}%
%BeginExpansion
{\displaystyle\sum\limits_{n=0}^{\infty}}
%EndExpansion
\frac{(-1)^{n}(n!)^{6}}{(2n+1)!}\frac{2(x-1)^{2}+6(n+1)(x-1)+5(n+1)^{2}%
}{\left[  x(x+1)\cdots(x+n)\right]  ^{4}}.
\]
Ramanujan \cite{Vpts1-err, Vpts2-err} discovered the series for $\zeta(3)$
attributed to Grosswald. Plouffe \cite{Plff-err} uncovered remarkable formulas
for $\pi^{2k+1}$ and $\zeta(2k+1)$, including
\[
\pi=72%
%TCIMACRO{\dsum \limits_{n=1}^{\infty}}%
%BeginExpansion
{\displaystyle\sum\limits_{n=1}^{\infty}}
%EndExpansion
\frac{1}{n(e^{\pi n}-1)}-96%
%TCIMACRO{\dsum \limits_{n=1}^{\infty}}%
%BeginExpansion
{\displaystyle\sum\limits_{n=1}^{\infty}}
%EndExpansion
\frac{1}{n(e^{2\pi n}-1)}+24%
%TCIMACRO{\dsum \limits_{n=1}^{\infty}}%
%BeginExpansion
{\displaystyle\sum\limits_{n=1}^{\infty}}
%EndExpansion
\frac{1}{n(e^{4\pi n}-1)},
\]%
\[
\pi^{3}=720%
%TCIMACRO{\dsum \limits_{n=1}^{\infty}}%
%BeginExpansion
{\displaystyle\sum\limits_{n=1}^{\infty}}
%EndExpansion
\frac{1}{n^{3}(e^{\pi n}-1)}-900%
%TCIMACRO{\dsum \limits_{n=1}^{\infty}}%
%BeginExpansion
{\displaystyle\sum\limits_{n=1}^{\infty}}
%EndExpansion
\frac{1}{n^{3}(e^{2\pi n}-1)}+180%
%TCIMACRO{\dsum \limits_{n=1}^{\infty}}%
%BeginExpansion
{\displaystyle\sum\limits_{n=1}^{\infty}}
%EndExpansion
\frac{1}{n^{3}(e^{4\pi n}-1)},
\]%
\[
\pi^{5}=7056%
%TCIMACRO{\dsum \limits_{n=1}^{\infty}}%
%BeginExpansion
{\displaystyle\sum\limits_{n=1}^{\infty}}
%EndExpansion
\frac{1}{n^{5}(e^{\pi n}-1)}-6993%
%TCIMACRO{\dsum \limits_{n=1}^{\infty}}%
%BeginExpansion
{\displaystyle\sum\limits_{n=1}^{\infty}}
%EndExpansion
\frac{1}{n^{5}(e^{2\pi n}-1)}+63%
%TCIMACRO{\dsum \limits_{n=1}^{\infty}}%
%BeginExpansion
{\displaystyle\sum\limits_{n=1}^{\infty}}
%EndExpansion
\frac{1}{n^{5}(e^{4\pi n}-1)},
\]%
\[
\zeta(3)=28%
%TCIMACRO{\dsum \limits_{n=1}^{\infty}}%
%BeginExpansion
{\displaystyle\sum\limits_{n=1}^{\infty}}
%EndExpansion
\frac{1}{n^{3}(e^{\pi n}-1)}-37%
%TCIMACRO{\dsum \limits_{n=1}^{\infty}}%
%BeginExpansion
{\displaystyle\sum\limits_{n=1}^{\infty}}
%EndExpansion
\frac{1}{n^{3}(e^{2\pi n}-1)}+7%
%TCIMACRO{\dsum \limits_{n=1}^{\infty}}%
%BeginExpansion
{\displaystyle\sum\limits_{n=1}^{\infty}}
%EndExpansion
\frac{1}{n^{3}(e^{4\pi n}-1)},
\]%
\[
\zeta(5)=24%
%TCIMACRO{\dsum \limits_{n=1}^{\infty}}%
%BeginExpansion
{\displaystyle\sum\limits_{n=1}^{\infty}}
%EndExpansion
\frac{1}{n^{5}(e^{\pi n}-1)}-\frac{259}{10}%
%TCIMACRO{\dsum \limits_{n=1}^{\infty}}%
%BeginExpansion
{\displaystyle\sum\limits_{n=1}^{\infty}}
%EndExpansion
\frac{1}{n^{5}(e^{2\pi n}-1)}-\frac{1}{10}%
%TCIMACRO{\dsum \limits_{n=1}^{\infty}}%
%BeginExpansion
{\displaystyle\sum\limits_{n=1}^{\infty}}
%EndExpansion
\frac{1}{n^{5}(e^{4\pi n}-1)},
\]%
\[
\zeta(7)=\frac{304}{13}%
%TCIMACRO{\dsum \limits_{n=1}^{\infty}}%
%BeginExpansion
{\displaystyle\sum\limits_{n=1}^{\infty}}
%EndExpansion
\frac{1}{n^{7}(e^{\pi n}-1)}-\frac{103}{4}%
%TCIMACRO{\dsum \limits_{n=1}^{\infty}}%
%BeginExpansion
{\displaystyle\sum\limits_{n=1}^{\infty}}
%EndExpansion
\frac{1}{n^{7}(e^{2\pi n}-1)}+\frac{19}{52}%
%TCIMACRO{\dsum \limits_{n=1}^{\infty}}%
%BeginExpansion
{\displaystyle\sum\limits_{n=1}^{\infty}}
%EndExpansion
\frac{1}{n^{7}(e^{4\pi n}-1)}.
\]
A claimed proof that $\zeta(5)$ is irrational awaits confirmation
\cite{YCKim-err}. Volchkov's formula (which is equivalent to the Riemann
hypothesis) was revisited in \cite{SkBrM1-err}; a new criterion
\cite{SkBrM2-err} has the advantage that it involves only integrals of
$\zeta(z)$ taken exclusively along the real axis. \ We mention a certain
alternating double sum \cite{BghTBij-err, Brdhrst-err}%
\begin{align*}%
%TCIMACRO{\dsum \limits_{i=2}^{\infty}}%
%BeginExpansion
{\displaystyle\sum\limits_{i=2}^{\infty}}
%EndExpansion
\sum_{j=1}^{i-1}\frac{(-1)^{i+j}}{i^{3}j}  &  =\frac{\pi^{4}}{180}+\frac
{\pi^{2}}{12}\ln(2)^{2}-\frac{1}{12}\ln(2)^{4}-2\operatorname*{Li}%
\nolimits_{4}\left(  \frac{1}{2}\right) \\
&  =-0.1178759996...
\end{align*}
and wonder about possible generalizations.

\textbf{1.7. Catalan's Constant.} Rivoal \&\ Zudilin \cite{RivZud-err,
Zud0-err} proved that there exist infinitely many integers $k$ for which
$\beta(2k)$ is irrational, and that at least one of the numbers $\beta(2)$,
$\beta(4)$, $\beta(6)$, $\beta(8)$, $\beta(10)$, $\beta(12)$ is irrational.
More double integrals (see section 1.5 earlier) include \cite{Beuk-err,
Adam-err, Zud1-err, So6-err}
\[%
\begin{array}
[c]{ccc}%
\zeta(3)=-\dfrac{1}{2}%
%TCIMACRO{\dint \limits_{0}^{1}}%
%BeginExpansion
{\displaystyle\int\limits_{0}^{1}}
%EndExpansion%
%TCIMACRO{\dint \limits_{0}^{1}}%
%BeginExpansion
{\displaystyle\int\limits_{0}^{1}}
%EndExpansion
\dfrac{\ln(xy)\,dx\,dy}{1-xy}, &  & G=\dfrac{1}{8}%
%TCIMACRO{\dint \limits_{0}^{1}}%
%BeginExpansion
{\displaystyle\int\limits_{0}^{1}}
%EndExpansion%
%TCIMACRO{\dint \limits_{0}^{1}}%
%BeginExpansion
{\displaystyle\int\limits_{0}^{1}}
%EndExpansion
\dfrac{dx\,dy}{(1-xy)\sqrt{x(1-y)}}.
\end{array}
\]
Zudilin \cite{Zud1-err} found the continued fraction expansion
\[
\frac{13}{2G}=7+\dfrac{\left.  1040\right\vert }{\left\vert 10699\right.
}+\dfrac{\left.  42322176\right\vert }{\left\vert 434871\right.  }%
+\dfrac{\left.  15215850000\right\vert }{\left\vert 4090123\right.  }+\cdots.
\]
where the partial numerators and partial denominators are generated according
to the polynomials $(2n-1)^{4}(2n)^{4}(20n^{2}-48n+29)(20n^{2}+32n+13)$ and
$3520n^{6}+5632n^{5}+2064n^{4}-384n^{3}-156n^{2}+16n+7$. \ See \cite{BMW-err}
for more on BBP-type formulas. \ A\ natural trilogarithmic extension of
$G=\operatorname{Im}(\operatorname*{Li}\nolimits_{2}(i))$ is%
\[
\operatorname{Im}\left(  \operatorname*{Li}\nolimits_{3}\left(  \frac{1+i}%
{2}\right)  \right)  =0.5700774070...
\]
which occurs in surprisingly many applications \cite{CLvN-err}, for example,
in evaluating
\[%
\begin{array}
[c]{ccccccc}%
%TCIMACRO{\dint \limits_{0}^{1}}%
%BeginExpansion
{\displaystyle\int\limits_{0}^{1}}
%EndExpansion
\tfrac{\arctan(x)\ln(1\pm x)}{x}{\small dx}, &  &
%TCIMACRO{\dint \limits_{0}^{1}}%
%BeginExpansion
{\displaystyle\int\limits_{0}^{1}}
%EndExpansion
{\small K(x)}\ln{\small (x)dx}, &  &
%TCIMACRO{\dint \limits_{0}^{1}}%
%BeginExpansion
{\displaystyle\int\limits_{0}^{1}}
%EndExpansion
\tfrac{\ln(x)\ln(1\pm x)}{1+x^{2}}{\small dx}, &  &
%TCIMACRO{\dint \limits_{0}^{\pi/4}}%
%BeginExpansion
{\displaystyle\int\limits_{0}^{\pi/4}}
%EndExpansion
\ln{\small (2}\sin{\small (x))}^{2}{\small dx}.
\end{array}
\]
A different dilogarithmic variation of $G$ is found in section 3.10.

\textbf{1.8. Khintchine-L\'evy Constants.} Let $m(n,x)$ denote the number of
partial denominators of $x$ correctly predicted by the first $n$ decimal
digits of $x$. Lochs' result is usually stated as \cite{DjKr-err}
\begin{align*}
\lim_{n\rightarrow\infty}\frac{m(n,x)}n  &  =\frac{6\ln(2)\ln(10)}{\pi^{2}%
}=0.9702701143...\\
\  &  =(1.0306408341...)^{-1}=\left[  (2)(0.5153204170...)\right]  ^{-1}%
\end{align*}
for almost all $x$. In words, an extra $3\%$ in decimal digits delivers the
required partial denominators. The constant $0.51532...$ appears in
\cite{Silva-err} and our entry [2.17]. A\ corresponding Central Limit Theorem
is stated in \cite{Fi19-err, Faivr-err}.

If $x$ is a quadratic irrational, then its continued fraction expansion is
periodic; hence $\lim_{n\rightarrow\infty}M(n,x)$ is easily found and is
algebraic. For example, $\lim_{n\rightarrow\infty}M(n,\varphi)=1$, where
$\varphi$ is the Golden mean. We study the set $\Sigma$ of values
$\lim_{n\rightarrow\infty}\ln(Q_{n})/n$ taken over all quadratic irrationals
$x$ in \cite{Fi15-err}. Additional references include \cite{Choe1-err,
Choe2-err, FsRvl-err}.

\textbf{1.9. Feigenbaum-Coullet-Tresser Constants.} Consider the unique
solution of $\varphi(x)=T_{2}[\varphi](x)$ as pictured in Figure 1.6:
\begin{align*}
\varphi(x)  &  =1-(1.5276329970...)x^{2}+(0.1048151947...)x^{4}\\
&  \ +(0.0267056705...)x^{6}-(0.0035274096...)x^{8}+-\cdots
\end{align*}
The Hausdorff dimension $D$ of the Cantor set $\{x_{k}\}_{k=1}^{\infty
}\subseteq[-1,1]$, defined by $x_{1}=1$ and $x_{k+1}=\varphi(x_{k})$, is known
to satisfy $0.53763<D<0.53854$. This set may be regarded as the simplest of
all strange attractors \cite{Lnfd-err, Gras-err, LePr-err}.

In two dimensions, Kuznetsov \& Sataev \cite{KKSv-err} computed parameters
$\alpha=2.502907875...$, $\beta=1.505318159...$, $\delta=4.669201609...$ for
the map
\[
\left(
\begin{array}
[c]{c}%
x_{n+1}\\
y_{n+1}%
\end{array}
\right)  =\left(
\begin{array}
[c]{c}%
1-c\,x_{n}^{2}\\
1-a\,y_{n}^{2}-b\,x_{n}^{2}%
\end{array}
\right)  ;
\]
$\alpha=1.90007167...$, $\beta=4.00815785...$, $\delta=6.32631925...$ for the
map
\[
\left(
\begin{array}
[c]{c}%
x_{n+1}\\
y_{n+1}%
\end{array}
\right)  =\left(
\begin{array}
[c]{c}%
1-a\,x_{n}^{2}+d\,x_{n}y_{n}\\
1-b\,x_{n}y_{n}%
\end{array}
\right)  ;
\]
and $\alpha=6.565350...$, $\beta=22.120227...$, $\delta=92.431263...$ for the
map
\[
\left(
\begin{array}
[c]{c}%
x_{n+1}\\
y_{n+1}%
\end{array}
\right)  =\left(
\begin{array}
[c]{c}%
a-x_{n}^{2}+b\,y_{n}\\
e\,y_{n}-x_{n}^{2}%
\end{array}
\right)  .
\]
``Certainly, this is only a little part of some great entire pattern'', they wrote.

Let us return to the familiar one-dimensional map $x\mapsto a\,x(1-x)$, but
focus instead on the region $a>a_{\infty}=3.5699456718...=4(0.8924864179...).$
We are interested in bifurcation of cycles whose periods are odd multiples of
two:
\[
\lambda(m,n)=%
\begin{array}
[c]{c}%
\text{the smallest value of }a\text{ for which a cycle of}\\
\text{period }(2m+1)2^{n}\text{ first appears.}%
\end{array}
\]
For any fixed $m\geq0$,
\[
\lim_{n\rightarrow\infty}\frac{\lambda(m,n)-\lambda(m,n-1)}{\lambda
(m,n+1)-\lambda(m,n)}=\delta=4.6692...
\]
which is perhaps unsurprising. A\ new constant emerges if we reverse the roles
of $m$ and $n$:
\[
\lim_{n\rightarrow\infty}\underset{\gamma_{n}}{\underbrace{\lim_{m\rightarrow
\infty}\frac{\lambda(m,n)-\lambda(m-1,n)}{\lambda(m+1,n)-\lambda(m,n)}}%
}=\gamma=2.9480...
\]
due to Geisel \&\ Nierwetberg \cite{GNwb-err} and Kolyada \& Sivak
\cite{KydS-err}. High-precision values of $\gamma_{0}$, $\gamma_{1}$,
$\gamma_{2}$, $\ldots$ would be good to see. A proof of the existence of
$\gamma$ is in \cite{SKSF-err}, but apart from mention in \cite{RPMt-err},
this constant has been unjustly neglected.

\textbf{1.10. Madelung's Constant.} The following \textquotedblleft near
miss\textquotedblright\ exact expression \cite{Tyagi-err}:
\begin{align*}
M_{3}  &  =-\frac{1}{8}+\frac{1}{2\sqrt{2}}-\frac{4\pi}{3}-\frac{\ln(2)}{4\pi
}+\frac{\Gamma(1/8)\Gamma(3/8)}{\pi^{3/2}\sqrt{2}}\\
&  -2%
%TCIMACRO{\dsum \limits_{i,j,k=-\infty}^{\infty}}%
%BeginExpansion
{\displaystyle\sum\limits_{i,j,k=-\infty}^{\infty}}
%EndExpansion
{}^{\prime}\frac{(-1)^{i+j+k}}{\sqrt{i^{2}+j^{2}+k^{2}}\left(  e^{8\pi
\sqrt{i^{2}+j^{2}+k^{2}}}-1\right)  }%
\end{align*}
is noteworthy because the series portion is rapidly convergent. See also
\cite{Knmt1-err, Knmt2-err, Knmt3-err}. \ Related to our function $f(z)$ is
the limit%
\[%
%TCIMACRO{\dsum \limits_{i,j=-n}^{n}}%
%BeginExpansion
{\displaystyle\sum\limits_{i,j=-n}^{n}}
%EndExpansion
{}^{\prime}\frac{1}{i^{2}+j^{2}}-2\pi\ln(n)\rightarrow\left[  4\ln(2)+3\ln
(\pi)+2\gamma-4\ln\left(  \Gamma(1/4)\right)  \right]  \pi-4G
\]
as $n\rightarrow\infty$, where $\gamma$ is Euler's constant and $G$ is
Catalan's constant \cite{Scott-err}. \ Another series%
\[%
%TCIMACRO{\dsum \limits_{i,j=-\infty}^{\infty}}%
%BeginExpansion
{\displaystyle\sum\limits_{i,j=-\infty}^{\infty}}
%EndExpansion
\frac{(-1)^{i+j}}{i^{2}+(3j+1)^{2}}=\frac{2\pi}{9}\ln\left[  2\left(  \sqrt
{3}-1\right)  \right]
\]
is only the first of many evaluations appearing in \cite{ZMcPC-err,
BLmbR-err}. \ Likewise%
\[
-%
%TCIMACRO{\dsum \limits_{i,j=-n}^{n}}%
%BeginExpansion
{\displaystyle\sum\limits_{i,j=-n}^{n}}
%EndExpansion
{}^{\prime}\ln\left(  i^{2}+j^{2}\right)  +%
%TCIMACRO{\dint \limits_{x,y=-n-\frac{1}{2}}^{n+\frac{1}{2}}}%
%BeginExpansion
{\displaystyle\int\limits_{x,y=-n-\frac{1}{2}}^{n+\frac{1}{2}}}
%EndExpansion
\ln\left(  x^{2}+y^{2}\right)  \,dx\,dy\rightarrow\ln\left(  \frac{2}{\pi
}\right)  -2\ln\left(  \frac{\Gamma(1/4)}{\Gamma(3/4)}\right)  +\frac{\pi}%
{6},
\]%
\[%
%TCIMACRO{\dsum \limits_{k=1}^{\infty}}%
%BeginExpansion
{\displaystyle\sum\limits_{k=1}^{\infty}}
%EndExpansion
(-1)^{k+1}\frac{\ln(2k+1)}{2k+1}=\frac{\pi}{4}\left\{  \gamma+\ln(2\pi
)-2\ln\left(  \frac{\Gamma(1/4)}{\Gamma(3/4)}\right)  \right\}  ,
\]%
\[%
%TCIMACRO{\dsum \limits_{k=0}^{\infty}}%
%BeginExpansion
{\displaystyle\sum\limits_{k=0}^{\infty}}
%EndExpansion
\left\{  \frac{\ln(3k+1)}{3k+1}-\frac{\ln(3k+2)}{3k+2}\right\}  =\frac{\pi
}{\sqrt{3}}\left\{  \ln\left(  \frac{\Gamma(1/3)}{\Gamma(2/3)}\right)
-\frac{1}{3}\left(  \gamma+\ln(2\pi)\right)  \right\}  ,
\]%
\[%
%TCIMACRO{\dsum \limits_{k=0}^{\infty}}%
%BeginExpansion
{\displaystyle\sum\limits_{k=0}^{\infty}}
%EndExpansion
(-1)^{k}\left\{  \frac{\ln(4k+1)}{4k+1}+\frac{\ln(4k+3)}{4k+3}\right\}
=\frac{\pi}{2\sqrt{2}}\left\{  \ln\left(  \frac{\Gamma(1/8)\Gamma(3/8)}%
{\Gamma(5/8)\Gamma(7/8)}\right)  -\left(  \gamma+\ln(2\pi)\right)  \right\}
\]
are just starting points for research reported in \cite{Shail1-err,
Shail2-err, Shail3-err, Mlmstn1-err, Mlmstn2-err}.

\textbf{1.11. Chaitin's Constant.} Ord \&\ Kieu \cite{OrKi-err} gave a
different Diophantine representation for $\Omega$; apparently Chaitin's
equation can be reduced to 2--3 pages in length \cite{Mtys-err}. A\ rough
sense of the type of equations involved can be gained from \cite{JSWW-err}.
Calude \&\ Stay \cite{CaSt-err} suggested that the uncomputability of bits of
$\Omega$ can be recast as an uncertainty principle.

\textbf{2.1. Hardy-Littlewood Constants.} In a breakthrough, Zhang
\cite{YZhng1-err, YZhng2-err, YZhng3-err, JMynrd-err} proved that the sequence
of gaps between consecutive primes has a finite liminf (an impressive step
toward confirming the Twin Prime Conjecture). In another breakthrough, Green
\&\ Tao \cite{GrnTao-err} proved that there are arbitrarily long arithmetic
progressions of primes. In particular, the number of prime triples
$p_{1}<p_{2}<p_{3}\leq x$ in arithmetic progression is
\[
\sim\frac{C_{\text{twin}}}{2}\frac{x^{2}}{\ln(x)^{3}}=(0.3300809079...)\frac
{x^{2}}{\ln(x)^{3}}%
\]
as $x\rightarrow\infty$, and the number of prime quadruples $p_{1}<p_{2}%
<p_{3}<p_{4}\leq x$ in arithmetic progression is likewise
\[
\sim\frac{D}{6}\frac{x^{2}}{\ln(x)^{4}}=(0.4763747659...)\frac{x^{2}}%
{\ln(x)^{4}}.
\]

Here is a different extension $C_{\text{twin}}=C_{2}^{\prime}$:
\[
P_{n}(p,p+2r)\sim2\underset{C_{2r}^{\prime}}{\underbrace{C_{\text{twin}}%
%TCIMACRO{\dprod \limits_{\substack{p|r \\p>2}}}%
%BeginExpansion
{\displaystyle\prod\limits_{\substack{p|r \\p>2}}}
%EndExpansion
\frac{p-1}{p-2}}}\frac{n}{\ln(n)^{2}},
\]
and $C_{2r}^{\prime}$ has mean value one in the sense that $\sum_{r=1}%
^{m}C_{2r}^{\prime}\sim m$ as $m\rightarrow\infty$. Further generalization is
possible \cite{BltKvr-err, Mathr4-err}.

Fix $\varepsilon>0$. Let $N(x,k)$ denote the number of positive integers
$n\leq x$ with $\Omega(n)=k$, where $k$ is allowed to grow with $x$. Nicolas
\cite{Nic1-err} proved that
\[
\lim_{x\rightarrow\infty}\frac{N(x,k)}{(x/2^{k})\ln(x/2^{k})}=\frac
1{4C_{\text{twin}}}=\frac14%
%TCIMACRO{\dprod \limits_{p>2}}%
%BeginExpansion
{\displaystyle\prod\limits_{p>2}}
%EndExpansion
\left(  1+\frac1{p(p-2)}\right)  =0.3786950320....
\]
under the assumption that $(2+\varepsilon)\ln(\ln(x))\leq k\leq\ln(x)/\ln(2)$.
More relevant results appear in \cite{Hwa1-err}; see also the next entry.

Let $L(x)$ denote the number of positive odd integers $n\leq x$ that can be
expressed in the form $2^{l}+p$, where $l$ is a positive integer and $p$ is a
prime. Then $0.09368\,x\leq L(x)<0.49095\,x$ for all sufficiently large $x$.
The lower bound can be improved to $0.2893\,x$ if the Hardy-Littlewood
conjectures in sieve theory are true \cite{XYZ1-err, XYZ2-err, XYZ3-err,
XYZ4-err, XYZ5-err}.

Let $Q(x)$ denote the number of integers $\leq x$ with prime factorizations
$p_{1}^{\alpha_{1}}p_{2}^{\alpha_{2}}\cdots p_{r}^{\alpha_{r}}$ satisfying
$\alpha_{1}\geq\alpha_{2}\geq\ldots\geq\alpha_{r}$. Extending results of Hardy
\&\ Ramanujan \cite{Hrdy0-err}, Richmond \cite{Rcmd0-err} deduced that
\[
\ln(Q(x))\sim\tfrac{2\pi}{\sqrt{3}}\left(  \tfrac{\ln(x)}{\ln(\ln(x))}\right)
^{1/2}\left(  1-\tfrac{2\ln(\pi)+12B/\pi^{2}-2}{2\ln(\ln(x))}-\tfrac
{\ln(3)-\ln(\ln(\ln(x)))}{2\ln(\ln(x))}\right)
\]
where
\[
B=-%
%TCIMACRO{\dint \limits_{0}^{\infty}}%
%BeginExpansion
{\displaystyle\int\limits_{0}^{\infty}}
%EndExpansion
\ln(1-e^{-y})\ln(y)\,dy=\zeta^{\prime}(2)-\tfrac{\pi^{2}}{6}\gamma.
\]
The Bateman-Horn conjecture arises unexpectedly in \cite{ChMgKo-err}. \ The
ternary Goldbach conjecture ($G^{\prime}$), finally, is proved
\cite{Hlfgtt-err}.

\textbf{2.2. Meissel-Mertens Constants.} See \cite{DKT-err} for more
occurrences of the constants $M$ and $M^{\prime}$, and \cite{Villr-err} for a
historical treatment. Higher-order asymptotic series for $\operatorname*{E}%
\nolimits_{n}(\omega)$, $\operatorname*{Var}\nolimits_{n}(\omega)$,
$\operatorname*{E}\nolimits_{n}(\Omega)$ and $\operatorname*{Var}%
\nolimits_{n}(\Omega)$ are given in \cite{Fi1-err}. The values $m_{1,3}%
=-0.3568904795...$ and $m_{2,3}=0.2850543590...$ are calculated in
\cite{Fi18-err}; of course, $m_{1,3}+m_{2,3}+1/3=M$. While $%
%TCIMACRO{\tsum \nolimits_{p}}%
%BeginExpansion
{\textstyle\sum\nolimits_{p}}
%EndExpansion
1/p$ is divergent, the following prime series is convergent \cite{HChn-err}:
\[%
%TCIMACRO{\dsum \limits_{p}}%
%BeginExpansion
{\displaystyle\sum\limits_{p}}
%EndExpansion
\left(  \frac{1}{p^{2}}+\frac{1}{p^{3}}+\frac{1}{p^{4}}+\cdots\right)  =%
%TCIMACRO{\dsum \limits_{p}}%
%BeginExpansion
{\displaystyle\sum\limits_{p}}
%EndExpansion
\frac{1}{p(p-1)}=0.7731566690....
\]
The same is true if we replace primes by semiprimes \cite{Mathr1-err}:
\[%
%TCIMACRO{\dsum \limits_{p,q}}%
%BeginExpansion
{\displaystyle\sum\limits_{p,q}}
%EndExpansion%
%TCIMACRO{\dsum \limits_{k=2}^{\infty}}%
%BeginExpansion
{\displaystyle\sum\limits_{k=2}^{\infty}}
%EndExpansion
\frac{1}{(pq)^{k}}=%
%TCIMACRO{\dsum \limits_{p,q}}%
%BeginExpansion
{\displaystyle\sum\limits_{p,q}}
%EndExpansion
\frac{1}{pq(pq-1)}=0.1710518929....
\]
Also, the reciprocal sum of semiprimes satisfies \cite{Saidak-err, Popa-err}
\[
\lim_{n\rightarrow\infty}\left(
%TCIMACRO{\dsum \limits_{pq\leq n}}%
%BeginExpansion
{\displaystyle\sum\limits_{pq\leq n}}
%EndExpansion
\frac{1}{pq}-\ln(\ln(n))^{2}-2M\ln(\ln(n))\right)  =\frac{\pi^{2}}{6}+M^{2}%
\]
and the corresponding analog of Mertens' product formula is
\[
\lim_{n\rightarrow\infty}\left(  \ln(n)\right)  ^{\ln(\ln(n))+2M}%
%TCIMACRO{\dprod \limits_{pq\leq n}}%
%BeginExpansion
{\displaystyle\prod\limits_{pq\leq n}}
%EndExpansion
\left(  1-\frac{1}{pq}\right)  =e^{-\pi^{2}/6-M^{2}-\Lambda}%
\]
where \cite{Mathr1-err}
\[
\Lambda=%
%TCIMACRO{\dsum \limits_{p,q}}%
%BeginExpansion
{\displaystyle\sum\limits_{p,q}}
%EndExpansion%
%TCIMACRO{\dsum \limits_{k=2}^{\infty}}%
%BeginExpansion
{\displaystyle\sum\limits_{k=2}^{\infty}}
%EndExpansion
\frac{1}{k\,(pq)^{k}}=-%
%TCIMACRO{\dsum \limits_{p,q}}%
%BeginExpansion
{\displaystyle\sum\limits_{p,q}}
%EndExpansion
\left(  \ln\left(  1-\frac{1}{pq}\right)  +\frac{1}{pq}\right)
=0.0798480403....
\]
We can think of $\pi^{2}/6+M^{2}+\Lambda$ as another two-dimensional
generalization of Euler's constant $\gamma$.

The second moment of $\operatorname*{Im}(\ln(\zeta(1/2+i\,t)))$ over an
interval $[0,T]$ involves asymptotically a constant \cite{Gldsn-err,
TszHo-err}
\[%
%TCIMACRO{\dsum \limits_{m=2}^{\infty}}%
%BeginExpansion
{\displaystyle\sum\limits_{m=2}^{\infty}}
%EndExpansion%
%TCIMACRO{\dsum \limits_{p}}%
%BeginExpansion
{\displaystyle\sum\limits_{p}}
%EndExpansion
\left(  \frac{1}{m}-\frac{1}{m^{2}}\right)  \frac{1}{p^{m}}=-%
%TCIMACRO{\dsum \limits_{p}}%
%BeginExpansion
{\displaystyle\sum\limits_{p}}
%EndExpansion
\left(  \ln\left(  1-\frac{1}{p}\right)  +\operatorname*{Li}\nolimits_{2}%
\left(  \frac{1}{p}\right)  \right)  =0.1762478124...
\]
as $T\rightarrow\infty$. \ This assumes, however, that a certain random matrix
model is applicable (asymptotics for the pair correlation of zeros).

If $Q_{k}$ denotes the set of positive integers $n$ for which $\Omega
(n)-\omega(n)=k$, then $Q_{1}=\tilde{S}$ and the asymptotic density
$\delta_{k}$ satisfies \cite{ARny-err, MKac-err, KMrs-err}
\[
\lim_{k\rightarrow\infty}2^{k}\delta_{k}=\frac{1}{4C_{\text{twin}}%
}=0.3786950320...;
\]
the expression $4C_{\text{twin}}$ also appears on pages 86 and 133--134, as
well as in the preceding entry.

Given a positive integer $n$, let $K(n)=%
%TCIMACRO{\tprod \nolimits_{p\mid n}}%
%BeginExpansion
{\textstyle\prod\nolimits_{p\mid n}}
%EndExpansion
p$ denote the square-free kernel of $n$ and $\rho_{n}=n/K(n)$. \ We say that
$n$ is \textbf{flat} if the ratio $\rho_{n}=1$. Define $R_{k}$ to be the set
of $n$ such that $\rho_{n}$ itself is flat and $\omega(\rho_{n})=k$. \ We have
$R_{1}=$\ $\tilde{S}$ and asymptotic densities for $R_{2}$, $R_{3}$ equal to
\cite{AdRvL-err}%
\[
\frac{6}{\pi^{2}}%
%TCIMACRO{\dsum \limits_{p<q}}%
%BeginExpansion
{\displaystyle\sum\limits_{p<q}}
%EndExpansion
\frac{1}{p(p+1)q(q+1)}=0.0221245744...,
\]%
\[
\frac{6}{\pi^{2}}%
%TCIMACRO{\dsum \limits_{p<q<r}}%
%BeginExpansion
{\displaystyle\sum\limits_{p<q<r}}
%EndExpansion
\frac{1}{p(p+1)q(q+1)r(r+1)}=0.0010728279....
\]
Averaging $\rho_{n}$ over all $n\geq1$ remains unsolved \cite{Fi34-err}.

Define $f_{k}(n)=\#\{p:p^{k}|n\}$ and $F_{k}(n)=\#\{p^{k+m}:p^{k+m}|n$ and
$m\geq0\}$; hence $f_{1}(n)=\omega(n)$ and $F_{1}(n)=\Omega(n)$. It is known
that, for $k\geq2$,
\[%
\begin{array}
[c]{ccc}%
%TCIMACRO{\dsum \limits_{n\leq x}}%
%BeginExpansion
{\displaystyle\sum\limits_{n\leq x}}
%EndExpansion
f_{k}(n)\sim x%
%TCIMACRO{\dsum \limits_{p}}%
%BeginExpansion
{\displaystyle\sum\limits_{p}}
%EndExpansion
\dfrac1{p^{k}}, &  &
%TCIMACRO{\dsum \limits_{n\leq x}}%
%BeginExpansion
{\displaystyle\sum\limits_{n\leq x}}
%EndExpansion
F_{k}(n)\sim x%
%TCIMACRO{\dsum \limits_{p}}%
%BeginExpansion
{\displaystyle\sum\limits_{p}}
%EndExpansion
\dfrac1{p^{k-1}(p-1)}%
\end{array}
\]
as $x\rightarrow\infty$. Also define $g_{k}(n)=\#\{p:p|n$ and $p^{k}\nmid n\}$
and $G_{k}(n)=\#\{p^{m}:p^{m}|n$, $p^{k}\nmid n$ and $m\geq1\}$. Then, for
$k\geq2$,
\[%
%TCIMACRO{\dsum \limits_{n\leq x}}%
%BeginExpansion
{\displaystyle\sum\limits_{n\leq x}}
%EndExpansion
g_{k}(n)\sim x\left(  \ln(\ln(x))+M-%
%TCIMACRO{\dsum \limits_{p}}%
%BeginExpansion
{\displaystyle\sum\limits_{p}}
%EndExpansion
\dfrac1{p^{k}}\right)  ,
\]
\[%
%TCIMACRO{\dsum \limits_{n\leq x}}%
%BeginExpansion
{\displaystyle\sum\limits_{n\leq x}}
%EndExpansion
G_{k}(n)\sim x\left(  \ln(\ln(x))+M+%
%TCIMACRO{\dsum \limits_{p}}%
%BeginExpansion
{\displaystyle\sum\limits_{p}}
%EndExpansion
\dfrac{p^{k-1}-kp+k-1}{p^{k}(p-1)}\right)
\]
as $x\rightarrow\infty$. Other variations on $k$-full and $k$-free prime
factors appear in \cite{JGrh-err}; the growth rate of $\sum_{n\leq x}%
1/\omega(n)$ and $\sum_{n\leq x}1/\Omega(n)$ is covered in \cite{Knnck-err} as well.

\textbf{2.3. Landau-Ramanujan Constant.} It is not hard to show that
$C_{2}=0.6093010224...$ \cite{Fi28-err}. The second-order constant
corresponding to non-hypotenuse numbers should be
\[
\tilde{C}=C+\frac{1}{2}\ln\left(  \frac{\pi^{2}e^{\gamma}}{2L^{2}}\right)
=0.7047534517...
\]
(numerically unchanged, but $\pi$ is replaced by $\pi^{2}$). Moree
\cite{PMoree4-err} expressed such constants somewhat differently:%
\[%
\begin{array}
[c]{ccc}%
1-2C=-0.1638973186..., &  & 1-2\tilde{C}=-0.4095069034...
\end{array}
\]
calling these Euler-Kronecker constants. \ His terminology is unfortunately
inconsistent with ours \cite{Fi32-err, Fi33-err}. \ 

Define $B_{3,j}(x)$ to be the number of positive integers $\leq x$, all of
whose prime factors are $\equiv j\operatorname*{mod}3$, where $j=1$ or $2$. We
have \cite{ELandau-err, EWirsing-err, PMoree0-err}
\[
\lim_{x\rightarrow\infty}\frac{\sqrt{\ln(x)}}{x}B_{3,1}(x)=\frac{\sqrt{3}%
}{9K_{3}}=0.3012165544...,
\]%
\[
\lim_{x\rightarrow\infty}\frac{\sqrt{\ln(x)}}{x}B_{3,2}(x)=\frac{2\sqrt
{3}K_{3}}{\pi}=0.7044984335....
\]
An analog of Mertens' theorem for primes $\equiv j\operatorname*{mod}3$
unsurprisingly involves $K_{3}$ as well \cite{Fi18-err}. Here is a more
complicated example (which arises in the theory of partitions). Let
\[
W(x)=\#\left\{  n\leq x:n=2^{h}p_{1}^{e_{1}}p_{2}^{e_{2}}\cdots p_{h}^{e_{h}%
},\;h\geq1,\;e_{k}\geq1,\;p_{k}\equiv3,5,6\operatorname*{mod}7\text{ for all
}k\right\}  ,
\]
then the Selberg-Delange method gives \cite{BsN1-err, BsN2-err}
\begin{align*}
\lim_{x\rightarrow\infty}\frac{\ln(x)^{3/4}}{x}W(x)  &  =\frac{1}{\Gamma
(1/4)}\left(  \frac{6}{\sqrt{7}\pi}\right)  ^{1/4}%
%TCIMACRO{\dprod \limits_{\substack{p\equiv3,5,6 \\\operatorname*{mod}7}}}%
%BeginExpansion
{\displaystyle\prod\limits_{\substack{p\equiv3,5,6 \\\operatorname*{mod}7}}}
%EndExpansion
\left(  1+\frac{1}{2(p-1)}\right)  \left(  1-\frac{1}{p}\right)  ^{1/4}\left(
1+\frac{1}{p}\right)  ^{-1/4}\\
\  &  =\frac{1}{\Gamma(1/4)}\left(  \frac{6}{\sqrt{7}\pi}\right)
^{1/4}(1.0751753443...)=0.2733451113...\\
\  &  =\frac{7}{24}(0.9371832387...).
\end{align*}
Other examples appear in \cite{BsN2-err} as well.

Define $Z_{3,j}(x)$ to be the number of positive integers $n\leq x$ for which
$\varphi(n)$ $\equiv j\operatorname*{mod}3$, where $\varphi$ is Euler's
totient and $j=1$ or $2$. We have \cite{ETnt1-err, ETnt2-err}%
\[
\lim_{x\rightarrow\infty}\frac{\sqrt{\ln(x)}}{x}Z_{3,j}(x)=\frac{\sqrt
{2\sqrt{3}}}{3\pi}\frac{2\xi+(-1)^{j+1}\eta}{\xi^{1/2}}=\left\{
\begin{array}
[c]{ccc}%
0.6109136202... &  & \text{if }j=1,\\
0.3284176245... &  & \text{if }j=2
\end{array}
\right.
\]
where%
\[
\xi=%
%TCIMACRO{\dprod \limits_{p\equiv2\operatorname*{mod}3}}%
%BeginExpansion
{\displaystyle\prod\limits_{p\equiv2\operatorname*{mod}3}}
%EndExpansion
\left(  1+\frac{1}{p^{2}-1}\right)  =1.4140643909...,
\]%
\[
\eta=%
%TCIMACRO{\dprod \limits_{p\equiv2\operatorname*{mod}3}}%
%BeginExpansion
{\displaystyle\prod\limits_{p\equiv2\operatorname*{mod}3}}
%EndExpansion
\left(  1-\frac{1}{(p+1)^{2}}\right)  =0.8505360177....
\]
Analogous results for $Z_{4,j}(x)$ with $j=1$ or $3$ are open, as far as is known.

Estermann \cite{Sqfr1-err, Sqfr2-err, Sqfr3-err} first examined the
asymptotics
\[
\hat{B}(x)=%
%TCIMACRO{\dsum \limits_{1\leq m\leq x}}%
%BeginExpansion
{\displaystyle\sum\limits_{1\leq m\leq x}}
%EndExpansion
\mu\left(  m^{2}+1\right)  ^{2}\sim\hat{K}\,x=(0.8948412245...)x
\]
as $x\rightarrow\infty$, where $\mu$ is the M\"{o}bius mu function. \ One
possible generalization is \cite{Sqfr4-err}%
\[%
%TCIMACRO{\dsum \limits_{1\leq m,n\leq x}}%
%BeginExpansion
{\displaystyle\sum\limits_{1\leq m,n\leq x}}
%EndExpansion
\mu\left(  m^{2}+n^{2}+1\right)  ^{2}\sim\hat{J}\,x^{2}%
\]
and a numerical value for $\hat{J}$ evidently remains open. \ See
\cite{Fi30-err} for another occurrence of $\hat{K}$.

Fix $h\geq2$. Define $N_{h}(x)$ to be the number of positive integers not
exceeding $x$ that can be expressed as a sum of two nonnegative integer
$h^{\text{th}}$ powers. Clearly $N_{2}(x)=B(x)$. Hooley \cite{Hoo1-err,
Hoo2-err} proved that
\[
\lim_{x\rightarrow\infty}x^{-2/h}N_{h}(x)=\frac1{4\,h}\frac{\Gamma(1/h)^{2}%
}{\Gamma(2/h)}%
\]
when $h$ is an odd prime, and Greaves \cite{Greav-err} proved likewise when
$h$ is the smallest composite $4$. It is possible that such asymptotics are
true for larger composites, for example, $h=6$.

While $N_{2}(x)$ also counts $n\leq x$ that can be expressed as a sum of two
\textit{rational} squares, it is not true that $N_{3}(x)$ does likewise for
sums of two rational cubes. See \cite{Fi26-err} for analysis of a related
family of elliptic curves (cubic twists of the Fermat equation $u^{3}+v^{3}%
=1$) and \cite{Fi27-err} for an unexpected appearance of the constant $K$.

The issue regarding counts of $x$ of the form $a^{3}+2\,b^{3}$ is addressed in
\cite{Hoo3-err}. We mention that products like \cite{MrtaChin-err}
\[%
%TCIMACRO{\dprod \limits_{p\equiv3\operatorname*{mod}4}}%
%BeginExpansion
{\displaystyle\prod\limits_{p\equiv3\operatorname*{mod}4}}
%EndExpansion
\left(  1-\frac{2p}{(p^{2}+1)(p-1)}\right)  =0.6436506796...,
\]
\[%
%TCIMACRO{\dprod \limits_{p\equiv2\operatorname*{mod}3}}%
%BeginExpansion
{\displaystyle\prod\limits_{p\equiv2\operatorname*{mod}3}}
%EndExpansion
\left(  1-\frac{2p}{(p^{2}+1)(p-1)}\right)  =0.1739771224...
\]
are evaluated to high precision in \cite{PMoree1-err, PMoree2-err} via special
values of Dirichlet L-series.

\textbf{2.4. Artin's Constant.} Other representations include
\cite{Chrdnk-err}
\[
\lim_{N\rightarrow\infty}\frac{\ln(N)}{N}%
%TCIMACRO{\dsum \limits_{p\leq N}}%
%BeginExpansion
{\displaystyle\sum\limits_{p\leq N}}
%EndExpansion
\frac{\varphi(p-1)}{p-1}=C_{\text{Artin}}=\lim_{N\rightarrow\infty}\frac{%
%TCIMACRO{\dsum \limits_{p\leq N}}%
%BeginExpansion
{\displaystyle\sum\limits_{p\leq N}}
%EndExpansion
\varphi(p-1)}{%
%TCIMACRO{\dsum \limits_{p\leq N}}%
%BeginExpansion
{\displaystyle\sum\limits_{p\leq N}}
%EndExpansion
(p-1)}.
\]
Stephens' constant $0.5759...$ and Matthews' constant $0.1473...$ actually
first appeared in \cite{Stphns-err}. Let $\iota(n)=1$ if $n$ is square-free
and $\iota(n)=0$ otherwise. Then \cite{Carl-err, HtBr1-err, Inghm-err,
Mrsky1-err, Mrsky2-err, Rrick-err, Pappa-err}
\begin{align*}
\lim_{N\rightarrow\infty}\frac{1}{N}%
%TCIMACRO{\dsum \limits_{n=1}^{N}}%
%BeginExpansion
{\displaystyle\sum\limits_{n=1}^{N}}
%EndExpansion
\iota(n)\iota(n+1)  &  =%
%TCIMACRO{\dprod \limits_{p}}%
%BeginExpansion
{\displaystyle\prod\limits_{p}}
%EndExpansion
\left(  1-\frac{2}{p^{2}}\right)  =0.3226340989...=-1+2(0.6613170494...)\\
\  &  =\frac{6}{\pi^{2}}%
%TCIMACRO{\dprod \limits_{p}}%
%BeginExpansion
{\displaystyle\prod\limits_{p}}
%EndExpansion
\left(  1-\frac{1}{p^{2}-1}\right)  =\frac{6}{\pi^{2}}(0.5307118205...),
\end{align*}
that is, the Feller-Tornier constant arises with regard to consecutive
square-free numbers and to other problems. Also, consider the cardinality
$N(X)$ of nontrivial primitive integer vectors $(x_{0},x_{1},x_{2},x_{3})$
that fall on Cayley's cubic surface
\[
x_{0}x_{1}x_{2}+x_{0}x_{1}x_{3}+x_{0}x_{2}x_{3}+x_{1}x_{2}x_{3}=0
\]
and satisfy $|x_{j}|\leq X$ for $0\leq j\leq3$. It is known that $N(X)\sim
cX(\ln(X))^{6}$ for some constant $c>0$ \cite{HtBr2-err, SwDy-err}; finding
$c$ remains an open problem.

\textbf{2.5. Hafner-Sarnak-McCurley Constant.} In the \textquotedblleft Added
In Press\textquotedblright\ section (pages 601--602), the asymptotics of
coprimality and of square-freeness are discussed for the Gaussian integers and
for the Eisenstein-Jacobi integers. Generalizations appear in \cite{DSncR-err,
BrdyCL-err}. Cai \&\ Bach \cite{CaiBch-err} and T\'{o}th \cite{LToth-err}
independently proved that the probability that $k$ positive integers are
\textit{pairwise} coprime is \cite{TnbmWu-err, DIToev-err}
\[%
%TCIMACRO{\dprod \limits_{p}}%
%BeginExpansion
{\displaystyle\prod\limits_{p}}
%EndExpansion
\left(  1-\frac{1}{p}\right)  ^{k-1}\left(  1+\frac{k-1}{p}\right)
=\lim_{N\rightarrow\infty}\frac{(k-1)!}{N\ln(N)^{k-1}}%
%TCIMACRO{\dsum \limits_{n=1}^{N}}%
%BeginExpansion
{\displaystyle\sum\limits_{n=1}^{N}}
%EndExpansion
k^{\omega(n)}.
\]
Freiberg \cite{Frbg-err, Hymn-err, JHu-err}, building on Moree's work
\cite{PMoree3-err}, determined the probability that three positive integers
are pairwise \textit{not} coprime to be $1-18/\pi^{2}+3P-Q=0.1742197830....$
The constant $Q$ also appears in \cite{FMSbh-err, CDFHP-err, HuaVg-err}. More
about sums involving $2^{\omega(n)}$ and $2^{-\omega(n)}$ appears in
\cite{Fi2-err}. The asymptotics of $%
%TCIMACRO{\tsum \nolimits_{n=1}^{N}}%
%BeginExpansion
{\textstyle\sum\nolimits_{n=1}^{N}}
%EndExpansion
3^{\Omega(n)}$, due to Tenenbaum, are mentioned in \cite{Fi1-err}. Also, we
have \cite{EChn-err}
\[%
%TCIMACRO{\dsum \limits_{n\leq N}}%
%BeginExpansion
{\displaystyle\sum\limits_{n\leq N}}
%EndExpansion
\kappa(n)^{\ell}\sim\frac{1}{\ell+1}\frac{\zeta(2\ell+2)}{\zeta(2)}N^{\ell
+1},
\]%
\[%
%TCIMACRO{\dsum \limits_{n\leq N}}%
%BeginExpansion
{\displaystyle\sum\limits_{n\leq N}}
%EndExpansion
K(n)^{\ell}\sim\frac{1}{\ell+1}\frac{\zeta(\ell+1)}{\zeta(2)}%
%TCIMACRO{\dprod \limits_{p}}%
%BeginExpansion
{\displaystyle\prod\limits_{p}}
%EndExpansion
\left(  1-\frac{1}{p^{\ell}(p+1)}\right)  \cdot N^{\ell+1}%
\]
as $N\rightarrow\infty$, for any positive integer $\ell$. In the latter
formula, the product for $\ell=1$ and $\ell=2$ appears in \cite{Fi2-err} with
regard to the number/sum of unitary square-free divisors; the product for
$\ell=2$ further is connected with class number theory \cite{Fi15-err}.

\textbf{2.6. Niven's Constant.} The quantity $C$ appears unexpectedly in
\cite{Pmrnc-err}. If we instead examine the mean of the exponents:
\[
L(m)=\left\{
\begin{array}
[c]{lll}%
1 &  & \text{if }m=1,\\
\dfrac{1}{k}%
%TCIMACRO{\dsum \limits_{j=1}^{k}}%
%BeginExpansion
{\displaystyle\sum\limits_{j=1}^{k}}
%EndExpansion
a_{j} &  & \text{if }m>1,
\end{array}
\right.
\]
then \cite{Duncan-err, HZCao-err}
\[%
%TCIMACRO{\dsum \limits_{m\leq n}}%
%BeginExpansion
{\displaystyle\sum\limits_{m\leq n}}
%EndExpansion
L(m)=n+C_{1}\frac{n}{\ln(\ln(n))}+C_{2}\frac{n}{\ln(\ln(n))^{2}}+O\left(
\frac{n}{\ln(\ln(n))^{3}}\right)
\]
as $n\rightarrow\infty$, where \cite{HChn-err}
\[
C_{1}=%
%TCIMACRO{\dsum \limits_{p}}%
%BeginExpansion
{\displaystyle\sum\limits_{p}}
%EndExpansion
\frac{1}{p(p-1)}=M^{\prime}-M=0.7731566690...,
\]%
\[
C_{2}=%
%TCIMACRO{\dsum \limits_{p}}%
%BeginExpansion
{\displaystyle\sum\limits_{p}}
%EndExpansion
\frac{1}{p^{2}(p-1)}-C_{1}M=C_{1}(1-M)-N=0.1187309349...,
\]
using notation defined on pages 94--95. The constant $C_{1}$ also appears in
our earlier entry [2.2]. A general formula for coefficients $c_{ij}$ was found
by Sinha \cite{Sinha-err} and gives two additional terms (involving $n^{1/6}$
and $n^{1/7}$) in the asymptotic estimate of $%
%TCIMACRO{\tsum \nolimits_{m=1}^{n}}%
%BeginExpansion
{\textstyle\sum\nolimits_{m=1}^{n}}
%EndExpansion
h(m)$.

Let $\tilde{N}_{2}(x)$ denote the number of positive integer primitive triples
$(i,j,k)$ with $i+j=k\leq x$ and $i,j,k$ square-full. It is conjectured that
\cite{BrwngV-err}%
\[
\tilde{N}_{2}(x)=\tilde{c}\,x^{1/2}\left(  1+o(1)\right)
\]
as $x\rightarrow\infty$, where $\tilde{c}=2.677539267...$ has a complicated
expression. \ Supporting evidence includes the inequality $\tilde{N}%
_{2}(x)\geq\tilde{c}\,x^{1/2}\left(  1+o(1)\right)  $ and $\tilde{N}%
_{2}(x)=O\left(  x^{3/5}\ln(x)^{12}\right)  $.

\textbf{2.7. Euler Totient Constants.} Let us clarify the third sentence:
$\varphi(n)$ is the number of generators in $\mathbb{Z}_{n}$, the additive
group of integers modulo $n$. It is also the number of elements in
$\mathbb{Z}_{n}^{\ast}$, the multiplicative group of invertible integers
modulo $n$.

The mean of $\varphi(n)/n$ is well known; a conjectured formula for the
variance awaits proof \cite{Ovrbk-err}.

Define $f(n)=n\varphi(n)^{-1}-e^{\gamma}\ln(\ln(n))$. Nicolas \cite{Nic2-err}
proved that $f(n)>0$ for infinitely many integers $n$ by the following
reasoning. Let $P_{k}$ denote the product of the first $k$ prime numbers. If
the Riemann hypothesis is true, then $f(P_{k})>0$ for all $k$. If the Riemann
hypothesis is false, then $f(P_{k})>0$ for infinitely many $k$ and
$f(P_{l})\leq0$ for infinitely many $l$.

Let $U(n)$ denote the set of values $\leq n$ taken by $\varphi$ and $v(n)$
denote its cardinality; for example \cite{Sloa-err},
$U(15)=\{1,2,4,6,8,10,12\}$ and $v(15)=7$. Let $\ln_{2}(x)=\ln(\ln(x))$ and
$\ln_{m}(x)=$ $\ln(\ln_{m-1}(x))$ for convenience. Ford \cite{Fo1-err} proved
that
\[
v(n)=\tfrac{n}{\ln(n)}\exp\left\{  C[\ln_{3}(n)-\ln_{4}(n)]^{2}+D\ln
_{3}(n)-[D+\tfrac{1}{2}-2C]\ln_{4}(n)+O(1)\right\}
\]
as $n\rightarrow\infty$, where
\[
C=-\tfrac{1}{2\ln(\rho)}=0.8178146464...,
\]%
\[
D=2C\left(  1+\ln(F^{\prime}(\rho))-\ln(2C)\right)  -\tfrac{3}{2}%
=2.1769687435...
\]%
\[
F(x)=%
%TCIMACRO{\dsum \limits_{k=1}^{\infty}}%
%BeginExpansion
{\displaystyle\sum\limits_{k=1}^{\infty}}
%EndExpansion
\left(  (k+1)\ln(k+1)-k\ln(k)-1\right)  x^{k}%
\]
and $\rho=0.5425985860...$ is the unique solution on $[0,1)$ of the equation
$F(\rho)=1$. Also,
\[
\lim_{n\rightarrow\infty}\frac{1}{v(n)\ln_{2}(n)}%
%TCIMACRO{\dsum \limits_{m\in U(n)}}%
%BeginExpansion
{\displaystyle\sum\limits_{m\in U(n)}}
%EndExpansion
\omega(m)=\frac{1}{1-\rho}=2.1862634648...
\]
which contrasts with a related result of Erd\H{o}s \&\ Pomerance
\cite{ErPo-err}:
\[
\lim_{n\rightarrow\infty}\frac{1}{n\ln_{2}(n)^{2}}%
%TCIMACRO{\dsum \limits_{m=1}^{n}}%
%BeginExpansion
{\displaystyle\sum\limits_{m=1}^{n}}
%EndExpansion
\omega(\varphi(n))=\frac{1}{2}.
\]
These two latter formulas hold as well if $\omega$ is replaced by $\Omega$.
See \cite{Fo2-err} for more on Euler's totient.

Define the \textit{reduced totient} or \textit{Carmichael function} $\psi(n) $
to be the size of the largest cyclic subgroup of $\mathbb{Z}_{n}^{*}$. We have
\cite{EPS-err}
\[
\frac1N%
%TCIMACRO{\dsum \limits_{n\leq N}}%
%BeginExpansion
{\displaystyle\sum\limits_{n\leq N}}
%EndExpansion
\psi(n)=\frac N{\ln(N)}\exp\left[  \frac{P\ln_{2}(N)}{\ln_{3}(N)}\left(
1+o(1)\right)  \right]
\]
as $N\rightarrow\infty$, where
\[
P=e^{-\gamma}%
%TCIMACRO{\dprod \limits_{p}}%
%BeginExpansion
{\displaystyle\prod\limits_{p}}
%EndExpansion
\left(  1-\frac1{(p-1)^{2}(p+1)}\right)  =0.3453720641....
\]
(note the similarity to a constant in \cite{Fi21-err}.) There is a set $S$ of
positive integers of asymptotic density $1$ such that, for $n\in S$,
\[
n\psi(n)^{-1}=\left(  \ln(n)\right)  ^{\ln_{3}(n)+Q+o(1)}%
\]
and
\[
Q=-1+%
%TCIMACRO{\dsum \limits_{p}}%
%BeginExpansion
{\displaystyle\sum\limits_{p}}
%EndExpansion
\dfrac{\ln(p)}{(p-1)^{2}}=0.2269688056...;
\]
it is not known whether $S=\mathbb{Z}^{+}$ is possible.

Let $X_{n}$ denote the $\gcd$ of two integers chosen independently from
Uniform $\{1,2,\ldots,n\}$ and $Y_{n}$ denote the $\operatorname*{lcm}$.
Diaconis \&\ Erd\H{o}s \cite{DiacE-err} proved that
\[%
\begin{array}
[c]{ccc}%
\operatorname*{E}(X_{n})=\dfrac{6}{\pi^{2}}\ln(n)+E+O\left(  \dfrac{\ln
(n)}{\sqrt{n}}\right)  , &  & \operatorname*{E}(Y_{n})=\dfrac{3\zeta(3)}%
{2\pi^{2}}n^{2}+O\left(  n\ln(n)\right)
\end{array}
\]
as $n\rightarrow\infty$, where
\[
E=%
%TCIMACRO{\dsum \limits_{k=1}^{\infty}}%
%BeginExpansion
{\displaystyle\sum\limits_{k=1}^{\infty}}
%EndExpansion
\tfrac{1}{k^{2}(k+1)^{2}}\left\{
%TCIMACRO{\tsum \limits_{j=1}^{k}}%
%BeginExpansion
{\textstyle\sum\limits_{j=1}^{k}}
%EndExpansion
\varphi(j)+2\left(  -\tfrac{3}{\pi^{2}}k^{2}+%
%TCIMACRO{\tsum \limits_{j=1}^{k}}%
%BeginExpansion
{\textstyle\sum\limits_{j=1}^{k}}
%EndExpansion
\varphi(j)\right)  k-\tfrac{6}{\pi^{2}}(2k+1)k\right\}  +\tfrac{12}{\pi^{2}%
}\left(  \gamma+\tfrac{1}{2}\right)  -\tfrac{1}{2}%
\]
but a vastly simpler expression
\[
E=\dfrac{6}{\pi^{2}}\left(  2\gamma-\frac{1}{2}-\frac{\pi^{2}}{12}-\dfrac
{6}{\pi^{2}}\zeta^{\prime}(2)\right)
\]
was found earlier by Cohen \cite{CohE-err, FrnE-err}; a reconcilation is needed.

\textbf{2.8. Pell-Stevenhagen Constants.} The constant $P$ is transcendental
via a general theorem on values of modular forms due to Nesterenko
\cite{Nest-err, Zud2-err}. Here is a constant similar to $P$: The number of
positive integers $n\leq N$, for which $2n-1$ is not divisible by $2^{p}-1$
for any prime $p$, is $\sim cN$, where
\[
c=%
%TCIMACRO{\dprod \limits_{p}}%
%BeginExpansion
{\displaystyle\prod\limits_{p}}
%EndExpansion
\left(  1-\frac1{2^{p}-1}\right)  =0.5483008312....
\]
A ring-theoretic analog of this statement, plus generalizations, appear in
\cite{BoSc-err}.

\textbf{2.9. Alladi-Grinstead Constant.} In the final paragraph, it should be
noted that the first product $1.7587436279...$ is $e^{C}/2$. See
\cite{Fi19-err} for another occurrence of $C$. It is a multiplicative analog
of Euler's constant $\gamma$ in the sense that \cite{Slzdrf-err}%
\[%
\begin{array}
[c]{ccc}%
\gamma=%
%TCIMACRO{\dint \limits_{1}^{\infty}}%
%BeginExpansion
{\displaystyle\int\limits_{1}^{\infty}}
%EndExpansion
\left(  \dfrac{1}{\left\lfloor x\right\rfloor }-\dfrac{1}{x}\right)  dx, &  &
C=%
%TCIMACRO{\dint \limits_{1}^{\infty}}%
%BeginExpansion
{\displaystyle\int\limits_{1}^{\infty}}
%EndExpansion
\left(  \dfrac{1}{\left\lfloor x\right\rfloor }\dfrac{1}{x}\right)  dx.
\end{array}
\]
It also has an alternating analog \cite{Boyz1-err, Boyz2-err}%
\[
\tilde{C}=%
%TCIMACRO{\dsum \limits_{n=1}^{\infty}}%
%BeginExpansion
{\displaystyle\sum\limits_{n=1}^{\infty}}
%EndExpansion
\frac{(-1)^{n-1}}{n}\ln\left(  1+\frac{1}{n}\right)  =0.552128322...
\]
with parallel formulas%
\[%
\begin{array}
[c]{ccc}%
%TCIMACRO{\dsum \limits_{k=1}^{\infty}}%
%BeginExpansion
{\displaystyle\sum\limits_{k=1}^{\infty}}
%EndExpansion
\dfrac{(-1)^{k-1}}{k}\zeta(k+1)=C, &  &
%TCIMACRO{\dsum \limits_{k=1}^{\infty}}%
%BeginExpansion
{\displaystyle\sum\limits_{k=1}^{\infty}}
%EndExpansion
\dfrac{(-1)^{k-1}}{k}\eta(k+1)=\tilde{C}%
\end{array}
\]
where $\eta(x)=(1-2^{1-x})\zeta(x)$ for $x\neq1$. \ Yet another constant%
\[%
%TCIMACRO{\dsum \limits_{n=1}^{\infty}}%
%BeginExpansion
{\displaystyle\sum\limits_{n=1}^{\infty}}
%EndExpansion
\frac{1}{n}\ln\left(  1+\frac{1}{n+1}\right)  =0.8606201928...
\]
awaits further study.

\textbf{2.10. Sierpinski's Constant.} Sierpinski's formulas for $\hat{S}$ and
$\tilde{S}$ contained a few errors:\ they should be \cite{Rama-err, KuNo-err,
BoCh-err, Stron-err, Frick-err, Fisch-err}
\[
\hat{S}=\gamma+S-\dfrac{12}{\pi^{2}}\zeta^{\prime}(2)+\dfrac{\ln(2)}%
{3}-1=1.7710119609...=\frac{\pi}{4}(2.2549224628...),
\]%
\[
\tilde{S}=2S-\dfrac{12}{\pi^{2}}\zeta^{\prime}(2)+\dfrac{\ln(2)}%
{3}-1=2.0166215457...=\frac{1}{4}(8.0664861829...).
\]
In the summation formula at the top of page 125, $D_{n}$ should be $D_{k}$.
Also, the divisor analog of Sierpinski's second series is \cite{Jager-err}
\[%
%TCIMACRO{\dsum \limits_{k=1}^{n}}%
%BeginExpansion
{\displaystyle\sum\limits_{k=1}^{n}}
%EndExpansion
d(k^{2})=\left(  \frac{3}{\pi^{2}}\ln(n)^{2}+\left(  \frac{18\gamma-6}{\pi
^{2}}-\frac{72}{\pi^{4}}\zeta^{\prime}(2)\right)  \ln(n)+c\right)  n+O\left(
n^{1/2+\varepsilon}\right)
\]
as $n\rightarrow\infty$, where the expression for $c$ is complicated. It is
easily shown that $d(n^{2})$ is the number of ordered pairs of positive
integers $(i,j)$ satisfying $\operatorname*{lcm}(i,j)=n$. \ 

The best known result for $r(n)$ is currently \cite{Hxly03-err}
\[
\sum_{k=1}^{n}r(k)=\pi\,n+O\left(  n^{\frac{131}{416}}\ln(n)^{\frac
{18627}{8320}}\right)  .
\]
Define $R(n)$ to be the number of representations of $n$ as a sum of three
squares, counting order and sign. Then
\[%
%TCIMACRO{\dsum \limits_{k=1}^{n}}%
%BeginExpansion
{\displaystyle\sum\limits_{k=1}^{n}}
%EndExpansion
R(k)=\frac{4\pi}3n^{3/2}+O\left(  n^{3/4+\varepsilon}\right)
\]
for all $\varepsilon>0$ and \cite{CKuO-err}
\[%
%TCIMACRO{\dsum \limits_{k=1}^{n}}%
%BeginExpansion
{\displaystyle\sum\limits_{k=1}^{n}}
%EndExpansion
R(k)^{2}=\frac{8\pi^{4}}{21\zeta(3)}n^{2}+O\left(  n^{14/9}\right)  .
\]
The former is the same as the number of integer ordered triples falling within
the ball of radius $\sqrt{n}$ centered at the origin; an extension of the
latter to sums of $m$ squares, when $m>3$, is also known \cite{CKuO-err}.

A claimed proof that
\[%
%TCIMACRO{\dsum \limits_{n\leq x}}%
%BeginExpansion
{\displaystyle\sum\limits_{n\leq x}}
%EndExpansion
d(n)=x\ln(x)+(2\gamma-1)x+O\left(  x^{1/4+\varepsilon}\right)
\]
as $x\rightarrow\infty$ awaits confirmation \cite{GMo-err}. Let $\delta(n) $
denote the number of square divisors of $n$, that is, all positive integers
$d$ for which $d^{2}|n$. It is known that \cite{ZhCao-err}
\[%
%TCIMACRO{\dsum \limits_{n\leq x}}%
%BeginExpansion
{\displaystyle\sum\limits_{n\leq x}}
%EndExpansion
\delta(n)\sim\zeta(2)x+\zeta(1/2)x^{1/2}%
\]
as $x\rightarrow\infty$. Analogous to various error-term formulas in
\cite{Fi17-err}, we have
\[%
%TCIMACRO{\dint \limits_{1}^{x}}%
%BeginExpansion
{\displaystyle\int\limits_{1}^{x}}
%EndExpansion
\left(
%TCIMACRO{\dsum \limits_{m\leq y}}%
%BeginExpansion
{\displaystyle\sum\limits_{m\leq y}}
%EndExpansion
\delta(m)-\zeta(2)y-\zeta(1/2)y^{1/2}\right)  ^{2}dy\sim C_{\delta}\,x^{4/3}%
\]
where
\[
C_{\delta}=\frac{2^{1/3}}{8\,\pi^{2}}%
%TCIMACRO{\dsum \limits_{n=1}^{\infty}}%
%BeginExpansion
{\displaystyle\sum\limits_{n=1}^{\infty}}
%EndExpansion
\left(
%TCIMACRO{\dsum \limits_{d^{2}|n}}%
%BeginExpansion
{\displaystyle\sum\limits_{d^{2}|n}}
%EndExpansion
\frac{d}{n^{5/6}}\right)  ^{2}.
\]
This supports a conjecture that the error in approximating $%
%TCIMACRO{\tsum \nolimits_{n\leq x}}%
%BeginExpansion
{\textstyle\sum\nolimits_{n\leq x}}
%EndExpansion
\delta(n)$ is $O(x^{1/6+\varepsilon})$. \ See also \cite{LcTh-err}.

\textbf{2.11. Abundant Numbers Density Constant.} An odd perfect number cannot
be less than $10^{1500}$ \cite{OchmR-err}. The definition of $A(x)$ should be
replaced by
\[
A(x)=\lim_{n\rightarrow\infty}\frac{\left\vert \left\{  k\leq n:\sigma(k)\geq
x\,k\right\}  \right\vert }{n}.
\]
Kobayashi \cite{abund1-err} proved that $0.24761<A(2)<0.24765$; see also
\cite{abund2-err, abund3-err, abund4-err, abund5-err, abund6-err}. \ If $K(x)$
is the number of all positive integers $m$ that satisfy $\sigma(m)\leq x$,
then \cite{Drsslr-err}
\begin{align*}
\lim_{x\rightarrow\infty}\frac{K(x)}{x}  &  =%
%TCIMACRO{\dprod \limits_{p}}%
%BeginExpansion
{\displaystyle\prod\limits_{p}}
%EndExpansion
\left(  1-\frac{1}{p}\right)  \left(  1+%
%TCIMACRO{\dsum \limits_{j=1}^{\infty}}%
%BeginExpansion
{\displaystyle\sum\limits_{j=1}^{\infty}}
%EndExpansion
\left(  1+%
%TCIMACRO{\dsum \limits_{i=1}^{j}}%
%BeginExpansion
{\displaystyle\sum\limits_{i=1}^{j}}
%EndExpansion
p^{i}\right)  ^{-1}\right) \\
&  =%
%TCIMACRO{\dprod \limits_{p}}%
%BeginExpansion
{\displaystyle\prod\limits_{p}}
%EndExpansion
\left(  1-\frac{1}{p}\right)  \left(  1+(p-1)%
%TCIMACRO{\dsum \limits_{j=1}^{\infty}}%
%BeginExpansion
{\displaystyle\sum\limits_{j=1}^{\infty}}
%EndExpansion
\frac{1}{p^{j+1}-1}\right)  .
\end{align*}

\textbf{2.12. Linnik's Constant.} In the definition of $L$, \textquotedblleft
lim\textquotedblright\ should be replaced by \textquotedblleft
limsup\textquotedblright. Clearly $L$ exists; the fact that $L<\infty$ was
Linnik's important contribution. Xylouris \cite{TXyl1-err} recently proved
that $L\leq5.18$; an unpublished\ proof that $L\leq5$ needs to be verified
\cite{TXyl2-err}.

\textbf{2.13. Mills' Constant.} Caldwell \&\ Cheng \cite{CwCg-err} computed
$C$ to high precision; Baillie \cite{Bailli-err} likewise examined $c$. \ The
question, \textquotedblleft Does there exist $\tilde{C}>1$ for which
$\left\lfloor \tilde{C}^{n}\right\rfloor $ is always prime?\textquotedblright,
remains open \cite{Farhi-err}. \ Replacing floor by ceiling, T\'{o}th
\cite{TothML-err} found the analog of $C$ to be $1.2405547052...$. \ Let
$q_{1}<q_{2}<\ldots<q_{k}$ denote the consecutive prime factors of an integer
$n>1$. Define
\[
F(n)=%
%TCIMACRO{\dsum \limits_{j=1}^{k-1}}%
%BeginExpansion
{\displaystyle\sum\limits_{j=1}^{k-1}}
%EndExpansion
\left(  1-\frac{q_{j}}{q_{j+1}}\right)  =\omega(n)-1-%
%TCIMACRO{\dsum \limits_{j=1}^{k-1}}%
%BeginExpansion
{\displaystyle\sum\limits_{j=1}^{k-1}}
%EndExpansion
\frac{q_{j}}{q_{j+1}}%
\]
if $k>1$ and $F(n)=0$ if $k=1$. Erd\H{o}s \&\ Nicolas \cite{ErNi-err}
demonstrated that there exists a constant $C^{\prime}=1.70654185...$ such
that, as $n\rightarrow\infty$, $F(n)\leq\sqrt{\ln(n)}-C^{\prime}+o(1)$, with
equality holding for infinitely many $n$. Further, $C^{\prime}=C^{\prime
\prime}+\ln(2)+1/2$, where \cite{ErNi-err, Sebah-err}
\[%
\begin{array}
[c]{ccc}%
C^{\prime\prime}=%
%TCIMACRO{\dsum \limits_{i=1}^{\infty}}%
%BeginExpansion
{\displaystyle\sum\limits_{i=1}^{\infty}}
%EndExpansion
\left\{  \ln\left(  \dfrac{p_{i+1}}{p_{i}}\right)  -\left(  1-\dfrac{p_{i}%
}{p_{i+1}}\right)  \right\}  =0.51339467..., &  &
%TCIMACRO{\dsum \limits_{i=1}^{\infty}}%
%BeginExpansion
{\displaystyle\sum\limits_{i=1}^{\infty}}
%EndExpansion
\left(  \dfrac{p_{i+1}}{p_{i}}-1\right)  ^{2}=1.65310351...,
\end{array}
\]
and $p_{1}=2$, $p_{2}=3$, $p_{3}=5$, ... is the sequence of all primes.

It now seems that $\operatorname*{liminf}_{n\rightarrow\infty}(p_{n+1}%
-p_{n})/\ln(p_{n})=0$ is a theorem \cite{GMPY-err, GGPY-err}, clarifying the
uncertainty raised in ``Added In Press'' (pages 601--602). More about small
prime gaps will surely appear soon; research concerning large prime gaps
continues as well \cite{MaiPm-err, Pintz-err}.

\textbf{2.14. Brun's Constant.} Klyve \cite{Klyv1-err} conjectured that
$1.90216054<B_{2}<1.90216063$ (consistent with Sebah's estimate); such bounds,
although very carefully established, are speculative. \ The best unconditional
interval is $1.8304<B_{2}<2.3467$. \ Assuming the truth of an extended Riemann
hypothesis \cite{Klyv1-err}, the upper limit can be reduced to $2.1754$.
\ Wolf \cite{MWolf-err} nonrigorously calculated that $\tilde{B}%
_{4}=1.1970449...$ and a meticulous check remains open.

\textbf{2.15. Glaisher-Kinkelin Constant.} A certain infinite product
\cite{FurdTylr-err}%
\[%
%TCIMACRO{\dprod \limits_{n=1}^{\infty}}%
%BeginExpansion
{\displaystyle\prod\limits_{n=1}^{\infty}}
%EndExpansion
\left(  \frac{n!}{\sqrt{2\pi n}(n/e)^{n}}\right)  ^{(-1)^{n-1}}=\frac{A^{3}%
}{2^{7/12}\pi^{1/4}}%
\]
features the ratio of $n!$ to its Stirling approximation. In the second
display for $D(x)$, $\exp(-x/2)$ should be replaced by $\exp(x/2)$. Another
proof of the formula for $D(1)$ is given in \cite{So0-err}; another special
case is \cite{KchTzr1-err, KchTzr2-err, BrwDyk-err}
\[
D(1/2)=\frac{2^{1/6}\sqrt{\pi}A^{3}}{\Gamma(1/4)}e^{G/\pi}.
\]
The two quantities
\[%
\begin{array}
[c]{ccc}%
G_{2}\left(  \frac{1}{2}\right)  =0.6032442812..., &  & G_{2}\left(  \frac
{3}{2}\right)  =\sqrt{\pi}G_{2}\left(  \frac{1}{2}\right)  =1.0692226492...
\end{array}
\]
play a role in a discussion of the limiting behavior of Toeplitz determinants
and the Fisher-Hartwig conjecture \cite{Basr1-err, Basr2-err}. Krasovsky
\cite{Krsv-err} and Ehrhardt \cite{Ehrd-err} proved Dyson's conjecture
regarding the asymptotic expansion of $E(s)$ as $s\rightarrow\infty$; a third
proof is given in \cite{DIKZ-err}. Also, the quantities
\[
G_{2}\left(  \tfrac{1}{2}\right)  ^{-1}=1.6577032408...=2^{-1/24}e^{-3/16}%
\pi^{1/4}(3.1953114860...)^{3/8}%
\]%
\[
G_{3}\left(  \tfrac{3}{2}\right)  ^{-1}=G_{2}\left(  \tfrac{1}{2}\right)
G_{3}\left(  \tfrac{1}{2}\right)  ^{-1}=0.9560900097...=\pi^{-1/2}%
(3.3388512141...)^{7/16}%
\]
appear in \cite{DkIm-err}. In the last paragraph on page 141, the polynomial
$q(x)$ should be assumed to have degree $n$. See \cite{FZ-err, Frmr-err} for
more on the GUE\ hypothesis.

Here is a sample result involving not random real polynomials, but a random
complex exponential. Let $a$, $b$ denote independent complex Gaussian
coefficients. The expected number of zeroes of $a+b\exp(z)$ satisfying $|z|<1$
is \cite{Mjvch-err}%
\[
\frac{1}{\pi}%
%TCIMACRO{\diint \limits_{x^{2}+y^{2}<1}}%
%BeginExpansion
{\displaystyle\iint\limits_{x^{2}+y^{2}<1}}
%EndExpansion
\frac{\exp(2x)}{\left(  1+\exp(2x)\right)  ^{2}}dx\,dy=0.2029189212...
\]
and higher-degree results are also known.

\textbf{2.16. Stolarsky-Harboth Constant.} The \textquotedblleft typical
growth\textquotedblright\ of $2^{b(n)}$ is $\approx n^{1/2}$ while the
\textquotedblleft average growth\textquotedblright\ of $2^{b(n)}$ is $\approx
n^{\ln(3/2)/\ln(2)}$; more examples are found in \cite{FSB1-err}. The
\textquotedblleft typical dispersion\textquotedblright\ of $2^{b(n)}$ is
$\approx n^{\ln(2)/4}$ while the \textquotedblleft average
dispersion\textquotedblright\ of $2^{b(n)}$ is $\approx n^{\ln(5/2)/\ln(2)} $;
more examples are found in \cite{FSB2-err}. A generating function%
\[%
%TCIMACRO{\dsum \limits_{n=1}^{\infty}}%
%BeginExpansion
{\displaystyle\sum\limits_{n=1}^{\infty}}
%EndExpansion
b(n)z^{n}=\frac{1}{1-z}%
%TCIMACRO{\dsum \limits_{\ell=0}^{\infty}}%
%BeginExpansion
{\displaystyle\sum\limits_{\ell=0}^{\infty}}
%EndExpansion
\frac{z^{2^{\ell}}}{z^{2^{\ell}}+1}%
\]
is interesting since it involves a Lambert series \cite{AWRsky-err,
VgntWk-err}. \ Coquet's 1983 result is discussed in \cite{Shevlv-err} and a
misprint is corrected. \ Ulam's 1-additive sequence $(1,2)$ surprisingly may
possess some rigid underlying structure \cite{Ulam-err}. \ The sequence
$\{0\}\cup\{c(n)\}_{n=0}^{\infty}$ is called Stern's diatomic sequence
\cite{Nrthshld-err} and our final question is answered in \cite{CnsTylr-err}:%
\[
\operatorname*{limsup}\limits_{n\rightarrow\infty}\frac{c(n)}{n^{\ln
(\varphi)/\ln(2)}}=\frac{\varphi}{\sqrt{5}}\left(  \frac{3}{2}\right)
^{\ln(\varphi)/\ln(2)}=\frac{\varphi^{\ln(3)/\ln(2)}}{\sqrt{5}}%
=0.9588541908....
\]

Given a positive integer $n$, define $s_{1}^{2}$ to be the largest square not
exceeding $n$. Then define $s_{2}^{2}$ to be the largest square not exceeding
$n-s_{1}^{2}$, and so forth. Hence $n=%
%TCIMACRO{\tsum \nolimits_{j=1}^{r}}%
%BeginExpansion
{\textstyle\sum\nolimits_{j=1}^{r}}
%EndExpansion
s_{j}^{2}$ for some $r$. We say that $n$ is a \textit{greedy sum of distinct
squares} if $s_{1}>s_{2}>\ldots>s_{r}$. Let $A(N)$ be the number of such
integers $n<N$, plus one. Montgomery \&\ Vorhauer \cite{MoVo-err} proved that
$A(N)/N$ does not tend to a constant, but instead that there is a continuous
function $f(x)$ of period $1$ for which
\[%
\begin{array}
[c]{ccc}%
\lim\limits_{k\rightarrow\infty}\dfrac{A(4\exp(2^{k+x}))}{4\exp(2^{k+x}%
)}=f(x), &  & \min\limits_{0\leq x\leq1}f(x)=0.50307...<\max\limits_{0\leq
x\leq1}f(x)=0.50964...
\end{array}
\]
where $k$ takes on only integer values. This is reminiscent of the behavior
discussed for digital sums.

Two simple examples, due to Hardy \cite{Hrdy1-err, Hrdy2-err} and Elkies
\cite{Elki-err}, involve the series
\[%
\begin{array}
[c]{ccc}%
\varphi(x)=%
%TCIMACRO{\dsum \limits_{k=0}^{\infty}}%
%BeginExpansion
{\displaystyle\sum\limits_{k=0}^{\infty}}
%EndExpansion
x^{2^{k}}, &  & \psi(x)=%
%TCIMACRO{\dsum \limits_{k=0}^{\infty}}%
%BeginExpansion
{\displaystyle\sum\limits_{k=0}^{\infty}}
%EndExpansion
(-1)^{k}x^{2^{k}}.
\end{array}
\]
As $x\rightarrow1^{-}$, the asymptotics of $\varphi(x)$ and $\psi(x)$ are
complicated by oscillating errors with amplitude
\[
\sup_{x\rightarrow1^{-}}\left|  \varphi(x)+\frac{\ln(-\ln(x))+\gamma}{\ln
(2)}-\frac32+x\right|  =(1.57...)\times10^{-6},
\]
\[
\sup_{x\rightarrow1^{-}}\left|  \psi(x)-\frac16-\frac13x\right|
=(2.75...)\times10^{-3}.
\]
The function $\varphi(x)$ also appears in what is known as Catalan's integral
(section 1.5.2) for Euler's constant $\gamma$. See \cite{GrbHwg-err,
AllaartK-err} as well.

\textbf{2.17. Gauss-Kuzmin-Wirsing Constant.} If $X$ is a random variable
following the Gauss-Kuzmin distribution, then its mean value is
\begin{align*}
\operatorname*{E}(X)  &  =\frac{1}{\ln(2)}%
%TCIMACRO{\dint \limits_{0}^{1}}%
%BeginExpansion
{\displaystyle\int\limits_{0}^{1}}
%EndExpansion
\frac{x}{1+x}dx=\frac{1}{\ln(2)}-1=0.4426950408...\\
&  =\frac{1}{\ln(2)}%
%TCIMACRO{\dint \limits_{0}^{1}}%
%BeginExpansion
{\displaystyle\int\limits_{0}^{1}}
%EndExpansion
\frac{\{1/x\}}{1+x}dx=\operatorname*{E}\left\{  \frac{1}{X}\right\}  .
\end{align*}
Further,
\begin{align*}
\operatorname*{E}\left(  \log_{10}\left(  X\right)  \right)   &  =\frac{1}%
{\ln(2)}%
%TCIMACRO{\dint \limits_{0}^{1}}%
%BeginExpansion
{\displaystyle\int\limits_{0}^{1}}
%EndExpansion
\frac{\log_{10}(x)}{1+x}dx=-\frac{\pi^{2}}{12\ln(2)\ln(10)}=-0.5153204170...\\
\  &  =\frac{1}{\ln(2)}%
%TCIMACRO{\dint \limits_{0}^{1}}%
%BeginExpansion
{\displaystyle\int\limits_{0}^{1}}
%EndExpansion
\frac{\log_{10}\{1/x\}}{1+x}dx=\operatorname*{E}\left(  \log_{10}\left\{
\frac{1}{X}\right\}  \right)  ,
\end{align*}
a constant that appears in \cite{Silva-err} and our earlier entry [1.8]. The
ratio conjecture involving eigenvalues of $G_{2}$ is now known to be true
\cite{Alksk-err}; moreover, the first two terms in the asymptotic series for
eigenvalues (involving $\varphi$ and $\zeta(3/2)$) are available. \ An attempt
to express $\lambda_{1}^{\prime\prime}(2)-\lambda_{1}^{\prime}(2)^{2}$ in
elementary terms appears in \cite{Fi19-err}.

A generalized Gauss transformation appears in \cite{PSun1-err, PSun2-err}.
\ The preprint math.NT/9908043 was withdrawn by the author without comment;
additional references on the Hausdorff dimension $0.5312805062...$ of real
numbers with partial denominators in $\{1,2\}$ include \cite{JnksPlc1-err,
OJnks-err, Lho1-err, Lho2-err, FlkNss-err, JnksPlc2-err}.

\textbf{2.18. Porter-Hensley Constants.} The formula for $H$ is wrong (by a
factor of $\pi^{6}$) and should be replaced by
\[
H=-\frac{\lambda_{1}^{\prime\prime}(2)-\lambda_{1}^{\prime}(2)^{2}}%
{\lambda_{1}^{\prime}(2)^{3}}=0.5160624088...=(0.7183748387...)^{2}.
\]
Lhote \cite{Lho1-err, Lho2-err} developed rigorous techniques for computing
$H$ and other constants to high precision. Ustinov \cite{Ustnv1-err,
Ustnv2-err} expressed Hensley's constant using some singular series:
\[
H=\frac{288\ln(2)^{2}}{\pi^{4}}\left(  \gamma-\frac{\zeta^{\prime}(2)}%
{\zeta(2)}-\frac{\ln(2)}2-1\right)  +\frac{24}{\pi^{2}}\left(  D+\frac
{3\ln(2)}2\right)
\]
where
\begin{align*}
D  &  =\ln\left(  \frac43\right)  -2\ln(2)^{2}+\\
&  \
%TCIMACRO{\dsum \limits_{n=2}^{\infty}}%
%BeginExpansion
{\displaystyle\sum\limits_{n=2}^{\infty}}
%EndExpansion
\left(
%TCIMACRO{\dsum \limits_{k,m=1}^{n}}%
%BeginExpansion
{\displaystyle\sum\limits_{k,m=1}^{n}}
%EndExpansion
\delta_{n}(k\,m+1)%
%TCIMACRO{\dint \limits_{0}^{1}}%
%BeginExpansion
{\displaystyle\int\limits_{0}^{1}}
%EndExpansion
\frac{d\xi}{(m\,\xi+n)\left[  \left(  \frac1n(k\,m+1)+m\right)  \xi
+(k+n)\right]  }\right.  +\\
&  \ \;\;\left.
%TCIMACRO{\dsum \limits_{k,m=1}^{n}}%
%BeginExpansion
{\displaystyle\sum\limits_{k,m=1}^{n}}
%EndExpansion
\delta_{n}(k\,m-1)%
%TCIMACRO{\dint \limits_{0}^{1}}%
%BeginExpansion
{\displaystyle\int\limits_{0}^{1}}
%EndExpansion
\frac{d\xi}{(m\,\xi+n)\left[  \left(  \frac1n(k\,m-1)+m\right)  \xi
+(k+n)\right]  }-2\ln(2)^{2}\frac{\varphi(n)}{n^{2}}\right)
\end{align*}
and $\delta_{n}(j)=1$ if $j\equiv0\operatorname*{mod}n$, $\delta_{n}(j)=0$ otherwise.

With regard to the binary GCD\ algorithm, Maze \cite{Maze-err} and Morris
\cite{Mrss-err} confirmed Brent's functional equation for a certain limiting
distribution \cite{Brent-err}
\[%
\begin{array}
[c]{ccc}%
g(x)=%
%TCIMACRO{\dsum \limits_{k\geq1}}%
%BeginExpansion
{\displaystyle\sum\limits_{k\geq1}}
%EndExpansion
2^{-k}\left(  g\left(  \dfrac{1}{1+2^{k}/x}\right)  -g\left(  \dfrac
{1}{1+2^{k}x}\right)  \right)  , &  & 0\leq x\leq1
\end{array}
\]
as well as important regularity properties including the formula
\[
2+\frac{1}{\ln(2)}%
%TCIMACRO{\dint \limits_{0}^{1}}%
%BeginExpansion
{\displaystyle\int\limits_{0}^{1}}
%EndExpansion
\frac{g(x)}{1-x}dx=\frac{2}{\kappa\ln(2)}=2.8329765709...=\frac{\pi
^{2}(0.3979226811...)}{2\ln(2)}.
\]

\textbf{2.19. Vall\'ee's Constant.} The $k^{\text{th}}$ \textit{circular}
continuant polynomial is the sum of monomials obtained from $x_{1}x_{2}\cdots
x_{k} $ by crossing out in all possible ways pairs of adjacent variables
$x_{j}x_{j+1}$, where $x_{k}x_{1}$ is now regarded as adjacent. For example
\cite{BBCrtn-err},
\[%
\begin{array}
[c]{ccccc}%
x_{1}x_{2}+2, &  & x_{1}x_{2}x_{3}+x_{1}+x_{3}+x_{2}, &  & x_{1}x_{2}%
x_{3}x_{4}+x_{1}x_{2}+x_{4}x_{1}+x_{3}x_{4}+x_{2}x_{3}+2
\end{array}
\]
are the cases for $k=2,3,4$.

\textbf{2.20. Erd\H{o}s' Reciprocal Sum Constants.} Recent work
\cite{ErdSA1-err, ErdSA2-err} gives $2.0654<S(A)<3.0752$; we have not yet
checked claims in \cite{ErdSA3-err, ErdSA4-err}. \ Improved bounds on the
reciprocal sums of Mian-Chowla and of Zhang were calculated in
\cite{Salvia-err}; the best lower estimate of $S(B_{2})$, however, still
appears to be $2.16086$ \cite{Lewis-err}. \ A sequence of positive integers
$b_{1}<b_{2}<\ldots<b_{m}$ is a $B_{h}$-sequence if all $h$-fold sums
$b_{i_{1}}+b_{i_{2}}+\cdots+b_{i_{h}}$, $i_{1}\leq i_{2}\leq\ldots\leq i_{h}$,
are distinct. Given $n$, choose a $B_{h}$-sequence $\{b_{i}\}$ so that
$b_{m}\leq n$ and $m $ is maximal; let $F_{h}(n)$ be this value of $m$. It is
known that $C_{h}=\operatorname*{limsup}\nolimits_{n\rightarrow\infty}%
n^{-1/h}F_{h}(n)$ is finite; we further have \cite{Lind-err, APLi-err,
Grah-err, Kol1-err, Cill-err, Gre1-err}
\[%
\begin{array}
[c]{ccccc}%
C_{2}=1, &  & 1\leq C_{3}\leq(7/2)^{1/3}, &  & 1\leq C_{4}\leq7^{1/4}.
\end{array}
\]
More generally, a sequence of positive integers $b_{1}<b_{2}<\ldots<b_{m}$ is
a $B_{h,g}$-sequence if, for every positive integer $k$, the equation
$x_{1}+x_{2}+\cdots+x_{h}=k$, $x_{1}\leq x_{2}\leq\ldots\leq x_{h}$, has at
most $g$ solutions with $x_{j}=b_{i_{j}}$ for all $j$. Defining $F_{h,g}(n)$
and $C_{h,g}$ analogously, we have \cite{Gre1-err, CRT-err, Pla1-err,
Pla2-err, Obry-err, GYu-err, CVnu-err, CRV-err}
\[%
\begin{array}
[c]{ccc}%
\dfrac{4}{\sqrt{7}}\leq C_{2,2}\leq\dfrac{\sqrt{21}}{2}, &  & 1.1509\leq
\lim\limits_{g\rightarrow\infty}\dfrac{C_{2,g}}{g^{1/2}}=\sqrt{\dfrac{2}{S}%
}\leq1.2525
\end{array}
\]
where the \textquotedblleft self-convolution constant\textquotedblright\ $S$
appears in \cite{Fi31-err} and satisfies $1.2748\leq S\leq1.5098$.

Here is a similar problem: for $k\geq1$, let $\nu_{2}(k)$ be the largest
positive integer $n$ for which there exists a set $S$ containing exactly $k$
nonnegative integers with
\[
\{0,1,2,\ldots,n-1\}\subseteq\{s+t:s\in S,t\in S\}.
\]
It is known that \cite{ABC1-err, ABC2-err, ABC3-err, ABC4-err, ABC5-err,
ABC6-err, ABC7-err, ABC8-err, ABC9-err, ABC0-err}
\[
0.2891\leq\operatorname*{liminf}_{k\rightarrow\infty}\frac{\nu_{2}(k)}{k^{2}%
}\leq\operatorname*{limsup}_{k\rightarrow\infty}\frac{\nu_{2}(k)}{k^{2}}%
\leq0.45504
\]
and likewise for $\nu_{j}(k)$ for $j\geq3$. See also \cite{Fi25-err}.

\textbf{2.21. Stieltjes Constants.} The number of recent articles is
staggering (see a list of references in \cite{Blgchn-err, BlgCpp-err}), more
than we can summarize here. If $d_{k}(n)$ denotes the number of sequences
$x_{1}$, $x_{2}$, ..., $x_{k}$ of positive integers such that $n=x_{1}%
x_{2}\cdots x_{k}$, then \cite{Titch-err, Sound-err, Cnry-err}
\[%
\begin{array}
[c]{ccc}%
%TCIMACRO{\dsum \limits_{n=1}^{N}}%
%BeginExpansion
{\displaystyle\sum\limits_{n=1}^{N}}
%EndExpansion
d_{2}(n)\sim N\ln(N)+(2\gamma_{0}-1)N &  & \text{(}d_{2}\text{ is the divisor
function),}%
\end{array}
\]%
\[%
%TCIMACRO{\dsum \limits_{n=1}^{N}}%
%BeginExpansion
{\displaystyle\sum\limits_{n=1}^{N}}
%EndExpansion
d_{3}(n)\sim\frac{1}{2}N\ln(N)^{2}+(3\gamma_{0}-1)N\ln(N)+(-3\gamma
_{1}+3\gamma_{0}^{2}-3\gamma_{0}+1)N,
\]%
\begin{align*}%
%TCIMACRO{\dsum \limits_{n=1}^{N}}%
%BeginExpansion
{\displaystyle\sum\limits_{n=1}^{N}}
%EndExpansion
d_{4}(n)  &  \sim\frac{1}{6}N\ln(N)^{3}+\frac{4\gamma_{0}-1}{2}N\ln
(N)^{2}+(-4\gamma_{1}+6\gamma_{0}^{2}-4\gamma_{0}+1)N\ln(N)\\
&  \ \ +(2\gamma_{2}-12\gamma_{1}\gamma_{0}+4\gamma_{1}+4\gamma_{0}%
^{3}-6\gamma_{0}^{2}+4\gamma_{0}-1)N
\end{align*}
as $N\rightarrow\infty$. More generally, $%
%TCIMACRO{\tsum \nolimits_{n=1}^{N}}%
%BeginExpansion
{\textstyle\sum\nolimits_{n=1}^{N}}
%EndExpansion
d_{k}(n)$ can be asymptotically expressed as $N$ times a polynomial of degree
$k-1$ in $\ln(N)$, which in turn can be described as the residue at $z=1$ of
$z^{-1}\zeta(z)^{k}N^{z}$. See \cite{Fi1-err} for an application of
$\{\gamma_{j}\}_{j=0}^{\infty}$ to asymptotic series for $\operatorname*{E}%
\nolimits_{n}(\omega)$ and $\operatorname*{E}\nolimits_{n}(\Omega)$,
\cite{Shrsk-err} for a generalization, and \cite{Mas1-err, Mas2-err, Mty1-err,
Mty2-err, Mtska-err, Voros-err} for connections to the Riemann hypothesis.

\textbf{2.22. Liouville-Roth Constants.} Zudilin \cite{Zud5-err} revisited the
Rhin-Viola estimate for the irrationality exponent for $\zeta(3)$.

\textbf{2.23. Diophantine Approximation Constants.} Which planar, symmetric,
bounded convex set $K$ has the worst packing density? If $K$ is a disk, the
packing density is $\pi/\sqrt{12}=0.9068996821...$, which surprisingly is
better than if $K$ is the smoothed octagon:
\[
\frac{8-4\sqrt{2}-\ln(2)}{2\sqrt{2}-1}=\frac{1}{4}%
(3.6096567319...)=0.9024141829....
\]
Do worse examples exist? A definitive answer is still at large
\cite{Scholl-err, Baez-err, Hales1-err, Hales2-err}. \ Improved lower bounds%
\[%
\begin{array}
[c]{ccc}%
\gamma_{5}\geq\frac{3}{46}\sqrt{\frac{3(9+5\sqrt{5})}{1166}}%
=0.01486...>0.0106, &  & \gamma_{6}\geq\frac{9+5\sqrt{5}}{11\sqrt{184607}%
}=0.00426...>0.0041
\end{array}
\]
are reported in \cite{Baslv0-err, Baslv1-err, Baslv2-err}.

\textbf{2.24. Self-Numbers Density Constant.} Choose $a$ to be any $r$-digit
integer expressed in base $10$ with not all digits equal. Let $a^{\prime}$ be
the integer formed by arranging the digits of $a$ in descending order, and
$a^{\prime\prime}$ be likewise with the digits in ascending order. Define
$T(a)=a^{\prime}-a^{\prime\prime}$. When $r=3$, iterates of $T$ converge to
the Kaprekar fixed point $495$; when $r=4$, iterates of $T$ converge to the
Kaprekar fixed point $6174$. For all other $r\geq2$, the situation is more
complicated \cite{Kaprekar-err, PrLudLap-err, MELines-err}. When $r=2$,
iterates of $T$ converge to the cycle $(09,81,63,27,45)$; when $r=5$, iterates
of $T$ converge to one of the following three cycles:
\[%
\begin{array}
[c]{ccccc}%
(74943,62964,71973,83952) &  & (63954,61974,82962,75933) &  & (53955,59994).
\end{array}
\]
We mention this phenomenon merely because it involves digit subtraction, while
self-numbers involved digit addition.

\textbf{2.25. Cameron's Sum-Free Set Constants.} Erd\H{o}s \cite{Erdos-err}
and Alon \&\ Kleitman \cite{AlonK-err} showed that any finite set $B$ of
positive integers must contain a sum-free subset $A$ such that $|A|>\frac
{1}{3}|B|$. See also \cite{AlonS-err, Kol2-err, Kol3-err}. The largest
constant $c$ such that $|A|>c|B|$ must satisfy $1/3\leq c<12/29$, but its
exact value is unknown. Using harmonic analysis, Bourgain \cite{Bour-err}
improved the original inequality to $|A|>\frac{1}{3}(|B|+2)$. Green
\cite{Gre2-err, Gre3-err} demonstrated that $s_{n}=O(2^{n/2})$, but the values
$c_{o}=6.8...$ and $c_{e}=6.0...$ await more precise computation.

Further evidence for the existence of complete aperiodic sum-free sets is
given in \cite{CFiFw-err, WeZhWu-err}.

\textbf{2.26. Triple-Free Set Constants.} The names for $\lambda\approx0.800$
and $\mu\approx0.613$ should be prepended by ``weakly'' and ``strongly'',
respectively. See \cite{Triple-err} for detailed supporting material. In
defining $\lambda$, the largest set $S$ such that $\forall x$
$\{x,2x,3x\}\not \subseteq S$ plays a role. The complement of $S$ in
$\{1,2,...,n\}$ is thus the smallest set $T$ such that $\forall x$
$T\cap\{x,2x,3x\}\neq\emptyset$. Clearly $T$ has size $n-p(n)$ and
$1-\lambda\approx0.199$ is the asymptotic ``hitting'' density.

\textbf{2.27. Erd\H{o}s-Lebensold Constant.} Certain variations of Erd\H{o}s'
conjecture for primitive sequences are false \cite{BFrh1-err, BFrh2-err} --
this has no bearing on the original, which remains open \cite{BksMtn-err} --
the Erd\H{o}s-Zhang conjecture for quasi-primitive sequences also requires
attention. Bounds on $M(n,k)/n$ for large $n$ and $k\geq3$ are given in
\cite{Vijay-err, Hegarty-err}. A more precise estimate $%
%TCIMACRO{\tsum }%
%BeginExpansion
{\textstyle\sum}
%EndExpansion
1/(q_{i}\ln(q_{i}))=2.0066664528...$ is now known \cite{Mathr2-err}, making
use of logarithmic integrals in \cite{HChn-err}.

\textbf{2.28. Erd\H{o}s' Sum-Distinct Set Constant.} Aliev \cite{Aliev-err}
proved that
\[
\alpha_{n}\geq\sqrt{\frac{3}{2\pi n}};
\]
Elkies \&\ Gleason's best lower bound (unpublished) is reported in
\cite{Aliev-err} to be $\sqrt{2/(\pi n)}$ rather than $\sqrt{1/n}$. Define
integer point sets $S$ and $T$ in $\mathbb{R}^{n}$ by
\[
S=\left\{  (s_{1},\ldots,s_{n}\}:s_{j}=0\text{ or }\pm1\text{ for each
}j\right\}  ,
\]%
\[
T=\left\{  (t_{1},\ldots,t_{n}\}:t_{j}=0\text{ or }\pm1\,\text{or }\pm2\text{
for each }j\right\}
\]
and let $H$ be a hyperplane in $\mathbb{R}^{n}$ such that $H\cap S$ consists
only of the origin $0$. Hence the normal vector $(a_{1},\ldots,a_{n})$ to $H
$, if each $a_{j}\in\mathbb{Z}^{+}$, has the property that $\{a_{1}%
,\ldots,a_{n}\}$ is sum-distinct. It is conjectured that \cite{AhAyK-err}
\[
\max_{H}\left\vert H\cap T\right\vert \sim c\cdot\beta^{n}%
\]
for some $c>0$ as $n\rightarrow\infty$, where $\beta=2.5386157635...$ is the
largest real zero of $x^{4}-2x^{3}-2x^{2}+2x-1$. See also \cite{BrwnM-err,
BaeCh-err}.

Fix a positive integer $n$. \ A sequence of nonnegative integers $a_{1}%
<a_{2}<\ldots<a_{k}$ is a difference basis with respect to $n$ if every
integer $0<\nu\leq n$ has a representation $a_{j}-a_{i}$; let $k(n)$ be the
minimum such $k$. The set is a \textit{restricted} difference basis if,
further, $a_{1}=0$ and $a_{k}=n$; let $\ell(n)$ be the minimum such $k$ under
these tighter constraints. \ We have \cite{DffBs1-err, DffBs2-err, DffBs3-err,
DffBs4-err, DffBs5-err}%
\[%
\begin{array}
[c]{ccc}%
2.4344\leq\lim\limits_{n\rightarrow\infty}\dfrac{k(n)^{2}}{n}\leq2.6571, &  &
2.4344\leq\lim\limits_{n\rightarrow\infty}\dfrac{\ell(n)^{2}}{n}\leq3;
\end{array}
\]
the latter may alternatively be recorded as \cite{DffBs6-err, DffBs7-err}
\[
\left(  c+o(1)\right)  \sqrt{n}\leq\ell(n)\leq\left(  \sqrt{3}+o(1)\right)
\sqrt{n}%
\]
where $c=1.5602779420...=\sqrt{2(1-\sin(\theta)/\theta}$ and $\theta$ is the
smallest positive zero of $\tan(\theta)-\theta$. \ Golay \cite{DffBs5-err}
wrote that the limiting ratio \textquotedblleft as $n\rightarrow\infty$ will,
undoubtedly, never be expressed in closed form\textquotedblright.

\textbf{2.29. Fast Matrix Multiplication Constants.} Efforts continue
\cite{FstMt1-err, FstMt2-err} to reduce the upper bound on $\omega$ to $2$.

\textbf{2.30. Pisot-Vijayaraghavan-Salem Constants.} Verger-Gaugry
\cite{VrgGgy-err} has proved that the set $T$ is bounded from below by
$1.08544$. \ Whether $\tau_{0}$ is the smallest Salem number remains open. \ 

The definition of Mahler's measure $M(\alpha)$ is unclear: It should be the
product of $\max\{1,|\alpha_{j}|\}$ over all $1\leq j\leq n$. Breusch
\cite{Breusch-err} gave a lower bound $>1$ for $M(\alpha)$ of non-reciprocal
algebraic integers $\alpha$, anticipating Smyth's stronger result by twenty years.

The sequence $\left\{  n^{1/2}\right\}  $ is uniformly distributed in $[0,1]$;
a fascinating side topic involves the gaps between adjacent points. \ A random
such gap is \textit{not} exponentially distributed but possesses a more
complicated density function. \ Elkies \&\ McMullen \cite{ElkMcM-err}
determined this density explicitly, which is piecewise analytic with phase
transistions at $1/2$ and $2$, and which has a heavy tail (implying that large
gaps are more likely than if the points were both uniform and independent).

Zudilin \cite{Zud4-err} improved Habsieger's lower bound on $(3/2)^{n}%
\operatorname*{mod}1$, progressing from $0.577^{n}$ to $0.5803^{n}$, and
similarly obtained estimates for $(4/3)^{n}\operatorname*{mod}1$ when $n$ is
suitably large. \ Concerning the latter, Pupyrev \cite{Pupy1-err, Pupy2-err}
obtained $(4/9)^{n}$ for every $n\geq2$, an important achievement.
\ Concerning the former, our desired bound $(3/4)^{n}$ for every $n\geq8$
seems out-of-reach.

Compare the sequence $\{(3/2)^{n}\}$, for which little is known, with the
recursion $x_{0}=0$, $x_{n}=\{x_{n-1}+\ln(3/2)/\ln(2)\}$, for which a musical
interpretation exists. If a guitar player touches a vibrating string at a
point two-thirds from the end of the string, its fundamental frequency is
dampened and a higher overtone is heard instead. This new pitch is a perfect
fifth above the original note. It is well-known that the \textquotedblleft
circle of fifths\textquotedblright\ never closes, in the sense that $2^{x_{n}%
}$ is never an integer for $n>0$. Further, the \textquotedblleft circle of
fifths\textquotedblright, in the limit as $n\rightarrow\infty$, fills the
continuum of pitches spanning the octave \cite{DuMc-err, HaJo-err}.

The Collatz function $f:\mathbb{Z}^{+}\rightarrow\mathbb{Z}^{+}$ is defined
by
\[
f(n)=\left\{
\begin{array}
[c]{lll}%
3n+1 &  & \text{if }n\text{ is odd}\\
n/2 &  & \text{if }n\text{ is even}%
\end{array}
\right.  .
\]
Let $f^{k}\,$denote the $k^{\text{th}}$ iterate of $f$. The $3x+1$ conjecture
asserts that, given any positive integer $n$, there exists $k$ such that
$f^{k}(n)=1$. Let $\sigma(n)$ be the first $k$ such that $f^{k}(n)<n$, called
the \textit{stopping time} of $n$. If we could demonstrate that every positive
integer $n$ has a finite stopping time, then the $3x+1$ conjecture would be
proved. Heuristic reasoning \cite{Wag-err, Cha-err, Roo-err} provides that the
average stopping time over all odd integers $1\leq n\leq N$ is asymptotically
\[
\lim_{N\rightarrow\infty}\operatorname*{E}\nolimits_{\text{odd}}(\sigma(n))=%
%TCIMACRO{\dsum \limits_{j=1}^{\infty}}%
%BeginExpansion
{\displaystyle\sum\limits_{j=1}^{\infty}}
%EndExpansion
\left\lfloor 1+\left(  1+\tfrac{\ln(3)}{\ln(2)}\right)  j\right\rfloor
c_{j}2^{-\left\lfloor \frac{\ln(3)}{\ln(2)}j\right\rfloor }=9.4779555565...
\]
where $c_{j}$ is the number of admissible sequences of order $j$. Such a
sequence $\{a_{k}\}_{k=1}^{m}$ satisfies $a_{k}=3/2$ exactly $j$ times,
$a_{k}=1/2$ exactly $m-j$ times, $%
%TCIMACRO{\tprod \nolimits_{k=1}^{m}}%
%BeginExpansion
{\textstyle\prod\nolimits_{k=1}^{m}}
%EndExpansion
a_{k}<1$ but $%
%TCIMACRO{\tprod \nolimits_{k=1}^{l}}%
%BeginExpansion
{\textstyle\prod\nolimits_{k=1}^{l}}
%EndExpansion
a_{k}>1$ for all $1\leq l<m$ \cite{Sln-err}. In contrast, the \textit{total}
stopping time $\sigma_{\infty}(n)$ of $n$, the first $k$ such that
$f^{k}(n)=1$, appears to obey
\[
\lim_{N\rightarrow\infty}\operatorname*{E}\left(  \frac{\sigma_{\infty}%
(n)}{\ln(n)}\right)  \sim\frac2{2\ln(2)-\ln(3)}=6.9521189935...=\frac
2{\ln(10)}(8.0039227796...).
\]

\textbf{2.31. Freiman's Constant.} New proofs of the Markov unicity conjecture
for prime powers $w$ appear in \cite{Frmn1-err, Frmn2-err, Frmn3-err,
Frmn4-err}. See \cite{Frmn5-err} for asymptotics for the number of admissible
triples of Diophantine equations such as
\[
u^{2}+v^{2}+2w^{2}=4uvw,
\]
\[
u^{2}+2v^{2}+3w^{2}=6uvw,
\]
\[
u^{2}+v^{2}+5w^{2}=5uvw
\]
and \cite{Frmn6-err} for mention of the constant $3.29304...$.

\textbf{2.32. De Bruijn-Newman Constant.} Ki, Kim \&\ Lee \cite{KiKmLe-err}
improved the inequality $\Lambda\leq1/2$ to $\Lambda<1/2$; building on
\cite{dBNwm1-err, dBNwm2-err, dBNwm3-err}, Rodgers \&\ Tao \cite{dBNwm4-err}
established that $\Lambda\geq0$ The constant $2\pi\,\Phi(0)=2.8066794017...$
appears in \cite{LgrsMnt-err}, in connection with a study of zeroes of the
integral of $\xi(z)$.

Further work regarding Li's criterion, which is equivalent to Riemann's
hypothesis and which involves the Stieltjes constants, appears in
\cite{Mas1-err, Mas2-err}. A different criterion is due to Matiyasevich
\cite{Mty1-err, Mty2-err}; the constant $-\ln(4\pi)+\gamma
+2=0.0461914179...=2(0.0230957089...)$ comes out as a special case. See also
\cite{Mtska-err, Vrs1-err, Vrs2-err}. As another aside, we mention the
unboundedness of $\zeta(1/2+i\,t)$ for $t\in(0,\infty)$, but that a precise
order of growth remains open \cite{Ktnk1-err, BsRc-err, Hxly05-err,
LfshtW-err}. In contrast, there is a conjecture that \cite{GrSd-err, Tchm-err,
Ktnk2-err}
\[
\max_{t\in\lbrack T,2T]}|\zeta(1+i\,t)|=e^{\gamma}\left(  \ln(\ln(T))+\ln
(\ln(\ln(T)))+C+o(1)\right)  ,
\]%
\[
\max_{t\in\lbrack T,2T]}\frac{1}{|\zeta(1+i\,t)|}=\frac{6e^{\gamma}}{\pi^{2}%
}\left(  \ln(\ln(T))+\ln(\ln(\ln(T)))+C+o(1)\right)
\]
as $T\rightarrow\infty$, where
\[
C=1-\ln(2)+%
%TCIMACRO{\dint \limits_{0}^{2}}%
%BeginExpansion
{\displaystyle\int\limits_{0}^{2}}
%EndExpansion
\frac{\ln(I_{0}(t))}{t^{2}}dt+%
%TCIMACRO{\dint \limits_{2}^{\infty}}%
%BeginExpansion
{\displaystyle\int\limits_{2}^{\infty}}
%EndExpansion
\frac{\ln(I_{0}(t))-t}{t^{2}}dt=-0.0893...
\]
and $I_{0}(t)$ is the zeroth modified Bessel function. These formulas have
implications for $|\zeta(i\,t)|$ and $1/|\zeta(i\,t)|$ as well by the analytic
continuation formula.

Looking at the sign of $\operatorname*{Re}(\zeta(1+i\,t))$ for $0\leq
t\leq10^{5} $ might lead one to conjecture that this quantity is always
positive. In fact, $t\approx682112.92$ corresponds to a negative value (the
first?) The problem can be generalized to $\operatorname*{Re}(\zeta(s+i\,t))$
for arbitrary fixed $s\geq1$. Van de Lune \cite{Lune-err, AdRBdL-err} computed
that
\[
\sigma=\sup\left\{  s\geq1:\operatorname*{Re}(\zeta(s+i\,t))<0\text{ for some
}t\geq0\right\}  =1.1923473371...
\]
is the unique solution of the equation
\[%
\begin{array}
[c]{ccc}%
%TCIMACRO{\dsum \limits_{p}}%
%BeginExpansion
{\displaystyle\sum\limits_{p}}
%EndExpansion
\arcsin\left(  p^{-\sigma}\right)  =\pi/2, &  & \sigma>1
\end{array}
\]
where the summation is over all prime numbers $p$. \ Also \cite{ReyLun-err},%
\[
x=\sup\left\{  \text{ real }s:\zeta(s+i\,t)=1\text{ for some real }t\right\}
=1.9401016837...
\]
is the unique solution $x>1$ of the equation $\zeta(x)=(2^{x}+1)/(2^{x}-1)$
and%
\[
y=\sup\left\{  \text{ real }s:\zeta^{\prime}(s+i\,t)=0\text{ for some real
}t\right\}  =2.8130140202...
\]
is the unique solution $y>1$ of the equation $\zeta^{\prime}(y)/\zeta
(y)=-2^{y+1}\ln(2)/(4^{y}-1)$.

\textbf{2.33. Hall-Montgomery Constant.} Let $\psi$ be the unique solution on
$(0,\pi)$ of the equation $\sin(\psi)-\psi\cos(\psi)=\pi/2$ and define
$K=-\cos(\psi)=0.3286741629....$ Consider any real multiplicative function $f$
whose values are constrained to $[-1,1]$. Hall \&\ Tenenbaum \cite{HT-err}
proved that, for some constant $C>0$,
\[%
\begin{array}
[c]{ccc}%
%TCIMACRO{\dsum \limits_{n=1}^{N}}%
%BeginExpansion
{\displaystyle\sum\limits_{n=1}^{N}}
%EndExpansion
f(n)\leq CN\exp\left\{  -K%
%TCIMACRO{\dsum \limits_{p\leq N}}%
%BeginExpansion
{\displaystyle\sum\limits_{p\leq N}}
%EndExpansion
\dfrac{1-f(p)}{p}\right\}  &  & \text{for sufficiently large }N,
\end{array}
\]
and that, moreover, the constant $K$ is sharp. (The latter summation is over
all prime numbers $p$.) This interesting result is a lemma used in
\cite{Hall-err}. A table of values of sharp constants $K\,$ is also given in
\cite{HT-err} for the generalized scenario where $f$ is complex, $|f|\leq1$
and, for all primes $p$, $f(p)$ is constrained to certain elliptical regions
in $\mathbb{C}$.

A fascinating coincidence involving $\delta_{0}$ is as follows. The limiting
probability that a random $n$-permutation has exactly $k$ cycles of length
exceeding $x\,n$ is \cite{Lugo1-err}%
\[
P_{0}(x)=\left\{
\begin{array}
[c]{lll}%
1-\frac{\pi^{2}}{12}+\ln(x)+\frac{1}{2}\ln(x)^{2}+\operatorname{Li}_{2}(x) &
& \text{if }\frac{1}{3}\leq x<\frac{1}{2},\\
1+\ln(x) &  & \text{if }\frac{1}{2}\leq x<1,
\end{array}
\right.
\]%
\[
P_{1}(x)=\left\{
\begin{array}
[c]{lll}%
\frac{\pi^{2}}{6}-\ln(x)-\ln(x)^{2}-2\operatorname{Li}_{2}(x) &  & \text{if
}\frac{1}{3}\leq x<\frac{1}{2},\\
-\ln(x) &  & \text{if }\frac{1}{2}\leq x<1,
\end{array}
\right.
\]%
\[
P_{2}(x)=\left\{
\begin{array}
[c]{lll}%
-\frac{\pi^{2}}{12}+\frac{1}{2}\ln(x)^{2}+\operatorname{Li}_{2}(x) &  &
\text{if }\frac{1}{3}\leq x<\frac{1}{2},\\
0 &  & \text{if }\frac{1}{2}\leq x<1
\end{array}
\right.
\]
as $n\rightarrow\infty$, where $k=0,1,2$. The value of $x$ that maximizes
$P_{1}(x)$ is $\xi=1/\left(  1+\sqrt{e}\right)  =0.3775406687...$; we have%
\[
P_{1}(\xi)=1-\delta_{0}=0.8284995068...,
\]
$P_{0}(\xi)=0.0987117544...$, $P_{2}(\xi)=0.0727887386...$ (which are
non-Poissonian). In particular, most $n$-permutations have \textit{exactly
one} cycle longer than $\xi\,n$.

\textbf{3.1. Shapiro-Drinfeld Constant}. A\ construction involving the
smallest concave down function $\geq$ prescribed data appears in
\cite{Fi36-err}.

\textbf{3.2. Carlson-Levin Constants}. Various generalizations appear in
\cite{LrsPrs1-err, LrsPrs2-err, LrsPrs3-err, Oskwsk-err}; analogous sharp
constants for finite series remain open, as for integrals over bounded regions.

\textbf{3.3. Landau-Kolmogorov Constants.} For $L_{2}(0,\infty)$, Bradley
\&\ Everitt \cite{BrEv-err} were the first to determine that
$C(4,2)=2.9796339059...=\sqrt{8.8782182137...}$; see also \cite{Russ-err,
Phng-err, Beyn-err}. Ditzian \cite{Ditz-err} proved that the constants for
$L_{1}(-\infty,\infty)$ are the same as those for $L_{\infty}(-\infty
,\infty).$ Ph\'ong \cite{Phng-err} obtained the following best possible
inequality in $L_{2}(0,1)$:
\[%
%TCIMACRO{\dint \limits_{0}^{1}}%
%BeginExpansion
{\displaystyle\int\limits_{0}^{1}}
%EndExpansion
\left|  f^{\prime}(x)\right|  ^{2}dx\leq(6.4595240299...)\left(
%TCIMACRO{\dint \limits_{0}^{1}}%
%BeginExpansion
{\displaystyle\int\limits_{0}^{1}}
%EndExpansion
\left|  f(x)\right|  ^{2}dx+%
%TCIMACRO{\dint \limits_{0}^{1}}%
%BeginExpansion
{\displaystyle\int\limits_{0}^{1}}
%EndExpansion
\left|  f^{\prime\prime}(x)\right|  ^{2}dx\right)
\]
where the constant is given by $\sec(2\theta)/2$ and $\theta$ is the unique
zero satisfying $0<\theta<\pi/4$ of
\begin{align*}
&  \ \ \ \sin(\theta)^{4}\left(  e^{2\sin(\theta)}-1\right)  ^{2}%
(e^{-2\sin(\theta)}-1)^{2}+\cos(\theta)^{4}[2-2\cos(2\cos(\theta))]^{2}\\
&  \ \ \ -\cos(2\theta)^{4}[1+e^{4\sin(\theta)}-2e^{2\sin(\theta)}\cos
(2\cos(\theta))][1+e^{-4\sin(\theta)}-2e^{-2\sin(\theta)}\cos(2\cos
(\theta))]\\
&  \ \ \ -2\cos(\theta)^{2}\sin(\theta)^{2}[2-2\cos(2\cos(\theta
))](1-e^{-2\sin(\theta)})\left(  e^{2\sin(\theta)}-1\right)  .
\end{align*}
We wonder about other such additive analogs of Landau-Kolmogorov inequalities.

\textbf{3.4. Hilbert's Constants.} Borwein \cite{JBrwn-err} mentioned the case
$p=q=4/3$ and $\lambda=1/2$, which evidently remains open. Peachey
\&\ Enticott \cite{TPchy-err} performed relevant numerical experiments. See
also \cite{AbHny-err}.

\textbf{3.5. Copson-de Bruijn Constant.} An English translation of
Ste\v{c}kin's paper is available \cite{Stckn-err}. Ackermans \cite{Ack-err}
studied the recurrence $\{u_{n}\}$ in greater detail. \ Let $\Omega$ be a
domain in $\mathbb{R}^{n}$ and let $p>1$. A multidimensional version of
Hardy's inequality is \cite{BFT-err}
\[%
%TCIMACRO{\dint \limits_{\Omega}}%
%BeginExpansion
{\displaystyle\int\limits_{\Omega}}
%EndExpansion
\left\vert \nabla f(x)\right\vert ^{p}dx\geq\left\vert \frac{n-p}%
{p}\right\vert ^{p}%
%TCIMACRO{\dint \limits_{\Omega}}%
%BeginExpansion
{\displaystyle\int\limits_{\Omega}}
%EndExpansion
\frac{|f(x)|^{p}}{|x|^{p}}dx
\]
and the constant is sharp. Let $\delta(x)$ denote the (shortest) distance
between $x$ and the boundary $\partial\Omega$ of $\Omega$. A variation of
Hardy's inequality is
\[%
%TCIMACRO{\dint \limits_{\Omega}}%
%BeginExpansion
{\displaystyle\int\limits_{\Omega}}
%EndExpansion
\left\vert \nabla f(x)\right\vert ^{p}dx\geq\left(  \frac{p-1}{p}\right)  ^{p}%
%TCIMACRO{\dint \limits_{\Omega}}%
%BeginExpansion
{\displaystyle\int\limits_{\Omega}}
%EndExpansion
\frac{|f(x)|^{p}}{\delta(x)^{p}}dx
\]
assuming $\Omega$ is a convex domain with smooth boundary. Again, the constant
is sharp. With regard to the latter inequality, let $n=2$, $p=2$ and
$\Omega=\Omega_{\alpha}$ be the nonconvex plane sector of angle $\alpha$:
\[
\Omega_{\alpha}=\left\{  r\,e^{i\,\theta}:0<r<1\text{ and }0<\theta
<\alpha\right\}  .
\]
Davies \cite{Dvs-err} demonstrated that the reciprocal of the best constant
is
\[
\left\{
\begin{array}
[c]{lll}%
4 &  & \text{if }0<\alpha<4.856...,\\
>4 &  & \text{if }4.856...<\alpha<2\pi,\\
4.869... &  & \text{if }\alpha=2\pi
\end{array}
\right.
\]
and Tidblom \cite{Tdb-err} found that the threshold angle is exactly
\begin{align*}
\alpha &  =\pi+4\arctan\left(  4\frac{\Gamma(3/4)^{2}}{\Gamma(1/4)^{2}}\right)
\\
\  &  =\pi+4\arctan\left(  \frac{1}{2}\frac{3^{2}-1}{3^{2}}\frac{5^{2}}%
{5^{2}-1}\frac{7^{2}-1}{7^{2}}\cdot\cdots\right)  =4.8560553209...
\end{align*}
A similar expression for $4.869...$ remains open.

\textbf{3.6. Sobolev Isoperimetric Constants.} In section 3.6.1,
$\sqrt{\lambda}=1$ represents the principal frequency of the sound we hear
when a string is plucked; in section 3.6.3, $\sqrt{\lambda}=\theta$ represents
likewise when a kettledrum is struck. (The square root was missing in both.)
The units of frequency, however, are not compatible between these two examples.

The \textquotedblleft rod \textquotedblright constant
$500.5639017404...=(4.7300407448...)^{4}$ appears in \cite{Rod1-err, Rod2-err,
Rod3-err}. It is the second term in a sequence $c_{1}$, $c_{2}$, $c_{3}$, ...
for which $c_{1}=\pi^{2}=9.869...$ (in connection with the \textquotedblleft
string\textquotedblright\ inequality) and $c_{3}=(2\pi)^{6}=61528.908...$; the
constant $c_{4}$ is the smallest eigenvalue of ODE
\[%
\begin{array}
[c]{ccc}%
f^{(viii)}(x)=\lambda\,f(x), &  & 0\leq x\leq1,
\end{array}
\]%
\[%
\begin{array}
[c]{ccc}%
f(0)=f^{\prime}(0)=f^{\prime\prime}(0)=f^{\prime\prime\prime}(0)=0, &  &
f(1)=f^{\prime}(1)=f^{\prime\prime}(1)=f^{\prime\prime\prime}(1)=0
\end{array}
\]
and was computed by Abbott \cite{Rod4-err} to be $(7.8187073432...)^{8}%
=(1.3966245157...)\times10^{7}$. Allied subjects include positive definite
Toeplitz matrices and conditioning of certain least squares problems.

If we eliminate the vanishing requirements on $f$ at $x=1$, then the
\textquotedblleft string\textquotedblright\ constant becomes $4/\pi^{2}$ and
the \textquotedblleft rod\textquotedblright\ constant becomes \cite{Rod5-err,
Rod6-err, Rod7-err}%
\[
\mu=\frac{1}{12.3623633683...}=(0.2844128718...)^{2}=\frac{1}{\theta^{4}}%
\]
where $\theta=1.8751040687...$ is the smallest positive root of%
\[
\cos(\theta)\cosh(\theta)=-1\text{.}%
\]
Under $f(0)=f^{\prime}(0)=f^{\prime\prime}(0)=0$, we have%
\[%
%TCIMACRO{\dint \limits_{0}^{1}}%
%BeginExpansion
{\displaystyle\int\limits_{0}^{1}}
%EndExpansion
f(x)^{2}dx\leq\mu%
%TCIMACRO{\dint \limits_{0}^{1}}%
%BeginExpansion
{\displaystyle\int\limits_{0}^{1}}
%EndExpansion
\left(  \frac{d^{3}f}{dx^{3}}\right)  ^{2}dx
\]
where%
\[
\mu=\frac{1}{121.2590618589...}=(0.0908119283...)^{2}=\frac{1}{\theta^{6}}%
\]
is best possible and $\theta=2.2247729764...$ is the smallest positive root of%
\[
8\cos(\theta)+\cos(2\theta)+16\cos(\theta/2)\cosh\left(  \sqrt{3}%
\theta/2\right)  +2\cos(\theta)\cosh\left(  \sqrt{3}\theta\right)  =-9.
\]
What is the corresponding equation when additionally $f(1)=f^{\prime
}(1)=f^{\prime\prime}(1)=0$?

Here is another example \cite{Boyd-err, SkrSt-err}: the best constant $K$ for
the inequality%
\[%
\begin{array}
[c]{ccc}%
%TCIMACRO{\dint \limits_{0}^{\pi}}%
%BeginExpansion
{\displaystyle\int\limits_{0}^{\pi}}
%EndExpansion
g(x)^{2}g^{\prime}(x)^{2}dx\leq\left(  \dfrac{\pi}{2}\right)  ^{2}K%
%TCIMACRO{\dint \limits_{0}^{\pi}}%
%BeginExpansion
{\displaystyle\int\limits_{0}^{\pi}}
%EndExpansion
g^{\prime}(x)^{4}dx, &  & g(0)=g(\pi)=0
\end{array}
\]
is $K=2/(L+1)^{2}=0.3461189656...$, where%
\[
L=%
%TCIMACRO{\dint \limits_{0}^{1}}%
%BeginExpansion
{\displaystyle\int\limits_{0}^{1}}
%EndExpansion
\frac{1}{1-\frac{2}{3}t^{2}}\,dt=\sqrt{\frac{3}{2}}\operatorname{arctanh}%
\left(  \sqrt{\frac{2}{3}}\right)  =1.4038219651....
\]
More relevant material is found in \cite{Fi3-err, Fi4-err, Fi13-err}. See
\cite{MrsPzc-err} for a variation involving the norm of a product $f\,g$,
bounded by the product of the norms of $f$ and $g$.

\textbf{3.7. Korn Constants.} \ A\ closed-form expression for even the
smallest Laplacian eigenvalue $7.1553391339...$ \cite{HHrr-err} over a regular
hexagon is unavailable.

\textbf{3.8. Whitney-Mikhlin Extension Constants}. For completeness' sake, we
mention that%
\[%
\begin{array}
[c]{lllll}%
\chi_{2}=\sqrt{\frac{1}{I_{1}(1)K_{0}(1)}}, &  & \chi_{4}=\sqrt{\frac
{1}{\left(  I_{0}(1)-2I_{1}(1)\right)  K_{1}(1)}}, &  & \chi_{6}=\sqrt
{\frac{1}{\left(  9I_{1}(1)-4I_{0}(1)\right)  \left(  2K_{1}(1)+K_{0}%
(1)\right)  }}%
\end{array}
\]
via recursions for modified Bessel functions.

\textbf{3.9. Zolotarev-Schur Constant}. Explicit expressions for $s_{n}$ \&
$Z_{n}(x)$ are known for $4\leq n\leq7$, thanks to Rack \cite{Rack5-err,
Rack6-err} and Rack \&\ Vajda \cite{Rack7-err, Rack8-err}.

Here is a different problem involving approximation over an ellipse $E$. \ We
assume that $E$ possesses foci $\pm1$ and sum of semi-axes equal to $1/q$,
where $0<q<1$. \ Let $f(z)$ be analytic in the interior of $E$, real-valued
along the major axis of $E$, and bounded in the sense that $\left\vert
\operatorname{Re}(f(z))\right\vert \leq1$ in the interior of $E$. \ Then the
best approximation of $f(z)$ on $[-1,1]$ by a polynomial of degree $n-1$ has
error at most%
\[
\frac{8}{\pi}%
%TCIMACRO{\dsum \limits_{k=0}^{\infty}}%
%BeginExpansion
{\displaystyle\sum\limits_{k=0}^{\infty}}
%EndExpansion
\frac{(-1)^{k}}{2k+1}\frac{q^{(2k+1)n}}{1+q^{2(2k+1)n}}.
\]
Further, there exists an $f(z)$ for which equality is attained, that is, the
Favard-like constant (in $q$) is sharp \cite{Achsr1-err, Achsr2-err,
Achsr3-err}.

\textbf{3.10. Kneser-Mahler Constants.} The constants $\ln(\delta)$ and
$\ln(\beta)$ appear in \cite{ABYZ-err}; we also have%
\[%
\begin{array}
[c]{ccc}%
\pi\ln(\delta)=2\operatorname{Im}\left(  \operatorname*{Li}_{2}\left(
i\right)  \right)  , &  & \pi\ln(\beta)=\operatorname{Im}\left(
\operatorname*{Li}_{2}\left(  \dfrac{1+i\sqrt{3}}{2}\right)  \right)
\end{array}
\]
from \cite{CLvN-err}. Conjectured L-series expressions for $M\left(  1+%
%TCIMACRO{\tsum \nolimits_{j=1}^{n}}%
%BeginExpansion
{\textstyle\sum\nolimits_{j=1}^{n}}
%EndExpansion
x_{j}\right)  $, due to Rodriguez-Villegas, are exhibited for $n=4$, $5$ in
\cite{Fi17-err}.

\textbf{3.11. Grothendieck's Constants}. \ It is now known \cite{BrMNr1-err,
BrMNr2-err} that $\kappa_{R}<\pi/\left(  2\ln(1+\sqrt{2})\right)
-\varepsilon$ for some explicit $\varepsilon>0$; a similar result for
$\kappa_{C}$ remains open. \ See \cite{RghvS-err, Pisier-err} for connections
with theoretical computer science and quantum physics.

\textbf{3.12. Du Bois Reymond's Constants.} The smallest positive solution
$4.4934094579...$ of the equation $\tan(x)=x$ appears in \cite{DffBs2-err}; it
is also the smallest positive local minimum of $\sin(x)/x$. The constant
$(\pi/\xi)^{2}$ is equal to the largest eigenvalue of the infinite symmetric
matrix $(a_{m,n})_{m\geq1,n\geq1}$ with elements $a_{m,n}=m^{-1}n^{-1}%
+m^{-2}\delta_{m,n}$, where $\delta_{m,n}=1$ if $m=n$ and $\delta_{m,n}=0$.
Boersma \cite{Boe-err} employed this fact to give an alternative proof of
Szeg\"{o}'s theorem. \ Let $\eta_{0}$ be the positive solution of
$\tanh(1/x)=x$ and $\eta_{1}$, $\eta_{2}$, $\eta_{3}$, ... be all positive
solutions of $\tan(1/x)=-x$. \ We have \cite{MhtNrm-err}%
\[%
\begin{array}
[c]{ccc}%
\eta_{0}^{4}+%
%TCIMACRO{\dsum \limits_{k=1}^{\infty}}%
%BeginExpansion
{\displaystyle\sum\limits_{k=1}^{\infty}}
%EndExpansion
\eta_{k}^{4}=\dfrac{1}{2}, &  & \eta_{0}^{6}-%
%TCIMACRO{\dsum \limits_{k=1}^{\infty}}%
%BeginExpansion
{\displaystyle\sum\limits_{k=1}^{\infty}}
%EndExpansion
\eta_{k}^{6}=\dfrac{1}{3}%
\end{array}
\]
and much more.

\textbf{3.13. Steinitz constants}. \ We hope to report on \cite{Barany-err,
Bnszyk-err} later.

\textbf{3.14. Young-Fej\'er-Jackson Constants}. The quantity $0.3084437795...
$, called Zygmund's constant, would be better named after
Littlewood-Salem-Izumi \cite{Zygm1-err, Zygm2-err, Zygm3-err, Zygm4-err,
Zygm5-err}.

\textbf{3.15. Van der Corput's Constant.} We examined only the case in which
$f$ is a real twice-continuously differentiable function on the interval
$[a,b]$; a generalization to the case where $f$ is $n$ times differentiable,
$n\geq2$, is discussed in \cite{AKC-err, Rog-err} with some experimental
numerical results for $n=3$.

\textbf{3.16. Tur\'{a}n's Power Sum Constants.} Recent work appears in
\cite{Tu1-err, Tu2-err, Tu3-err, Tu4-err, Tu5-err, Tu6-err, Tu7-err, Tu8-err},
to be reported on later.

\textbf{4.1. Gibbs-Wilbraham Constant}. On the one hand, Gibbs' constant for a
jump discontinuity for Fourier-Bessel partial sums seems to be numerically
equal to that for ordinary Fourier partial sums (a proof is not given in
\cite{FayKlp-err}). On the other hand, the analog of $(2/\pi)G$ corresponding
to de la Vall\'{e}e Poussin sums is%
\[%
%TCIMACRO{\dint \limits_{0}^{2\pi/3}}%
%BeginExpansion
{\displaystyle\int\limits_{0}^{2\pi/3}}
%EndExpansion
\frac{\cos(\theta)-\cos(2\theta)}{\theta^{2}}d\theta=1.1427281269...
\]
which is slightly less than $1.1789797444...$\cite{ByrGoh-err}. \ It is
possible to generalize the classical case to piecewise smooth functions $f$
for which the jump discontinuity occurs not for $f$, but rather for its
derivative. \ The lowest undershooting corresponding to such `kinks' is
$\cos(\xi)=-0.3482010120...$ where $\xi=1.9264476603...$ is the smallest
positive root of
\[
x%
%TCIMACRO{\dint \limits_{x}^{\infty}}%
%BeginExpansion
{\displaystyle\int\limits_{x}^{\infty}}
%EndExpansion
\frac{\cos(u)}{u^{2}}du=\cos(x).
\]
This phenomenon, although more subtle than the usual scenario, deserves to be
better known \cite{ByrGoh-err}.

\textbf{4.2. Lebesgue Constants}. \ Asymptotic expansions (in terms of
negative integer powers of $n+1$) for $G_{n}$ and $L_{n/2}$ appear in
\cite{LLZh-err, LLNe-err, LLCh-err}. \ If, for $n=4$, we restrict $x_{1}=-1$,
$x_{4}=1$ and $x_{2}=-x_{3}$, then the smallest $\Lambda_{4}$ corresponds to
$x_{3}^{\ast}=0.4177913013...$ with minimal polynomial $25z^{6}+17z^{4}%
+2z^{2}-1$; it also has value $\Lambda_{4}^{\ast}=1.4229195732...$ with
minimal polynomial $43w^{3}-93w^{2}+53w-11$. \ In contrast, $\Lambda_{2}%
^{\ast}=1$ and $\Lambda_{3}^{\ast}=5/4$ trivially, but $\Lambda_{5}^{\ast
}=1.5594902098...$ nontrivially with minimal polynomial of degree $73$
\cite{Rack1-err, Rack2-err, Rack3-err, Rack4-err}.

\textbf{4.3. Achieser-Krein-Favard Constants.} An English translation of
Nikolsky's work is available \cite{Nklsk1-err}. While on the subject of
trigonometric polynomials, we mention Littlewood's conjecture \cite{Erde-err}.
Let $n_{1}<n_{2}<\ldots<n_{k}$ be integers and let $c_{j}$, $1\leq j\leq k $,
be complex numbers with $|c_{j}|\geq1$. Konyagin \cite{Kony-err} and McGehee,
Pigno \&\ Smith \cite{MPS-err} proved that there exists $C>0$ so that the
inequality
\[%
%TCIMACRO{\dint \limits_{0}^{1}}%
%BeginExpansion
{\displaystyle\int\limits_{0}^{1}}
%EndExpansion
\left\vert
%TCIMACRO{\dsum \limits_{j=1}^{k}}%
%BeginExpansion
{\displaystyle\sum\limits_{j=1}^{k}}
%EndExpansion
c_{j}e^{2\pi in_{j}\xi}\right\vert d\xi\geq C\ln(k)
\]
always holds. It is known that the smallest such constant $C$ satisfies
$C\leq4/\pi^{2}$; Stegeman \cite{Steg-err} demonstrated that $C\geq0.1293$ and
Yabuta \cite{Yabu-err} improved this slightly to $C\geq0.129590$. What is the
true value of $C$?

\textbf{4.4. Bernstein's Constant}. \ Consider more generally the case
$f(x)=|x|^{s}$ and $B(s)=\lim_{n\rightarrow\infty}n^{s}E_{n}(f)$ for $s>0$,
where the error is quantified in $L_{\infty}[-1,1]$.\ Although we know $B(1) $
to high precision, no explicit expression for it (or for $B(s)$ when $s\neq1$)
is known. \ In contrast, the $L_{1}$ and $L_{2}$ analogs of $B(s)$ are
\cite{Nklsk2-err, Raitsin-err, Gnzbrg-err, Lbnsky-err}%
\[%
\begin{array}
[c]{lll}%
(8/\pi)|\sin(s\pi/2)|\Gamma(s+1)\beta(s+2), &  & \left(  2/\sqrt{\pi}\right)
|\sin(s\pi/2)|\Gamma(s+1)\sqrt{1/(2s+1)}%
\end{array}
\]
respectively, where $\beta(z)$ is Dirichlet's beta function. \ Also
\cite{Stahl0-err}%
\[
\lim_{n\rightarrow\infty}e^{\pi\sqrt{sn}}E_{n,n}(f)=4^{1+s/2}|\sin(s\pi/2)|
\]
which reduces to $8$ in special circumstance $s=1$. \ 

\textbf{4.5. The \textquotedblleft One-Ninth\textquotedblright\ Constant.}
Zudilin \cite{Zud3-err} deduced that $\Lambda$ is transcendental by use of
Theorem 4 in \cite{Bert-err}. See also \cite{Magnus-err, Aptkrv-err,
BrgErk-err}.

\textbf{4.6. Frans\'{e}n-Robinson Constant}. \ For thoroughness' sake, we give
moments%
\[%
\begin{array}
[c]{lll}%
\dfrac{1}{I}%
%TCIMACRO{\dint \limits_{0}^{\infty}}%
%BeginExpansion
{\displaystyle\int\limits_{0}^{\infty}}
%EndExpansion
\dfrac{x}{\Gamma(x)}dx=1.9345670421..., &  & \dfrac{1}{I}%
%TCIMACRO{\dint \limits_{0}^{\infty}}%
%BeginExpansion
{\displaystyle\int\limits_{0}^{\infty}}
%EndExpansion
\dfrac{x^{2}}{\Gamma(x)}dx=4.8364859746...
\end{array}
\]
of the reciprocal gamma distribution (not to be confused with the
\textit{inverse} gamma distribution).

\textbf{4.7. Berry-Esseen Constant.} The upper bound for $C$ can be improved
to $0.4785$ when $X_{1}$, $X_{2}$, $\ldots$, $X_{n}$ are identically
distributed \cite{ITyurin-err, Shvtk1-err} and to $0.5600$ when
non-identically distributed \cite{Shvtk2-err}. A different form of the
inequality is found in \cite{Fi5-err}.

\textbf{4.8. Laplace Limit Constant.} The quantity $\lambda=0.6627434193...$
appears in \cite{Knuth-err} with regard to Knuth's whirlpool permutations and
\cite{Fi14-err} wrt Plateau's problem for two circular rings dipped in soap
solution; $\mu=\sqrt{\lambda^{2}+1}$ appears in \cite{Sluciak-err} wrt solving
an exponential equation. Definite integral expressions include
\cite{Nicolau-err, Cantrell-err}
\[
\mu=1+\frac{%
%TCIMACRO{\tint \nolimits_{0}^{2\pi}}%
%BeginExpansion
{\textstyle\int\nolimits_{0}^{2\pi}}
%EndExpansion
\frac{e^{2i\,\theta}d\theta}{\coth(e^{i\,\theta}+1)-e^{i\,\theta}-1}}{%
%TCIMACRO{\tint \nolimits_{0}^{2\pi}}%
%BeginExpansion
{\textstyle\int\nolimits_{0}^{2\pi}}
%EndExpansion
\frac{e^{i\,\theta}d\theta}{\coth(e^{i\,\theta}+1)-e^{i\,\theta}-1}}%
=\sqrt{\frac{1-\frac{1}{2}%
%TCIMACRO{\tint \nolimits_{-1}^{1}}%
%BeginExpansion
{\textstyle\int\nolimits_{-1}^{1}}
%EndExpansion
\frac{t^{2}dt}{(t-\operatorname*{arctanh}(t))^{2}+\pi^{2}/4}}{1-\frac{1}{2}%
%TCIMACRO{\tint \nolimits_{-1}^{1}}%
%BeginExpansion
{\textstyle\int\nolimits_{-1}^{1}}
%EndExpansion
\frac{dt}{(t-\operatorname*{arctanh}(t))^{2}+\pi^{2}/4}}}.
\]
Also, $\sinh(\mu)=1.5088795615...$ occurs in asymptotic combinatorics and as
an extreme result in complex analysis \cite{Lap01-err, Lap02-err, Lap03-err,
Lap04-err}; $\sinh(\mu)/\mu=1.2577364561...$ occurs when minimizing the
maximum tension of a heavy cable spanning two points of equal height
\cite{Lap05-err}.

Let $c>0$. The boundary value problem
\[%
\begin{array}
[c]{ccc}%
y^{\prime\prime}(x)+c\,e^{y(x)}=0, &  & y(0)=y(1)=0
\end{array}
\]
has zero, one or two solutions when $c>\gamma$, $c=\gamma$ and $c<\gamma$,
respectively; the critical threshold
\[
\gamma=8\lambda^{2}=3.5138307191...=4(0.8784576797...)
\]
was found by Bratu \cite{Lplc1-err, Lplc2-err} and Frank-Kamenetskii
\cite{Lplc3-err, Lplc4-err}. Another way of expressing this is that the
largest $\beta>0$ for which
\[%
\begin{array}
[c]{ccc}%
y^{\prime\prime}(x)+e^{y(x)}=0, &  & y(0)=y(\beta)=0
\end{array}
\]
possesses a solution is $\beta=\sqrt{8}\lambda=1.8745214640....$ Under the
latter circumstance, it follows that
\[
y^{\prime}(0)=\sqrt{2}\sinh(\mu)=2.1338779399...=\sqrt{2(\delta-1)}%
\]
where $\delta=\cosh(\mu)^{2}=3.2767175312...$. These differential equations
are useful in modeling thermal ignition and combustion \cite{Lplc5-err,
Lplc6-err, Lplc7-err, Lplc8-err}; see \cite{Fi24-err} for similar equations
arising in astrophysics.

\textbf{4.9. Integer Chebyshev Constant}. \ The bounds $0.4213<\chi
(0,1)<0.422685$ are currently best known \cite{Prtskr1-err, Mchsnr-err,
Flmng1-err, Flmng2-err}. \ Other values of $\chi(a,b)$ and various techniques
are studied in \cite{Hare0-err}. If the integer polynomials are constrained to
be monic, then a different line of research emerges \cite{BPPrt-err,
HrSmy-err, Hlmr0-err}. Consider instead the class $S_{n}$~of all integer
polynomials of the exact degree $n$ and all $n$ zeroes both in $[-1,1]$ and
simple.\ Let%
\[%
\begin{array}
[c]{ccccc}%
%TCIMACRO{\dsum \limits_{k=0}^{n}}%
%BeginExpansion
{\displaystyle\sum\limits_{k=0}^{n}}
%EndExpansion
a_{k,n}x^{n}\in S_{n}, &  & a_{n,n}\neq0, &  & n=1,2,3,\ldots
\end{array}
\]
be an arbitrary sequence $R$ of polynomials. \ Building on work of Schur
\cite{Schur0-err}, Pritsker \cite{Prtskr2-err} demonstrated that%
\[
1.5381<\frac{1}{\sqrt{\chi(0,1)}}\leq\inf_{R}\operatorname*{liminf}%
\limits_{n\rightarrow\infty}\left\vert a_{n,n}\right\vert ^{1/n}<1.5417
\]
(his actual lower bound $1.5377$ used $\chi(0,1)<0.42291334$ from
\cite{Flmng1-err}; we use the refined estimate from \cite{Flmng2-err}).
\ A\ follow-up essay on real transfinite diameter is \cite{Fi37-err}.

\textbf{5.1. Abelian Group Enumeration Constants.} Asymptotic expansions for $%
%TCIMACRO{\tsum \nolimits_{n\leq N}}%
%BeginExpansion
{\textstyle\sum\nolimits_{n\leq N}}
%EndExpansion
a(n)^{m}$ are possible for any integer $m\geq2$ \cite{ZhLuZh-err, TothAB-err}.
For a finite abelian group $G$, let $r(G)$ denote the minimum number of
generators of $G$ and let $E(G)$ denote the expected number of random elements
from $G$, drawn independently and uniformly, to generate $G$. Define
$e(G)=E(G)-r(G)$, the \textit{excess} of $G$. Then \cite{Pmrnc-err}
\[
e_{r}=\sup\left\{  e(G):r(G)=r\right\}  =1+%
%TCIMACRO{\dsum \limits_{j=1}^{\infty}}%
%BeginExpansion
{\displaystyle\sum\limits_{j=1}^{\infty}}
%EndExpansion
\left(  1-%
%TCIMACRO{\dprod \limits_{k=1}^{r}}%
%BeginExpansion
{\displaystyle\prod\limits_{k=1}^{r}}
%EndExpansion
\zeta(j+k)^{-1}\right)  ;
\]
in particular, $e_{1}=1.7052111401...$ (Niven's constant) for the cyclic case
and
\[
\sigma=\lim_{r\rightarrow\infty}e_{r}=1+%
%TCIMACRO{\dsum \limits_{j=2}^{\infty}}%
%BeginExpansion
{\displaystyle\sum\limits_{j=2}^{\infty}}
%EndExpansion
\left(  1-%
%TCIMACRO{\dprod \limits_{k=j}^{\infty}}%
%BeginExpansion
{\displaystyle\prod\limits_{k=j}^{\infty}}
%EndExpansion
\zeta(k)^{-1}\right)  =2.118456563...
\]
in general. It is remarkable that this limit is finite! Let also
\[
\tau=%
%TCIMACRO{\dsum \limits_{j=1}^{\infty}}%
%BeginExpansion
{\displaystyle\sum\limits_{j=1}^{\infty}}
%EndExpansion
\left(  1-\left(  1-2^{-j}\right)
%TCIMACRO{\dprod \limits_{k=j+1}^{\infty}}%
%BeginExpansion
{\displaystyle\prod\limits_{k=j+1}^{\infty}}
%EndExpansion
\zeta(k)^{-1}\right)  =1.742652311...,
\]
then for the multiplicative group $\mathbb{Z}_{n}^{\ast}$ of integers
relatively prime to $n$,
\[
\sup\left\{  e(G):G=\mathbb{Z}_{n}^{\ast}\text{ and }2<n\equiv
l\operatorname*{mod}8\right\}  =\left\{
\begin{array}
[c]{lll}%
\sigma &  & \text{if }l=1,3,5\text{ or }7,\\
\sigma-1 &  & \text{if }l=2\text{ or }6,\\
\tau &  & \text{if }l=4,\\
\tau+1 &  & \text{if }l=0.
\end{array}
\right.
\]
We emphasize that $l$, not $n$, is fixed in the supremum (as according to the
right-hand side). The constant $A_{1}^{-1}=0.4357570767...$ makes a small
appearence (as a certain \textquotedblleft best probability\textquotedblright%
\ corresponding to finite nilpotent groups).

Let $\mathbb{Z}^{n}$ denote the additive group of integer $n$-vectors (free
abelian group of rank $n$) and $M_{n}(\mathbb{Z})$ denote the ring of integer
$n\times n$ matrices. From a different point of view, we have
\cite{KrvMzP-err}%
\[
\operatorname*{P}\left\{  m\text{ random }n\text{-vectors generate }%
\mathbb{Z}^{n}\right\}  =\left\{
\begin{array}
[c]{ccc}%
0 &  & \text{if }m=n,\\
\tfrac{1}{\zeta(m-n+1)}\tfrac{1}{\zeta(m-n+2)}\cdots\tfrac{1}{\zeta(m)} &  &
\text{if }m>n,
\end{array}
\right.
\]%
\[
\operatorname*{P}\left\{  m\text{ random }2\times2\text{ matrices generate
}M_{2}(\mathbb{Z})\right\}  =\left\{
\begin{array}
[c]{ccc}%
0 &  & \text{if }m=2,\\
\frac{1}{\zeta(m-1)\zeta(m)} &  & \text{if }m>2,
\end{array}
\right.
\]%
\[
\operatorname*{P}\left\{  2\text{ random }3\times3\text{ matrices generate
}M_{3}(\mathbb{Z})\right\}  =\tfrac{1}{\zeta(2)^{2}\zeta(3)},
\]%
\[
\operatorname*{P}\left\{  3\text{ random }3\times3\text{ matrices generate
}M_{3}(\mathbb{Z})\right\}  =\tfrac{1}{\zeta(2)\zeta(3)\zeta(4)}%
%TCIMACRO{\dprod \limits_{p}}%
%BeginExpansion
{\displaystyle\prod\limits_{p}}
%EndExpansion
\left(  1+\tfrac{1}{p^{2}}+\tfrac{1}{p^{3}}-\tfrac{1}{p^{5}}\right)  .
\]
It is surprising that two $2\times2$ matrices differ from two $3\times3$
matrices in this regard (the former probability is zero but the latter is
positive!) See \cite{BEoB-err, CDoB-err} for more on nonabelian group enumeration.

\textbf{5.2. Pythagorean Triple Constants.} Improvements in estimates for
$P_{a}(n)$ and $P_{p}(n)$ are found in \cite{Pythg1-err, Pythg2-err}. Let
$P_{\ell}(n)$ denote the number of primitive Pythagorean triangles under the
constraint that the two legs are both $\leq n$; then \cite{Pythg3-err}%
\[
P_{\ell}(n)=\frac{4}{\pi^{2}}\ln\left(  1+\sqrt{2}\right)  n+O\left(  \sqrt
{n}\right)
\]
as $n\rightarrow\infty$. The quantity $H_{h}(n)$ should be defined as the
number of primitive Heronian triangles under the constraint that all three
sides are $\leq n$. A better starting point for studying $H_{a}^{\prime}(n)$
might be \cite{Heron1-err, Heron2-err, Heron3-err, Heron4-err}.

\textbf{5.3. R\'{e}nyi's Parking Constant.} Expressions similar to those for
$M(x)$, $m$ and $v$ appear in the analysis of a certain stochastic
fragmentation process \cite{Agch-err}. \ More constants appear in the jamming
limit of arbitrary graphs; for example, $0.3641323...$ and $0.3791394...$
correspond respectively to the square and hexagonal lattices \cite{BJkMy-err}.

Consider monomers on $1\times\infty$ that exclude $s$ neighbors on both right
and left sides. \ The expected density of cars parked on the lattice is
\cite{Seat1-err, Seat2-err, Seat3-err, Seat4-err}%
\[%
\begin{array}
[c]{ccccc}%
\tfrac{1-\mu(2)}{2}=\tfrac{1-e^{-2}}{2}=m_{1}, &  & \tfrac{1-\mu(3)}%
{3}=0.2745509877..., &  & \tfrac{1-\mu(4)}{4}=0.2009733699...
\end{array}
\]
for $s=1$, $2$, $3$. \ On the one hand, the expected density $m_{2}%
=(2-e^{-1})/4$ for $2\times\infty$ and $s=1$ is verified in \cite{Seat3-err}.
\ On the other hand, the expected density for $3\times\infty$ is reported as
$\approx0.3915$ (via Monte Carlo simulation), inconsistent with $m_{3}=1/3$.
\ This issue awaits resolution. \ An interesting asymptotics problem appears
in \cite{Seat5-err}, as well as a constant $%
%TCIMACRO{\tsum \nolimits_{\ell=0}^{\infty}}%
%BeginExpansion
{\textstyle\sum\nolimits_{\ell=0}^{\infty}}
%EndExpansion
2^{-\ell(\ell+1)/2}=1.6416325606....$

Call an $n$-bit binary word \textbf{legal} if every $1$ has an adjacent $0$.
\ For example, if $n=6$, the only legal words with maximal set of $1$s are%
\[%
\begin{array}
[c]{ccccccccccccc}%
010101, &  & 010110, &  & 011001, &  & 011010, &  & 100110, &  & 101010, &  &
101101.
\end{array}
\]
Imagine cars ($1$s) parking one-by-one at random on $000000$, satisfying
legality at all times and stopping precisely when maximality is fulfilled.
\ This process endows the seven words with probabilities%
\[%
\begin{array}
[c]{ccccccccccccc}%
\frac{5}{48}, &  & \frac{7}{60}, &  & \frac{5}{48}, &  & \frac{7}{60}, &  &
\frac{5}{48}, &  & \frac{5}{48}, &  & \frac{7}{20}%
\end{array}
\]
respectively (by tree analysis) and the mean density of cars is%
\[
\tfrac{1}{6}\left[  3\left(  4\cdot\tfrac{5}{48}+2\cdot\tfrac{7}{60}\right)
+4\left(  \tfrac{7}{20}\right)  \right]  =\tfrac{67}{120}.
\]
In the limit as $n\rightarrow\infty$, the mean density $\rightarrow0.598...$
via simulation \cite{Rnyi1-err}. \ Conceivably this constant is exactly $3/5$,
but a proof may be difficult. \ Several variations on a discrete parking theme
appear in \cite{Rnyi1-err, Rnyi2-err}.

\textbf{5.4. Golomb-Dickman Constant.} Let $P^{+}(n)$ denote the largest prime
factor of $n$ and $P^{-}(n)$ denote the smallest prime factor of $n$. We
mentioned that
\[%
\begin{array}
[c]{ccc}%
%TCIMACRO{\dsum \limits_{n=2}^{N}}%
%BeginExpansion
{\displaystyle\sum\limits_{n=2}^{N}}
%EndExpansion
\ln(P^{+}(n))\sim\lambda N\ln(N)-\lambda(1-\gamma)N, &  &
%TCIMACRO{\dsum \limits_{n=2}^{N}}%
%BeginExpansion
{\displaystyle\sum\limits_{n=2}^{N}}
%EndExpansion
\ln(P^{-}(n))\sim e^{-\gamma}N\ln(\ln(N))+cN
\end{array}
\]
as $N\rightarrow\infty$, but did not give an expression for the constant $c$.
Tenenbaum \cite{Ten-err} found that
\[
c=e^{-\gamma}(1+\gamma)+%
%TCIMACRO{\dint \limits_{1}^{\infty}}%
%BeginExpansion
{\displaystyle\int\limits_{1}^{\infty}}
%EndExpansion
\frac{\omega(t)-e^{-\gamma}}{t}dt+%
%TCIMACRO{\dsum \limits_{p}}%
%BeginExpansion
{\displaystyle\sum\limits_{p}}
%EndExpansion
\left\{  e^{-\gamma}\ln\left(  1-\frac{1}{p}\right)  +\frac{\ln(p)}{p-1}%
%TCIMACRO{\dprod \limits_{q\leq p}}%
%BeginExpansion
{\displaystyle\prod\limits_{q\leq p}}
%EndExpansion
\left(  1-\frac{1}{q}\right)  \right\}  ,
\]
where the sum over $p$ and product over $q$ are restricted to primes. A
numerical evaluation is still open. Another integral \cite{FKnygnL-err}
\[%
%TCIMACRO{\dint \limits_{1}^{\infty}}%
%BeginExpansion
{\displaystyle\int\limits_{1}^{\infty}}
%EndExpansion
\frac{\rho(x)}{x}dx=(1.916045...)^{-1}%
\]
deserves closer attention (when the denominator is replaced by $x^{2}$,
$1-\lambda$ emerges). \ A variation of permutation, called cyclation, appears
in \cite{Pipp-err}. \ Similar constants arise in the distribution of cycle
lengths, given a random $n$-cyclation:%
\[%
\begin{array}
[c]{c}%
\text{expected }\\
\text{longest cycle}%
\end{array}
\sim\left(
%TCIMACRO{\dint \limits_{0}^{\infty}}%
%BeginExpansion
{\displaystyle\int\limits_{0}^{\infty}}
%EndExpansion
e^{-x+\operatorname*{Ei}(-x)/2}\,dx\right)  n=(0.7578230112...)n,
\]%
\[%
\begin{array}
[c]{c}%
\text{expected }\\
\text{shortest cycle}%
\end{array}
\sim\left(  \frac{\sqrt{\pi}}{2}%
%TCIMACRO{\dint \limits_{0}^{\infty}}%
%BeginExpansion
{\displaystyle\int\limits_{0}^{\infty}}
%EndExpansion
e^{-x-\operatorname*{Ei}(-x)/2}\,dx\right)  \sqrt{n}=(1.4572708792...)\sqrt{n}%
\]
as $n\rightarrow\infty$. The former coefficient is the Flajolet-Odlyzko
constant; the analogous growth rate of the latter for permutations is only
$\ln(n)$.

The longest tail $L(\varphi)$, given a random mapping $\varphi:\{1,2,\ldots
,n\}\rightarrow\{1,2,\ldots,n\}$, is called the \textit{height} of $\varphi$
in \cite{Pros-err, Kolc-err, AP1-err} and satisfies
\[
\lim_{n\rightarrow\infty}\operatorname*{P}\left(  \frac{L(\varphi)}{\sqrt{n}%
}\leq x\right)  =%
%TCIMACRO{\dsum \limits_{k=-\infty}^{\infty}}%
%BeginExpansion
{\displaystyle\sum\limits_{k=-\infty}^{\infty}}
%EndExpansion
(-1)^{k}\exp\left(  -\frac{k^{2}x^{2}}{2}\right)
\]
for fixed $x>0$. For example,
\[
\lim_{n\rightarrow\infty}\operatorname*{Var}\left(  \frac{L(\varphi)}{\sqrt
{n}}\right)  =\frac{\pi^{2}}{3}-2\pi\ln(2)^{2}.
\]
The longest rho-path $R(\varphi)$ is called the \textit{diameter} of $\varphi$
in \cite{AP2-err} and has moments
\[
\lim_{n\rightarrow\infty}\operatorname*{E}\left[  \left(  \frac{R(\varphi
)}{\sqrt{n}}\right)  ^{j}\right]  =\frac{\sqrt{\pi}j}{2^{j/2}\Gamma((j+1)/2)}%
%TCIMACRO{\dint \limits_{0}^{\infty}}%
%BeginExpansion
{\displaystyle\int\limits_{0}^{\infty}}
%EndExpansion
x^{j-1}(1-e^{\operatorname*{Ei}(-x)-I(x)})\,dx
\]
for fixed $j>0$. Complicated formulas for the distribution of the largest tree
$P(\varphi)$ also exist \cite{Kolc-err, AP1-err, Step-err}. \ From a series of
papers \cite{Fi39-err, Fi40-err, Fi41-err, Fi42-err, Fi43-err}, we quote just
one constant (of many):%
\[
\lim_{n\rightarrow\infty}\frac{\operatorname*{Var}(M(\varphi))}{n}%
=2H(1,1)+\left(  2-\frac{\pi}{2}\right)  G(1,1)^{2}=0.2411140734...,
\]
corresponding to the longest cycle length. \ Also, letting $\tilde{Q}%
(\varphi)$ denote the length of the (unique) cycle contained within the
largest component of $\varphi$, we have \cite{MutaFc-err}
\[
\lim_{n\rightarrow\infty}\frac{\operatorname*{E}(\tilde{Q}(\varphi))}{\sqrt
{n}}=\frac{1}{\sqrt{2}}%
%TCIMACRO{\dint \limits_{0}^{\infty}}%
%BeginExpansion
{\displaystyle\int\limits_{0}^{\infty}}
%EndExpansion
\frac{e^{-x+\operatorname*{Ei}(-x)/2}}{\sqrt{x}}\,dx=0.6884050874...;
\]
the analysis of such \textit{deepest} cycles is feasible, unlike the situation
for \textit{richest} components (containing the longest cycle of $\varphi$),
which remains unresolved.

A permutation $p\in S_{n}$ is an \textit{involution} if $p^{2}=1$ in $S_{n}$.
Equivalently, $p$ does not contain any cycles of length $>2$: it consists
entirely of fixed points and transpositions. Let $t_{n}$ denote the number of
involutions on $S_{n}$. Then $t_{n}=t_{n-1}+(n-1)t_{n-2}$ and \cite{CHM51-err,
WmpZ-err}
\[
t_{n}\sim\frac1{2^{1/2}e^{1/4}}\left(  \frac ne\right)  ^{n/2}e^{\sqrt{n}}%
\]
as $n\rightarrow\infty$. The equation $p^{d}=1$ for $d\geq3$ has also been
studied \cite{MW55-err}.

A permutation $p\in S_{n}$ is a \textit{square} if $p=q^{2}$ for some $q\in
S_{n}$; it is a \textit{cube} if $p=r^{3}$ for some $r\in S_{n}$. For
convenience, let $\omega=(-1+i\sqrt{3})/2$ and
\[
\Psi(x)=\frac13\left(  \exp(x)+2\exp(-x/2)\cos(\sqrt{3}x/2)\right)  .
\]
The probability that a random $n$-permutation is a square is \cite{Blum-err,
Bndr-err, Wilf-err, Pynn-err, FGDP-err}
\begin{align*}
&  \sim\frac{2^{1/2}}{\Gamma(1/2)}\frac1{n^{1/2}}%
%TCIMACRO{\dprod \limits_{1\leq m\equiv0\operatorname*{mod}2}}%
%BeginExpansion
{\displaystyle\prod\limits_{1\leq m\equiv0\operatorname*{mod}2}}
%EndExpansion
\frac{e^{1/m}+e^{-1/m}}2=\sqrt{\frac2{\pi\,n}}%
%TCIMACRO{\dprod \limits_{k=1}^{\infty}}%
%BeginExpansion
{\displaystyle\prod\limits_{k=1}^{\infty}}
%EndExpansion
\cosh\left(  \frac1{2k}\right) \\
&  =\sqrt{\frac2{\pi\,n}}(1.2217795151...)=(0.9748390118...)n^{-1/2}%
\end{align*}
as $n\rightarrow\infty$; the probability that it is a cube is \cite{Pynn-err,
FGDP-err}
\begin{align*}
&  \sim\frac{3^{1/3}}{\Gamma(2/3)}\frac1{n^{1/3}}%
%TCIMACRO{\dprod \limits_{1\leq m\equiv0\operatorname*{mod}3}}%
%BeginExpansion
{\displaystyle\prod\limits_{1\leq m\equiv0\operatorname*{mod}3}}
%EndExpansion
\frac{e^{1/m}+e^{\omega/m}+e^{\omega^{2}/m}}3\\
&  =\frac{3^{5/6}\Gamma(1/3)}{2\pi\,n^{1/3}}%
%TCIMACRO{\dprod \limits_{k=1}^{\infty}}%
%BeginExpansion
{\displaystyle\prod\limits_{k=1}^{\infty}}
%EndExpansion
\Psi\left(  \frac1{3k}\right)  =(1.0729979443...)n^{-1/3}.
\end{align*}

Two permutations $p,q\in S_{n}$ are of the \textit{same cycle type} if their
cycle decompositions are identical (in the sense that they possess the same
number of cycles of length $l$, for each $l\geq1$). The probability that two
independent, random $n$-permutations have the same cycle type is
\cite{FGDP-err}
\[
\sim\frac1{n^{2}}%
%TCIMACRO{\dprod \limits_{k=1}^{\infty}}%
%BeginExpansion
{\displaystyle\prod\limits_{k=1}^{\infty}}
%EndExpansion
I_{0}\left(  \frac2k\right)  =\left(  4.2634035141...\right)  n^{-2}%
\]
as $n\rightarrow\infty$, where $I_{0}$ is the zeroth modified Bessel function.

A mapping $\varphi$ on $\{1,2,\ldots,n\}$ has \textit{period} $\theta$ if
$\theta$ is the least positive integer for which iterates $\varphi^{m+\theta
}=\varphi^{m}$ for all sufficiently large $m$. It is known that
\cite{Shmtz1-err}
\[
\ln(\operatorname*{E}(\theta(\varphi)))=K\sqrt[3]{\frac{n}{\ln(n)^{2}}}\left(
1+o(1)\right)
\]
as $n\rightarrow\infty$, where $K=(3/2)(3\,b)^{2/3}=3.3607131721...$. A
typical mapping $\varphi$ satisfies $\ln(\theta(\varphi))\sim\frac{1}{8}%
\ln(n)^{2}$. When restricting the average to permutations $\pi$ only, we have
\[
\ln(\operatorname*{E}(\theta(\pi)))=B\sqrt{\frac{n}{\ln(n)}}\left(
1+o(1)\right)  ,
\]
where $B=2\sqrt{2b}=2.9904703993...$ (this corrects the error term on p. 287).
See \cite{DxPn-err, Shmtz2-err} for additional appearances of $B$. More on the
Erd\H{o}s-Tur\'{a}n constant is found in \cite{Bdsch-err, Lucht-err}.

Let $W(\pi)$ denote the number of factorizations of an $n$-permutation $\pi$
into two $n$-involutions. For example, if $\chi$ is an $n$-cycle, then
$W(\chi)=n$:%
\begin{align*}
(1\,2\,3\,4)  &  =(1\,2)(3\,4)\circ(1)(2\,4)(3)\\
&  =(1\,3)(2)(4)\circ(1\,2)(3\,4)\\
&  =(1\,4)(2\,3)\circ(1\,3)(2)(4)\\
&  =(1)(2\,4)(3)\circ(1\,4)(2\,3).
\end{align*}
If $\pi$ is chosen uniformly at random, then it is known that \cite{Lugo2-err}%
\[
\operatorname{E}\left(  W(\pi)\right)  \sim\frac{1}{\sqrt{8\pi e}}%
\frac{e^{2\sqrt{n}}}{\sqrt{n}}%
\]
as $n\rightarrow\infty$, and conjectured that%
\[
\lim_{n\rightarrow\infty}\operatorname{P}\left(  \frac{\ln(W(\pi))-\frac{1}%
{2}\ln(n)^{2}}{c\,\ln(n)^{3}}\leq x\right)  =\frac{1}{\sqrt{2\pi}}%
%TCIMACRO{\dint \limits_{-\infty}^{x}}%
%BeginExpansion
{\displaystyle\int\limits_{-\infty}^{x}}
%EndExpansion
\exp\left(  -\frac{t^{2}}{2}\right)  dt
\]
where $c\approx0.16$ is a constant.

\textbf{5.5. Kalm\'{a}r's Composition Constant.} See \cite{KlzL-err} for
precise inequalities involving $m(n)$ and $\rho=1.7286472389....$ The number
of factors in a random ordered factorization of $n\leq N$ into
$2,3,4,5,6,\ldots$ is asymptotically normal with mean \cite{Hwa2-err,
Hwa3-err}
\[
\sim\frac{-1}{\zeta^{\prime}(\rho)}\ln(N)=(0.5500100054...)\ln(N)
\]
and variance
\[
\sim\frac{-1}{\zeta^{\prime}(\rho)}\left(  \frac{\zeta^{\prime\prime}(\rho
)}{\zeta^{\prime}(\rho)^{2}}-1\right)  \ln(N)=(0.3084034446...)\ln(N)
\]
as $N\rightarrow\infty$. In contrast, the number of \textit{distinct} factors
in the same has mean
\[
\sim\frac{-1}{\rho}\Gamma\left(  \frac{-1}{\rho}\right)  \left(  \frac
{-1}{\zeta^{\prime}(\rho)}\right)  ^{1/\rho}\ln(N)^{1/\rho}%
=(1.4879159716...)\ln(N)^{1/\rho},
\]
hence on average there are many small factors occurring with high frequencies.
Also, the number of factors in a random ordered factorization of $n\leq N$
into $2,3,5,7,11,\ldots$ is asymptotically normal with mean $0.5776486251...$
and variance $0.4843965045...$ (with $\eta=1.3994333287...$ and $%
%TCIMACRO{\tsum \nolimits_{p}}%
%BeginExpansion
{\textstyle\sum\nolimits_{p}}
%EndExpansion
p^{-s}$ playing the roles of $\rho$ and $\zeta(s)-1$).

A\ \textit{Carlitz composition} of size $n$ is an additive composition
$n=x_{1}+x_{2}+\cdots+x_{k}$ such that $x_{j}\neq x_{j+1}$ for any $1\leq j<k
$. We call $k$ the \textit{number of parts} and
\[
d=1+%
%TCIMACRO{\dsum \limits_{i=2}^{k}}%
%BeginExpansion
{\displaystyle\sum\limits_{i=2}^{k}}
%EndExpansion
\left\{
\begin{array}
[c]{lll}%
1 &  & \text{if }x_{i}\neq x_{j}\text{ for all }1\leq j<i,\\
0 &  & \text{otherwise}%
\end{array}
\right.
\]
the \textit{number of distinct part sizes}. The number $a_{\text{c}}(n)$ of
Carlitz compositions is \cite{Crltz0-err, Crltz1-err, Ktsvc-err, FljSdg-err}
\[
a_{\text{c}}(n)\sim\frac{1}{\sigma\,F^{\prime}(\sigma)}\left(  \frac{1}%
{\sigma}\right)  ^{n}=(0.4563634740...)(1.7502412917...)^{n}%
\]
where $\sigma=0.5713497931...$ is the unique solution of the equation
\[%
\begin{array}
[c]{ccc}%
F(x)=%
%TCIMACRO{\dsum \limits_{j=1}^{\infty}}%
%BeginExpansion
{\displaystyle\sum\limits_{j=1}^{\infty}}
%EndExpansion
(-1)^{j-1}\dfrac{x^{j}}{1-x^{j}}=1, &  & 0\leq x\leq1.
\end{array}
\]
The expected number of parts is asymptotically
\[%
\begin{array}
[c]{ccccc}%
\dfrac{G(\sigma)}{\sigma\,F^{\prime}(\sigma)}n\sim(0.350571...)n &  &
\text{where} &  & G(x)=%
%TCIMACRO{\dsum \limits_{j=1}^{\infty}}%
%BeginExpansion
{\displaystyle\sum\limits_{j=1}^{\infty}}
%EndExpansion
(-1)^{j-1}\dfrac{j\,x^{j}}{1-x^{j}}%
\end{array}
\]
(by contrast, an unrestricted composition has $(n+1)/2$ parts on average). The
expected size of the largest part is
\[
\frac{-\ln(n)}{\ln(\sigma)}+\left(  \frac{\ln(F^{\prime}(\sigma))+\ln
(1-\sigma)-\gamma}{\ln(\sigma)}+\frac{1}{2}\right)  +\varepsilon
(n)=(1.786500...)\ln(n)+0.643117...+\varepsilon(n)
\]
where $\gamma$ is Euler's constant and $\varepsilon(n)$ is a small-amplitude
zero-mean periodic function. The expected number of distinct part sizes is
\cite{Crltz2-err}
\[
\frac{-\ln(n)}{\ln(\sigma)}+\left(  \frac{\ln(F^{\prime}(\sigma))+\gamma}%
{\ln(\sigma)}+\frac{1}{2}\right)  +\delta(n)=(1.786500...)\ln
(n)-2.932545...+\delta(n)
\]
where $\delta(n)$ is likewise negligible. (By contrast, an unrestricted
composition has a largest part of size roughly $\ln(n)/\ln(2)+0.332746...$ and
roughly $\ln(n)/\ln(2)-0.667253...$ distinct part sizes on average:\ see
\cite{Crltz3-err, Crltz4-err, Crltz5-err}, as well as the bottom of page 340.)
We wonder about the multiplicative analog of these results. See also
\cite{Crltz6-err}.

Another equation involving the Riemann zeta function: \cite{AiCLu-err}
\[
\zeta(x-2)-2\zeta(x-1)=0
\]
arises in random graph theory and its solution $x=3.4787507857...$ serves to
separate one kind of qualitative behavior (the existence of a giant component)
from another. \ The same expression (with $x$ replaced by $y+1$) appears in%
\[%
%TCIMACRO{\dsum \limits_{k=0}^{\infty}}%
%BeginExpansion
{\displaystyle\sum\limits_{k=0}^{\infty}}
%EndExpansion
\frac{1}{r(k)^{y}}=\frac{\zeta(y-1)}{2\zeta(y)-\zeta(y-1)},
\]
where $r(k)$ is the number of representations of an integer $k$ as a sum of
distinct Fibonacci numbers \cite{ZQBF1-err}. \ A conjectured limit involving
the number $c(n)$ of primitive subsets of $\{1,2,\ldots,n\}$ indeed exists
\cite{RAngl-err} but its precise value remains open.

\textbf{5.6. Otter's Tree Enumeration Constants. }Higher-order asymptotic
series for $T_{n}$, $t_{n}$ and $B_{n}$ are given in \cite{Fi1-err}. Analysis
of series-parallel posets \cite{Fi6-err} is similar to that of trees. By
Stirling's formula, another way of writing the asymptotics for labeled mobiles
is \cite{FljSdg-err}%
\[
\frac{\hat{M}_{n}}{n!}\sim\frac{\hat{\eta}}{\sqrt{2\pi}}\left(  e\,\hat{\xi
}\right)  ^{n}n^{-3/2}\sim(0.1857629435...)\left(  3.1461932206...\right)
^{n}n^{-3/2}%
\]
as $n\rightarrow\infty$. \ See \cite{FGL-err, LLL-err} for more about
$k$-gonal $2$-trees, as well as a new formula for $\alpha$ in terms of
rational expressions involving $e$.

The generating function $L(x)$ of leftist trees satisfies a simpler functional
equation than previously thought:
\[
L(x)=x+L\left(  x\,L(x)\right)
\]
which involves an unusual nested construction. The radius of convergence
$\rho=0.3637040915...=(2.7494879027...)^{-1}$ of $L(x)$ satisfies
\[
\rho\,L^{\prime}\left(  \rho\,L(\rho)\right)  =1
\]
and the coefficient of $\rho^{-n}n^{-3/2}$ in the asymptotic expression for
$L_{n}$ is
\[
\sqrt{\frac1{2\pi\rho^{2}}\frac{\rho+L(\rho)}{L^{\prime\prime}\left(
\rho\,L(\rho)\right)  }}=0.2503634293...=(0.6883712204...)\rho.
\]
The average height of $n$-leaf leftist trees is asymptotically
$(1.81349371...)\sqrt{\pi n}$ and the average depth of vertices belonging to
such trees is asymptotically $(0.90674685...)\sqrt{\pi n}$. Nogueira
\cite{Nogu-err} conjectured that the ratio of the two coefficients is exactly
$2$, but his only evidence is numerical (to over $1000$ decimal digits). Let
the $d$-number of an ordered binary tree $\tau$ be
\[
d(\tau)=\left\{
\begin{array}
[c]{lll}%
1 &  & \text{if }\tau_{L}=\emptyset\text{ or }\tau_{R}=\emptyset,\\
1+\min(d(\tau_{L}),d(\tau_{R})) &  & \text{otherwise.}%
\end{array}
\right.
\]
Such a tree is leftist if and only if for every subtree $\sigma$ of $\tau$
with $\sigma_{L}\neq\emptyset$ and $\sigma_{R}\neq\emptyset$, the inequality
$d(\sigma_{L})>d(\sigma_{R})$ holds. Another relevant constant,
$0.6216070079...$, is involved in a distribution law for leftist trees in
terms of their $d$-number \cite{Nogu-err}.

For the following, we consider only unordered forests whose connected
components are (strongly) ordered binary trees. Let $F_{n}$ denote the number
of such forests with $2n-1$ vertices; then the generating function
\[
\Phi(x)=1+%
%TCIMACRO{\dsum \limits_{n=1}^{\infty}}%
%BeginExpansion
{\displaystyle\sum\limits_{n=1}^{\infty}}
%EndExpansion
F_{n}x^{n}=1+x+2x^{2}+4x^{3}+10x^{4}+26x^{5}+77x^{6}+\cdots
\]
satisfies
\[
\Phi(x)=\exp\left(
%TCIMACRO{\dsum \limits_{k=1}^{\infty}}%
%BeginExpansion
{\displaystyle\sum\limits_{k=1}^{\infty}}
%EndExpansion
\frac{1-\sqrt{1-4x^{k}}}{2k}\right)  =%
%TCIMACRO{\dprod \limits_{m=1}^{\infty}}%
%BeginExpansion
{\displaystyle\prod\limits_{m=1}^{\infty}}
%EndExpansion
(1-x^{m})^{-\frac{1}{m}\tbinom{2m-2}{m-1}}.
\]
It can be shown that \cite{FGDP-err}
\[
F_{n}\sim\frac{\Phi(1/4)}{\sqrt{\pi}}\frac{4^{n-1}}{n^{3/2}}=\frac
{1.7160305349...}{4\sqrt{\pi}}\frac{4^{n}}{n^{3/2}}%
\]
as $n\rightarrow\infty$. The constant $1.716...$ also plays a role in the
asymptotic analysis of the probability that a random forest has no two
components of the same size.

A phylogenetic tree of size $n$ is a strongly binary tree whose $n$ leaves are
labeled. \ The number of such trees is $1\cdot3\cdots(2n-3)$ and two such
trees are isomorphic if removing their labels will associate them to the same
unlabeled tree. \ The probability that two uniformly-selected phylogenetic
trees are isomorphic is asymptotically \cite{BonFlj-err}%
\[
(3.17508...)(2.35967...)^{-n}n^{3/2}%
\]
as $n\rightarrow\infty$, where the growth rate is $4\rho$ and $\rho
=0.5899182714...$ is the radius of convergence of a certain radical expansion%
\[
1-\sqrt{\frac{3}{2}-2z-\frac{1}{2}\sqrt{\frac{15}{8}-2z^{2}-\frac{7}{8}%
\sqrt{\frac{255}{128}-2z^{4}-\frac{127}{128}\sqrt{\ldots}}}}.
\]

An arithmetic formula is an expression involving only the number $1$ and
operations $+$ and $\cdot$, with multiplication by $1$ disallowed. For
example, $4$ has exactly six arithmetic formulas:%
\[%
\begin{array}
[c]{ccccc}%
1+(1+(1+1)), &  & 1+((1+1)+1), &  & (1+(1+1))+1,\\
((1+1)+1)+1, &  & (1+1)+(1+1), &  & (1+1)\cdot(1+1).
\end{array}
\]
Let $f(n)$ denote the number of arithmetic formulas for $n$ and $F(x)=%
%TCIMACRO{\tsum \nolimits_{n=1}^{\infty}}%
%BeginExpansion
{\textstyle\sum\nolimits_{n=1}^{\infty}}
%EndExpansion
f(n)x^{n}$, then define $\xi$ to be the smallest positive solution of the
equation
\[
\frac{1}{4}=x+%
%TCIMACRO{\dsum \limits_{k=2}^{\infty}}%
%BeginExpansion
{\displaystyle\sum\limits_{k=2}^{\infty}}
%EndExpansion
f(k)\left(  F\left(  x^{k}\right)  -x^{k}\right)
\]
and $\eta=1/\xi$ to be the growth rate. \ A binary tree-like argument yields
that $f(n)$ is asymptotically \cite{Gnang1-err, Gnang2-err}%
\[
(0.1456918546...)(4.0765617852...)^{n}n^{-3/2}%
\]
as $n\rightarrow\infty$. \ Suppose moreover that exponentiation is included
but that $1$ again is disallowed; thus $(1+1)^{(1+1)}$ also counts. An analog
holds for counting arithmetic exponential formulas but with a larger
$\eta=4.1307352951...$.

\textbf{5.7. Lengyel's Constant.} Constants of the form $%
%TCIMACRO{\tsum \nolimits_{k=-\infty}^{\infty}}%
%BeginExpansion
{\textstyle\sum\nolimits_{k=-\infty}^{\infty}}
%EndExpansion
2^{-k^{2}}$ and $%
%TCIMACRO{\tsum \nolimits_{k=-\infty}^{\infty}}%
%BeginExpansion
{\textstyle\sum\nolimits_{k=-\infty}^{\infty}}
%EndExpansion
2^{-(k-1/2)^{2}}$ appear in \cite{Fi7-err, Fi8-err}. We discussed the
refinement of $B_{n}$ given by $S_{n,k}$, which counts partitions of
$\{1,2,\ldots,n\}$ possessing exactly $k$ blocks. Another refinement of
$B_{n}$ is based jointly on the maximal $i$ such that a partition has an
$i$-crossing and the maximal $j$ such that the partition has a $j$-nesting
\cite{CDSY-err}. The cardinality of partitions avoiding $2$-crossings is the
$n^{\text{th}}$ Catalan number; see \cite{BsqX-err} for partitions avoiding
$3$-crossings and \cite{JQiR-err} for what are called $3$-noncrossing
\textit{braids}.

\textbf{5.8. Takeuchi-Prellberg Constant}. Knuth's recursive formula should be
replaced by
\[
T_{n+1}=%
%TCIMACRO{\dsum \limits_{k=0}^{n-1}}%
%BeginExpansion
{\displaystyle\sum\limits_{k=0}^{n-1}}
%EndExpansion
\left[  2\tbinom{n+k}k-\tbinom{n+k+1}k\right]  T_{n-k}+%
%TCIMACRO{\dsum \limits_{k=1}^{n+1}}%
%BeginExpansion
{\displaystyle\sum\limits_{k=1}^{n+1}}
%EndExpansion
\tbinom{2k}k\frac1{k+1}.
\]

\textbf{5.9. P\'olya's Random Walk Constants.} Properties of the gamma
function lead to a further simplification \cite{JZ-err}:
\[
m_{3}=\frac1{32\pi^{3}}\left(  \sqrt{3}-1\right)  \left[  \Gamma\left(
\frac1{24}\right)  \Gamma\left(  \frac{11}{24}\right)  \right]  ^{2}%
\]
Consider a variation in which the drunkard performs a random walk starting
from the origin with $2^{d}$ equally probable steps, each of the form
$(\pm1,\pm1,\ldots,\pm1)$. The number of walks that end at the origin after
$2n $ steps is
\[
\tilde U_{d,0,2n}=\dbinom{2n}n^{d}%
\]
and the number of such walks for which $2n$ is the time of \textit{first
return} to the origin is $\tilde V_{d,0,2n}$, where \cite{FFK-err}
\[
2^{-n}\tilde V_{1,0,2n}=\frac1{n2^{2n-1}}\dbinom{2n-2}{n-1}\sim\frac
1{2\sqrt{\pi}n^{3/2}},
\]
\[
2^{-2n}\tilde V_{2,0,2n}=\frac\pi{n(\ln(n))^{2}}-2\pi\frac{\gamma+\pi B}%
{n(\ln(n))^{3}}+O\left(  \frac1{n(\ln(n))^{4}}\right)  ,
\]
\[
2^{-3n}\tilde V_{3,0,2n}=\frac1{\pi^{3/2}C^{2}n^{3/2}}+O\left(  \frac1{n^{2}%
}\right)
\]
as $n\rightarrow\infty$, where
\[
B=1+%
%TCIMACRO{\dsum \limits_{k=1}^{\infty}}%
%BeginExpansion
{\displaystyle\sum\limits_{k=1}^{\infty}}
%EndExpansion
\left[  2^{-4k}\dbinom{2k}k^{2}-\frac1{\pi k}\right]  =\frac{4\ln(2)}%
\pi=0.8825424006...,
\]
\[
C=%
%TCIMACRO{\dsum \limits_{k=0}^{\infty}}%
%BeginExpansion
{\displaystyle\sum\limits_{k=0}^{\infty}}
%EndExpansion
2^{-6k}\dbinom{2k}k^{3}=\frac1{4\pi^{3}}\Gamma\left(  \frac14\right)
^{4}=1.3932039296....
\]

The quantity $W_{d,n}$ is often called the \textit{average range} of the
random walk (equal to $\operatorname*{E}(\max\omega_{j}-\min\omega_{j})$ when
$d=1$). The corresponding variance is
\[
\sim4\left(  \ln(2)-\frac2\pi\right)  n=(0.2261096327...)n
\]
if $d=1$ \cite{Fllr-err} and is
\[
\sim8\pi^{2}\left(  \frac32L_{-3}(2)+\frac12-\frac{\pi^{2}}{12}\right)
\frac{n^{2}}{\ln(n)^{4}}=8\pi^{2}\left(  0.8494865859...\right)  \frac{n^{2}%
}{\ln(n)^{4}}%
\]
if $d=2$ \cite{JnPt-err}. Various representations include
\[
\frac32L_{-3}(2)=1.1719536193...=-%
%TCIMACRO{\dint \limits_{0}^{1}}%
%BeginExpansion
{\displaystyle\int\limits_{0}^{1}}
%EndExpansion
\frac{\ln(x)}{1-x+x^{2}}dx=\frac2{\sqrt{3}}(1.0149416064...),
\]
the latter being Lobachevsky's constant (p. 233). Exact formulas for the
corresponding distribution, for any $n$, are available when $d=1$
\cite{Vlls-err}.

More on the constant $\rho$ appears in \cite{Wjsmn-err, Tuentr-err}. It turns
out that the constant $\sigma$, given by an infinite series, has a more
compact integral expression \cite{CtMj-err, MjCZ-err}:
\[
\sigma=\frac1\pi%
%TCIMACRO{\dint \limits_{0}^{\infty}}%
%BeginExpansion
{\displaystyle\int\limits_{0}^{\infty}}
%EndExpansion
\frac1{x^{2}}\ln\left[  \frac6{x^{2}}\left(  1-\frac{\sin(x)}x\right)
\right]  dx=-0.2979521902...=\frac{-0.5160683318...}{\sqrt{3}}%
\]
and surprisingly appears in both 3D statistical mechanics \cite{Ziff-err} and
1D probabilistic algorithmics \cite{CFFH-err}.

Here is a problem about stopping times for certain one-dimensional walks. Fix
a large integer $n$. At time $0$, start with a total of $n+1$ particles, one
at each integer site in $[0,n]$. At each positive integer time, randomly
choose one of the particles remaining in $[1,n]$ and move it $1$ step to the
left, coalescing with any particle that might already occupy the site. Let
$T_{n}$ denote the time at which only one particle is left (at $0$). An exact
expression for the mean of $T_{n}$ is known \cite{LrsLy-err}:
\[
\operatorname*{E}(T_{n})=\frac{2n(2n+1)}3\binom{2n}n\frac1{2^{2n}}\sim
\frac4{3\sqrt{\pi}}n^{3/2}=(0.7522527780...)n^{3/2}%
\]
and the variance is conjectured to satisfy
\[%
\begin{array}
[c]{ccc}%
\operatorname*{Var}(T_{n})\sim C\,n^{5/2}, &  & 0<C\leq\dfrac8{15\sqrt{\pi}%
}<0.301.
\end{array}
\]
Simulation suggests that $C\sim0.026$ and that a Central Limit Theorem holds
\cite{PDykl-err}.

\textbf{5.10. Self-Avoiding Walk Constants.} A conjecture due to Jensen
\&\ Guttmann \cite{JnsGtm-err}%
\[
\mu=\sqrt{\frac{7+\sqrt{30261}}{26}}%
\]
for the square lattice seems completely unmotivated yet numerically
reasonable; in contrast, a proposal%
\[
\mu=\sqrt{2+\sqrt{2}}%
\]
for the hexgonal lattice is now a theorem \cite{DCSmv-err, Grmmtt-err}. \ If
we examine SAPs rather than SAWs, it seems that $\gamma=-3/2$ and
$A=0.56230129...$\cite{Jens2-err, Jens3-err}. \ Fascinating complications
arise if such are restricted to be \textit{prudent}, that is, never take a
step towards an already occupied vertex \cite{BtnFG-err}.

Hueter \cite{Hu1-err, Hu2-err} claimed a proof that $\nu_{2}=3/4$ and that
$7/12\leq\nu_{3}\leq2/3$, $1/2\leq v_{4}\leq5/8$ (if the mean square
end-to-end distance exponents $\nu_{3}$, $v_{4}$ exist; otherwise the bounds
apply for
\[%
\begin{array}
[c]{ccc}%
\underline{\nu}_{d}=\operatorname*{liminf}\limits_{n\rightarrow\infty}%
\dfrac{\ln(r_{n})}{2\ln(n)}, &  & \overline{\nu}_{d}=\operatorname*{limsup}%
\limits_{n\rightarrow\infty}\dfrac{\ln(r_{n})}{2\ln(n)}%
\end{array}
\]
when $d=3,4$). She confirmed that the same exponents apply for the mean square
radius of gyration $s_{n}$ for $d=2,3,4$; the results carry over to
self-avoiding trails as well. Burkhardt \&\ Guim \cite{BuGu-err} adjusted the
estimate for $\lim_{k\rightarrow\infty}p_{k,k}^{1/k^{2}}$ to $1.743...$; this
has now further been improved to $1.74455...$ \cite{BMGI-err}.

\textbf{5.11. Feller's Coin Tossing Constants}. \ The cubic irrational
$1.7548776662...$ turns out to be the square of the Plastic constant $\psi$
and has infinite radical expression%
\[
\psi^{2}=1+\frac{1}{\sqrt{1+\frac{1}{\sqrt{1+\frac{1}{\sqrt{1+\cdots}}}}}%
}=1+\dfrac{\left.  1\right\vert }{\sqrt{1}}+\dfrac{\left.  1\right\vert
}{\sqrt{1}}+\dfrac{\left.  1\right\vert }{\sqrt{1}}+\cdots,
\]
an observation due to Knuth \cite{Grdnr-err}. \ Additional references on
oscillatory phenomena in probability theory include \cite{FGH1-err, FGH2-err,
FGH3-err, FGH4-err, FGH5-err}; see also our earlier entry [5.5]. Consider $n$
independent non-homogeneous Bernoulli random variables $X_{j}$ with
$\operatorname*{P}(X_{j}=1)=p_{j}=\operatorname*{P}($heads$)$ and
$\operatorname*{P}(X_{j}=0)=1-p_{j}=\operatorname*{P}($tails$)$. If all
probabilities $p_{j}$ are equal, then
\[
\sqrt{%
%TCIMACRO{\tsum \limits_{j=1}^{n}}%
%BeginExpansion
{\textstyle\sum\limits_{j=1}^{n}}
%EndExpansion
p_{j}(1-p_{j})}\operatorname*{P}(X_{1}+X_{2}+\cdots+X_{n}=k)\leq\frac{1}%
{\sqrt{2e}}=0.4288819424...
\]
for all integers $k$ and the bound is sharp. If there exist at least two
distinct values $p_{i}$, $p_{j}$, then \cite{CmVs-err}
\[
\sqrt{%
%TCIMACRO{\tsum \limits_{j=1}^{n}}%
%BeginExpansion
{\textstyle\sum\limits_{j=1}^{n}}
%EndExpansion
p_{j}(1-p_{j})}\operatorname*{P}(X_{1}+X_{2}+\cdots+X_{n}=k)\leq
M=0.4688223554...
\]
for all integers $k$ and the bound is sharp, where
\[
M=\max_{u\geq0}\sqrt{2u}e^{-2u}%
%TCIMACRO{\dsum \limits_{\ell=0}^{\infty}}%
%BeginExpansion
{\displaystyle\sum\limits_{\ell=0}^{\infty}}
%EndExpansion
\left(  \frac{u^{\ell}}{\ell!}\right)  ^{2}%
\]
and the maximizing argument is $u=0.3949889297...$.

\textbf{5.12. Hard Square Entropy Constant.} McKay \cite{McK-err} observed the
following asymptotic behavior:
\[
F(n)\sim(1.06608266...)(1.0693545387...)^{2n}(1.5030480824...)^{n^{2}}%
\]
based on an analysis of the terms $F(n)\,$ up to $n=19$. He emphasized that
the form of right hand side is conjectural, even though the data showed quite
strong convergence to this form. \ Counting \textit{maximal} independent
vertex subsets of the $n\times n$ grid graph is more difficult
\cite{Euler-err}: we have $1$, $2$, $10$, $42$, $358$ for $1\leq n\leq5$ but
nothing yet for $n\geq6$. \ By \textquotedblleft maximal\textquotedblright, we
mean with respect to set-inclusion. \ There is a natural connection with
discrete parking (see section 5.3.1). \ Asymptotics remain open here.

To calculate entropy constants of more complicated planar examples, such as
the 4-8-8 and triangular Kagom\'{e} lattices, requires more intricate
analysis. \ The former has numerical value
$1.54956010...=(5.76545652...)^{1/4}$; the latter evidently still remains open
\cite{ZhgMS-err}. \ A nonplanar example is the square lattice with crossed
diagonal bonds, which has entropy constant between $1.34254$ and $1.34265$.

Let $L(m,n)$ denote the number of legal positions on an $m\times n$ Go board
(a popular game). Then \cite{JTrmp1-err}
\[
\lim_{n\rightarrow\infty}L(1,n)^{1/n}=1+\frac{1}{3}\left(  \left(
27+3\sqrt{57}\right)  ^{1/3}+\left(  27-3\sqrt{57}\right)  ^{1/3}\right)
=2.7692923542...,
\]%
\[
\lim_{n\rightarrow\infty}L(n,n)^{1/n^{2}}=2.9757341920...
\]
and, subject to a plausible conjecture,
\[
L(m,n)\sim(0.8506399258...)(0.96553505933...)^{m+n}(2.9757341920...)^{mn}%
\]
as $\min\{m,n\}\rightarrow\infty$.

The number $Q(n)$ of configurations of $n$ non-attacking Queens on an $n\times
n$ chessboard satisfies%
\[
\lim_{n\rightarrow\infty}\frac{Q(n)^{1/n}}{n}=e^{-\alpha}%
\]
where $\alpha=1.94400...$ is called the $n$-queens constant \cite{Que1-err,
Que2-err, Que3-err}. \ This estimate emerges via numerical solution of a
difficult convex optimization problem. \ Counting configurations of $n$
non-attacking Bishops or of $n$ non-attacking Rooks is simpler \cite{Que4-err,
Que5-err}. \ A Queen, of course, possesses precisely the capabilities of both
a Bishop and a Rook.

\textbf{5.13. Binary Search Tree Constants.} The variance for the number of
comparisons in a successful search (odd $x$) should be%
\[
\operatorname*{Var}(f(x,V))=\left(  2+\frac{10}{n}\right)  H_{n}-4\left(
1+\frac{1}{n}\right)  \left(  H_{n}^{(2)}+\frac{H_{n}^{2}}{n}\right)  +4,
\]
that is, the denominator for $H_{n}^{2}$ is not $4$ but rather $n$. In the
subsequent two asymptotic expressions, $\pi^{3}$ should be replaced by
$\pi^{2}$.\ Also, the total (internal) path length satisfies
\begin{align*}
\operatorname*{E}\left[
%TCIMACRO{\dsum \limits_{\text{odd }x=1}^{2n-1}}%
%BeginExpansion
{\displaystyle\sum\limits_{\text{odd }x=1}^{2n-1}}
%EndExpansion
\left(  f(x,V)-1\right)  \right]   &  =2(n+1)H_{n}-4n\\
&  =2n\ln(n)+2(\gamma-2)n+2\ln(n)+(2\gamma+1)+o(1),
\end{align*}%
\begin{align*}
\operatorname*{Var}\left[
%TCIMACRO{\dsum \limits_{\text{odd }x=1}^{2n-1}}%
%BeginExpansion
{\displaystyle\sum\limits_{\text{odd }x=1}^{2n-1}}
%EndExpansion
\left(  f(x,V)-1\right)  \right]   &  =7n^{2}-4(n+1)^{2}H_{n}^{(2)}%
-2(n+1)H_{n}+13n\\
&  =\left(  7-\frac{2\pi^{2}}{3}\right)  n^{2}-2n\ln(n)+\left(  17-2\gamma
-\frac{4\pi^{2}}{3}\right)  n\\
&  -2\ln(n)+\left(  5-2\gamma-\frac{2\pi^{2}}{3}\right)  +o(1)
\end{align*}
as $n\rightarrow\infty$.

The random permutation model for generating weakly binary trees (given an
$n$-vector of distinct integers, construct $T$ via insertions) does
\textit{not} provide equal weighting on the $\tbinom{2n}{n}/(n+1)$ possible
trees. For example, when $n=3$, the permutations $(2,1,3)$ and $(2,3,1)$ both
give rise to the same tree $S$, which hence has probability $q(S)=1/3$ whereas
$q(T)=1/6$ for the other four trees. Fill \cite{FFK-err, Fil1-err, Fil2-err}
asked how the numbers $q(T)$ themselves are distributed, for fixed $n$. If the
trees are endowed with the uniform distribution, then
\begin{align*}
\frac{-\operatorname*{E}\left[  \ln(q(T))\right]  }{n}  &  \rightarrow%
%TCIMACRO{\dsum \limits_{k=1}^{\infty}}%
%BeginExpansion
{\displaystyle\sum\limits_{k=1}^{\infty}}
%EndExpansion
\frac{\ln(k)}{(k+1)4^{k}}\dbinom{2k}{k}\\
\  &  =-\gamma-%
%TCIMACRO{\dint \limits_{0}^{1}}%
%BeginExpansion
{\displaystyle\int\limits_{0}^{1}}
%EndExpansion
\frac{\ln(\ln(1/t))}{\sqrt{1-t}\left(  1+\sqrt{1-t}\right)  ^{2}%
}dt=2.0254384677...
\end{align*}
as $n\rightarrow\infty$. If, instead, the trees follow the distribution $q$,
then
\begin{align*}
\frac{-\operatorname*{E}\left[  \ln(q(T))\right]  }{n}  &  \rightarrow2%
%TCIMACRO{\dsum \limits_{k=1}^{\infty}}%
%BeginExpansion
{\displaystyle\sum\limits_{k=1}^{\infty}}
%EndExpansion
\frac{\ln(k)}{(k+1)(k+2)}\\
\  &  =-\gamma-2%
%TCIMACRO{\dint \limits_{0}^{1}}%
%BeginExpansion
{\displaystyle\int\limits_{0}^{1}}
%EndExpansion
\frac{\left(  (t-2)\ln(1-t)-2t\right)  \ln(\ln(1/t))}{t^{3}}%
dt=1.2035649167....
\end{align*}
The maximum value of $-\ln(q(T))$ is $\sim n\ln(n)$ and the minimum value is
$\sim c\,n$, where
\[
c=\ln(4)+%
%TCIMACRO{\dsum \limits_{k=1}^{\infty}}%
%BeginExpansion
{\displaystyle\sum\limits_{k=1}^{\infty}}
%EndExpansion
2^{-k}\ln(1-2^{-k})=0.9457553021....
\]
See also \cite{HtrOch-err, ShppZZ-err} for more on random sequential bisections.

\textbf{5.14. Digital Search Tree Constants.} Letting $Q_{\ell}$ denote the
$\ell^{\text{th}}$ partial product of $Q$ and%
\[
g(x)=\left\{
\begin{array}
[c]{lll}%
\dfrac{x-\ln(x)-1}{(x-1)^{2}} &  & \text{if }x\neq1,\\
\dfrac{1}{2} &  & \text{if }x=1,
\end{array}
\right.
\]
the total (internal) path length satisfies \cite{KPSz-err, HwFZ-err}%
\begin{align*}
\operatorname*{E}\left[
%TCIMACRO{\dsum \limits_{x=m_{i}}}%
%BeginExpansion
{\displaystyle\sum\limits_{x=m_{i}}}
%EndExpansion
\left(  f(x,M,1)-1\right)  \right]   &  =\frac{n\ln(n)}{\ln(2)}+n\left(
\frac{\gamma-1}{\ln(2)}+\frac{1}{2}-\alpha+\delta_{1}\left(  n\right)
\right)  +\frac{\ln(n)}{\ln(2)}\\
&  +\left(  \frac{2\gamma-1}{2\ln(2)}+\frac{5}{2}-\alpha\right)  +\delta
_{2}\left(  n\right)  +O\left(  \frac{\ln(n)}{n}\right)  ,
\end{align*}%
\[
\operatorname*{Var}\left[
%TCIMACRO{\dsum \limits_{x=m_{i}}}%
%BeginExpansion
{\displaystyle\sum\limits_{x=m_{i}}}
%EndExpansion
\left(  f(x,M,1)-1\right)  \right]  =n\left(  C+\delta_{3}\left(  n\right)
\right)  +O\left(  \frac{\ln(n)^{2}}{n}\right)
\]
where the sum is taken over all rows of $M$ and%
\[
C=\frac{Q}{\ln(2)}%
%TCIMACRO{\dsum \limits_{j,k,\ell\geq0}}%
%BeginExpansion
{\displaystyle\sum\limits_{j,k,\ell\geq0}}
%EndExpansion
\frac{(-1)^{j}}{Q_{j}Q_{k}Q_{\ell}}2^{-j(j+1)/2-k-\ell}g\left(  2^{-j-k}%
+2^{-j-\ell}\right)  =0.2660036454....
\]

Erd\H{o}s' 1948 irrationality proof is discussed in \cite{Vndhy-err}. The
constant $Q$ is transcendental via a general theorem on values of modular
forms due to Nesterenko \cite{Nest-err, Zud2-err}. A correct formula for
$\theta$ is
\[
\theta=%
%TCIMACRO{\dsum \limits_{k=1}^{\infty}}%
%BeginExpansion
{\displaystyle\sum\limits_{k=1}^{\infty}}
%EndExpansion
\frac{k2^{k(k-1)/2}}{1\cdot3\cdot7\cdots(2^{k}-1)}%
%TCIMACRO{\dsum \limits_{j=1}^{k}}%
%BeginExpansion
{\displaystyle\sum\limits_{j=1}^{k}}
%EndExpansion
\frac{1}{2^{j}-1}=7.7431319855...
\]
(the exponent $k(k-1)/2$ was mistakenly given as $k+1$ in \cite{FS-err}, but
the numerical value is correct). The constants $\alpha$, $\beta$ and $Q^{-1}$
appear in \cite{LiZw-err}. \ Also, $\alpha$ appears in \cite{KuWa-err},
$Q^{-1}$ in \cite{Fi8-err} and%
\[%
%TCIMACRO{\dprod \limits_{n=1}^{\infty}}%
%BeginExpansion
{\displaystyle\prod\limits_{n=1}^{\infty}}
%EndExpansion
\left(  1-\frac{1}{2^{n/2}}\right)  =0.0375130167...
\]
in \cite{FuQiZ1-err, FuQiZ2-err, FuQiZ3-err}. \ What is the variance for the
count of internal nodes of a random trie under a Poisson($n$) model? \ This
turns out to be asymptotically linear in $n$ with slope \cite{RgJcq-err}%
\[
\frac{1}{\ln(2)}\left(  -\frac{1}{2}+2%
%TCIMACRO{\dsum \limits_{\ell=0}^{\infty}}%
%BeginExpansion
{\displaystyle\sum\limits_{\ell=0}^{\infty}}
%EndExpansion
\frac{1}{2^{\ell}+1}\right)  =2.9272276041...=\frac{1}{2\ln(2)}\left[
-1+4(1.2644997803...)\right]  ,
\]
omitting the fluctuation term. \ The value $2\lambda$ should be $3+\sqrt{5}$;
the subseries of Fibonacci terms with odd subscripts%
\[%
%TCIMACRO{\dsum \limits_{k=0}^{\infty}}%
%BeginExpansion
{\displaystyle\sum\limits_{k=0}^{\infty}}
%EndExpansion
\frac{1}{f_{2k+1}}=\frac{\sqrt{5}}{4}\left(
%TCIMACRO{\dsum \limits_{n=-\infty}^{\infty}}%
%BeginExpansion
{\displaystyle\sum\limits_{n=-\infty}^{\infty}}
%EndExpansion
\frac{1}{\lambda^{(n+1/2)^{2}}}\right)  ^{2}=1.8245151574...
\]
involves a Jacobi theta function $\vartheta_{2}(q)$ squared, where
$q=1/\lambda$. It turns out that $\nu$ and $\chi$ are linked via $\nu-1=\chi$;
we have \cite{Prdg-err, ChMd-err, Jnsn-err}
\[%
%TCIMACRO{\dsum \limits_{j=1}^{\infty}}%
%BeginExpansion
{\displaystyle\sum\limits_{j=1}^{\infty}}
%EndExpansion
\frac{(-1)^{j-1}}{j\left(  2^{j}-1\right)  }=%
%TCIMACRO{\dsum \limits_{k=1}^{\infty}}%
%BeginExpansion
{\displaystyle\sum\limits_{k=1}^{\infty}}
%EndExpansion
\ln\left(  1+2^{-k}\right)  =0.8688766526...=\frac{7.2271128245...}{12\ln
(2)}.
\]
Finally, a random variable $X$ with density $e^{-x}(e^{-x}-1+x)/(1-e^{-x}%
)^{2}$, $x\geq0$, has mean $\operatorname*{E}(X)=\pi^{2}/6$ and mean
fractional part \cite{Jnsn-err}
\[
\operatorname*{E}\left(  X-\left\lfloor X\right\rfloor \right)  =\frac{11}%
{24}+%
%TCIMACRO{\dsum \limits_{m=1}^{\infty}}%
%BeginExpansion
{\displaystyle\sum\limits_{m=1}^{\infty}}
%EndExpansion
\frac{\pi^{2}}{\sinh(2\pi^{2}m)^{2}}=\frac{11}{24}+(2.825535...)\times
10^{-16}.
\]
The distribution of $X$ is connected with the \textit{random assignment
problem} \cite{Alds-err, Prvn-err}\textit{.}

\textbf{5.15. Optimal Stopping Constants.} When discussing the expected rank
$R_{n}$, we assumed that no applicant would ever refuse a job offer! If each
applicant only accepts an offer with known probability $p$, then
\cite{Tamk-err}
\[
\lim_{n\rightarrow\infty}R_{n}=%
%TCIMACRO{\dprod \limits_{i=1}^{\infty}}%
%BeginExpansion
{\displaystyle\prod\limits_{i=1}^{\infty}}
%EndExpansion
\left(  1+\frac{2}{i}\frac{1+pi}{2-p+pi}\right)  ^{\frac{1}{1+pi}}%
\]
which is $6.2101994550...$ in the event that $p=1/2$. \ The same expression in
an integer parameter $p\geq2$ arises if instead we interview $p$ independent
streams of applicants; $\lim_{n\rightarrow\infty}R_{n}=2.6003019563...$ is
found for the bivariate case \cite{Gvndrj-err, Samu1-err}.

When discussing the full-information problem for Uniform $[0,1]$ variables, we
assumed that the number of applicants is known. If instead this itself is a
uniformly distributed variable on $\{1,2,\ldots,n\}$, then for the
\textquotedblleft nothing but the best objective\textquotedblright, the
asymptotic probability of success is \cite{Poro1-err, Poro2-err}
\[
(1-e^{a})\operatorname*{Ei}(-a)-(e^{-a}+a\operatorname*{Ei}(-a))(\gamma
+\ln(a)-\operatorname*{Ei}(a))=0.4351708055...
\]
where $a=2.1198244098...$ is the unique positive solution of the equation
\[
e^{a}(1-\gamma-\ln(a)+\operatorname*{Ei}(-a))-(\gamma+\ln
(a)-\operatorname*{Ei}(a))=1.
\]
It is remarkable that these constants occur in other, seemingly unrelated
versions of the secretary problem \cite{Petru-err, FrgHT-err, Samu2-err,
MzlvT-err}. Another relevant probabililty is \cite{MzlvT-err}%
\[
e^{-b}-\left(  e^{b}-b-1\right)  \operatorname*{Ei}(-b)=0.4492472188...
\]
where $b=1.3450166170...$ is the unique positive solution of the equation%
\[
\operatorname*{Ei}(-b)-\gamma-\ln(b)=-1.
\]
The corresponding full-information expected rank problem is called
\textit{Robbins' problem} \cite{Bruss-err, Swan-err}.

Suppose that you view successively terms of a sequence $X_{1}$, $X_{2}$,
$X_{3}$, ... of independent random variables with a common distribution
function $F$. You know the function $F$, and as $X_{k}$ is being viewed, you
must either stop the process or continue. If you stop at time $k$, you receive
a payoff $(1/k)%
%TCIMACRO{\tsum \nolimits_{j=1}^{k}}%
%BeginExpansion
{\textstyle\sum\nolimits_{j=1}^{k}}
%EndExpansion
X_{j}$. Your objective is to maximize the expected payoff. An optimal strategy
is to stop at the first $k$ for which $%
%TCIMACRO{\tsum \nolimits_{j=1}^{k}}%
%BeginExpansion
{\textstyle\sum\nolimits_{j=1}^{k}}
%EndExpansion
X_{j}\geq\alpha_{k}$, where $\alpha_{1} $, $\alpha_{2}$, $\alpha_{3}$, ... are
certain values depending on $F$. Shepp \cite{Shep-err, Walk-err} proved that
$\lim_{k\rightarrow\infty}\alpha_{k}/\sqrt{k}$ exists and is independent of
$F$ as long as $F$ has zero mean and unit variance; further,
\[
\lim_{k\rightarrow\infty}\frac{\alpha_{k}}{\sqrt{k}}%
=x=0.8399236756...=2(0.4199618378...)
\]
is the unique zero of $2x-\sqrt{2\pi}\left(  1-x^{2}\right)  \exp\left(
x^{2}/2\right)  \left(  1+\operatorname*{erf}(x/\sqrt{2})\right)  $. We wonder
if Shepp's constant can be employed to give a high-precision estimate of the
Chow-Robbins constant $2(0.7929535064...)-1=0.5859070128...$ \cite{Wsmn-err,
HggsW-err}, the value of the expected payoff for $F(-1)=F(1)=1/2.$

Consider a random binary string $Y_{1}Y_{2}Y_{3}\ldots Y_{n}$ with
$\operatorname*{P}(Y_{k}=1)=1-\operatorname*{P}(Y_{k}=0)$ independent of $k$
and $Y_{k}$ independent of the other $Y$s. Let $H$ denote the pattern
consisting of the digits
\[%
\begin{array}
[c]{ccccc}%
\underset{l}{\underbrace{1000...0}} &  & \text{or} &  & \underset
{l}{\underbrace{0111...1}}%
\end{array}
\]
and assume that its probability of occurrence for each $k$ is
\[
\operatorname*{P}\left(  Y_{k+1}Y_{k+2}Y_{k+3}\ldots Y_{k+l}=H\right)
=\frac1l\left(  1-\frac1l\right)  ^{l-1}\sim\frac1{el}=\frac{0.3678794411...}%
l.
\]
You observe sequentially the digits $Y_{1}$, $Y_{2}$, $Y_{3}$, ... one at a
time. You know the values $n$ and $l$, and as $Y_{k}$ is being observed, you
must either stop the process or continue. Your objective is to stop at the
final appearance of $H$ up to $Y_{n}$. Bruss \&\ Louchard \cite{BrLo-err}
determined a strategy that maximizes the probability of meeting this goal. For
$n\geq\beta l$, this success probability is
\[
\frac2{135}e^{-\beta}\left(  4-45\beta^{2}+45\beta^{3}\right)
=0.6192522709...
\]
as $l\rightarrow\infty$, where $\beta=3.4049534663...$ is the largest zero of
the cubic $45\beta^{3}-180\beta^{2}+90\beta+4$. Further, the interval
$[0.367...,0.619...]$ constitutes ``typical'' asymptotic bounds on success
probabilities associated with a wide variety of optimal stopping problems in strings.

Suppose finally that you view a sequence $Z_{1}$, $Z_{2}$, ..., $Z_{n}$ of
independent Uniform $[0,1]$ variables and that you wish to stop at a value of
$Z$ as large as possible. If you are a prophet (meaning that you have complete
foresight), then you know $Z_{n}^{\ast}=\max\{Z_{1},\ldots,Z_{n}\}$ beforehand
and clearly $\operatorname*{E}(Z_{n}^{\ast})\sim1-1/n$ as $n\rightarrow\infty
$. If you are a $1$-mortal (meaning that you have $1$ opportunity to choose a
$Z$ via stopping rules) and if you proceed optimally, then the value
$Z_{1}^{\ast}$ obtained satisfies $\operatorname*{E}(Z_{1}^{\ast})\sim1-2/n$.
If you are a $2$-mortal (meaning that you have $2 $ opportunities to choose
$Z$s and then take the maximum of these) and if you proceed optimally, then
the value $Z_{2}^{\ast}$ obtained satisfies $\operatorname*{E}(Z_{2}^{\ast
})\sim1-c/n$, where \cite{Assaf-err}
\[
c=\frac{2\xi}{\xi+2}=1.1656232877...
\]
and $\xi=2.7939976526...$ is the unique positive solution of the equation
\[
\left(  \frac{2}{\xi}+1\right)  \ln\left(  \frac{\xi}{2}+1\right)  =\frac
{3}{2}.
\]
The performance improvement in having two choices over just one is
impressive:\ $c$ is much closer to $1$ than $2$! See also \cite{Fi22-err,
OStp1-err, OStp2-err, OStp3-err}.

\textbf{5.16. Extreme Value Constants.} The median of the Gumbel distribution
is $-\ln(\ln(2))=0.3665129205...$.

\textbf{5.17. Pattern-Free Word Constants.} We now have improved bounds
$1.30173<S<1.30178858$ and $1.457567<C<1.45757921$ \cite{UGrm-err, XSun-err,
CRUG-err, OcRx-err, Klpkv-err, GrmHeu-err} and precise estimates%
\[%
\begin{array}
[c]{ccc}%
T_{L}=\dfrac{1}{11}\dfrac{\ln\left(  \rho(A\,B^{10})\right)  }{\ln
(2)}=1.273553265..., &  & T_{U}=\dfrac{1}{2}\dfrac{\ln\left(  \rho
(A\,B)\right)  }{\ln(2)}=1.332240491...
\end{array}
\]
where $A$, $B$ are known $20\times20$ integer matrices and $\rho$ denotes
spectral radius \cite{Berstel-err, JPB-err, GuPr-err}. \ The set of quaternary
words avoiding abelian squares grows exponentially (although $h(n)^{1/n}$ is
not well understood as length $n\rightarrow\infty$); the set of binary words
avoiding abelian fourth powers likewise is known to grow exponentially
\cite{Currie-err}.

\textbf{5.18. Percolation Cluster Density Constants.} Approximating $p_{c}$
for site percolation on the square lattice continues to draw attention
\cite{MNP-1, MNP-2, MNP-3, MNP-4, MNP-5, MNP-6}; for the hexagonal lattice,
$p_{c}=0.697043...$ improves upon the estimate given on p. 373. More about
mean cluster densities can be found in \cite{Tggmnn-err, ChngSh-err}. An
integral similar to that for $\kappa_{B}(p_{c})$ on the triangular lattice
appears in \cite{Ken-err}.

Hall's bounds for $\lambda_{c}$ on p. 375 can be written as $1.642<4\pi
\,\lambda_{c}<10.588$ and the best available estimate is $4\pi\,\lambda
_{c}=4.51223...$ \cite{BlbRdn-err, BBWlts-err}. Older references on 2D and 3D
continuum percolation include \cite{PQR0-err, PQR1-err, PQR2-err, PQR3-err,
PQR4-err, PQR5-err}. See also \cite{ZSclrd-err, RdnWlt-err, FeDgBt-err,
NBall-err, BJRdn-err}.

Two infinite $0$-$1$ sequences $X$, $Y$ are called \textit{compatible} if $0$s
can be deleted from $X$ and/or from $Y$ in such a way that the resulting $0
$-$1$ sequences $X^{\prime}$, $Y^{\prime}$ never have a $1$ in the same
position. For example, the sequences $X=000110\ldots$ and $Y=110101\ldots$ are
not compatible. Assume that $X$ and $Y$ are randomly generated with each
$X_{i}$, $Y_{j}$ independent and $\operatorname*{P}(X_{i}=1)=\operatorname*{P}%
(Y_{j}=1)=p$. Intuition suggests that $X$ and $Y$ are compatible with positive
probability if and only if $p$ is suitably small. What is the supremum $p^{*}$
of such $p $? It is known \cite{PWnklr-err, PGacs1-err, PGacs2-err,
IPeters-err} that $100^{-400}<p^{*}<1/2$; simulation indicates
\cite{JTrmp2-err} that $0.3<p^{*}<0.305$.

Consider what is called \textit{bootstrap percolation} on the $d$-dimensional
cubic lattice with $n^{d}$ vertices: starting from a random set of initially
``infected'' sites, new sites become infected at each time step if they have
at least $d$ infected neighbors and infected sites remain infected forever.
Assume that vertices of the initial set were chosen independently, each with
probability $p$. What is the critical probability $p_{c}(n,d)$ for which the
likelihood that the entire lattice is subsequently infected exceeds $1/2$?
Holroyd \cite{Hlryd-err} and Balogh, Bollob\'as \&\ Morris \cite{BBsM-err}
proved that
\[%
\begin{array}
[c]{ccc}%
p_{c}(n,2)=\dfrac{\pi^{2}/18+o(1)}{\ln(n)}, &  & p_{c}(n,3)=\dfrac{\mu
+o(1)}{\ln(\ln(n))}%
\end{array}
\]
as $n\rightarrow\infty$, where
\[
\mu=-%
%TCIMACRO{\dint \limits_{0}^{\infty}}%
%BeginExpansion
{\displaystyle\int\limits_{0}^{\infty}}
%EndExpansion
\ln\left(  \frac12-\frac{e^{-2x}}2+\frac12\sqrt{1+e^{-4\,x}-4e^{-3x}+2e^{-2x}%
}\right)  dx=0.4039127202....
\]
A closed-form expression for $\mu$ remains open.

\textbf{5.19. Klarner's Polyomino Constant.} A new estimate $4.0625696...$ for
$\alpha$ is reported in \cite{Jens-err} and a new rigorous lower bound of
$3.980137...$ in \cite{BMRR-err}. The number $\bar{A}(n)$ of row-convex $n
$-ominoes satisfies \cite{Hic-err}
\[%
\begin{array}
[c]{ccc}%
\bar{A}(n)=5\bar{A}(n-1)-7\bar{A}(n-2)+4\bar{A}(n-3), &  & n\geq5,
\end{array}
\]
with $\bar{A}(1)=1$, $\bar{A}(2)=2$, $\bar{A}(3)=6$ and $\bar{A}(4)=19$; hence
$\bar{A}(n)\sim u\,v^{n}$ as $n\rightarrow\infty$, where $v=3.2055694304...$
is the unique real zero of $x^{3}-5x^{2}+7x-4$ and $u=(41v^{2}%
-129v+163)/944=0.1809155018...$. While the multiplicative constant for
parallelogram $n$-ominoes is now known to be $0.2974535058...$, corresponding
improved accuracy for convex $n$-ominoes evidently remains open. A Central
Limit Theorem applies to the perimeter of a random parallelogram $n$-omino
$S$, which turns out to be normal with mean $(0.8417620156...)n$ and standard
deviation $(0.4242065326...)\sqrt{n}$ in the limit as $n\rightarrow\infty$.
\ Hence $S$ is expected to resemble a slanted stack of fairly short rods
\cite{FljSdg-err}. Again, corresponding quantities for a random convex
$n$-omino are not known. More on coin fountains and the constant
$0.5761487691...$ can be found in \cite{ParLeSc-err, BHShSn-err, MSprgnl-err,
DucFGS-err}.

\textbf{5.20. Longest Subsequence Constants.} Regarding common subsequences,
Lueker \cite{Luekr-err, KwSt-err} showed that $0.7880\leq\gamma_{2}\leq
0.8263$. The Sankoff-Mainville conjecture that $\lim_{k\rightarrow\infty
}\gamma_{k}k^{1/2}=2$ was proved by Kiwi, Loebl \&\ Matousek \cite{KwLM-err};
the constant $2$ arises from a connection with increasing subsequences. A
deeper connection with the Tracy-Widom distribution from random matrix theory
has now been confirmed \cite{MajN-err}:
\[%
\begin{array}
[c]{ccc}%
\operatorname*{E}(\lambda_{n,k})\sim2k^{-1/2}n+c_{1}k^{-1/6}n^{1/3}, &  &
\operatorname*{Var}(\lambda_{n,k})\sim c_{0}k^{-1/3}n^{2/3}%
\end{array}
\]
where $k\rightarrow\infty$, $n\rightarrow\infty$ in such a way that
$n/k^{1/2}\rightarrow0$.

Define $\lambda_{n,k,r}$ to be the length of the longest common subsequence
$c$ of $a$ and $b$ subject to the constraint that, if $a_{i}=b_{j}$ are paired
when forming $c$, then $|i-j|\leq r$. Define as well $\gamma_{k,r}%
=\lim_{n\rightarrow\infty}\operatorname*{E}(\lambda_{n,k,r})/n$. It is not
surprising \cite{Blsk-err} that $\lim_{r\rightarrow\infty}\gamma_{k,r}%
=\gamma_{k}$. Also, $\gamma_{2,1}=7/10$, but exact values for $\gamma_{3,1}$,
$\gamma_{4,1}$, $\gamma_{2,2}$ and $\gamma_{2,3}$ remain open.

Here is a geometric formulation \cite{AmBry-err}.\ Given $N$ independent
uniform random points $\{z_{j}\}_{j=1}^{N}$ in the unit square $S$, an
\textit{increasing chain} is a polygonal path that links the southwest and
northeast corners of $S$ and whose other vertices are $\{z_{j_{i}}\}_{i=1}%
^{k}$, $0\leq k\leq N$, assuming both $\operatorname{Re}(z_{j_{i}})$ and
$\operatorname{Im}(z_{j_{i}})$ are strictly increasing with $i$. \ The
\textit{length} of the chain is simply $k$. \ A variation of this requires
that $\operatorname{Re}(z_{j_{i}})>\operatorname{Im}(z_{j_{i}})$ always
(equivalently, the path never leaves the lower isosceles right triangle).
\ If, further, the region bounded by the path and the diagonal (hypotenuse) is
convex, then the path is a \textit{convex chain}. \ Under such circumstances,
it seems likely that the length $L_{N}^{\prime}$ of longest convex chains
satisfies
\[
\lim_{N\rightarrow\infty}N^{-1/3}\operatorname*{E}(L_{N}^{\prime})=3
\]
(we know that the limit exists and lies between $1.5772$ and $3.4249$). \ This
result seems to be true as well for chains that link two corners of arbitrary
(non-isosceles) triangles.

The Tracy-Widom distribution (specifically, $F_{\text{GOE}}(x)$ as described
in \cite{Fi12-err}) seems to play a role in other combinatorial problems
\cite{Srnk-err, Nvkf-err, Bber-err}, although the data is not conclusive. See
also \cite{Stnly1-err, Stnly2-err, Stnly3-err}.

\textbf{5.21. }$k$\textbf{-Satisfiability Constants.} On the one hand, the
lower bound for $r_{c}(3)$ was improved to 3.42 in \cite{KKL-err} and further
improved to 3.52 in \cite{HaSo-err}. On the other hand, the upper bound 4.506
for $r_{c}(3)$ in \cite{DBM-err} has not been confirmed; the preceding two
best upper bounds were 4.596 \cite{JSV-err} and 4.571 \cite{KKSVZ-err}. See
\cite{IbrKKM-err} for recent work on XOR-SAT.

\textbf{5.22. Lenz-Ising Constants.} Improved estimates for $K_{c}=0.11392...
$, $0.09229...$, $0.077709...$ when $d=5$, $6$, $7$ appear in
\cite{Ising0-err}. Define Ising susceptibility integrals
\[
D_{n}=\frac{4}{n!}%
%TCIMACRO{\dint \limits_{0}^{\infty}}%
%BeginExpansion
{\displaystyle\int\limits_{0}^{\infty}}
%EndExpansion%
%TCIMACRO{\dint \limits_{0}^{\infty}}%
%BeginExpansion
{\displaystyle\int\limits_{0}^{\infty}}
%EndExpansion
\cdots%
%TCIMACRO{\dint \limits_{0}^{\infty}}%
%BeginExpansion
{\displaystyle\int\limits_{0}^{\infty}}
%EndExpansion
\frac{%
%TCIMACRO{\tprod \nolimits_{i<j}}%
%BeginExpansion
{\textstyle\prod\nolimits_{i<j}}
%EndExpansion
\left(  \frac{x_{i}-x_{j}}{x_{i}+x_{j}}\right)  ^{2}}{\left(
%TCIMACRO{\tsum \nolimits_{k=1}^{n}}%
%BeginExpansion
{\textstyle\sum\nolimits_{k=1}^{n}}
%EndExpansion
(x_{k}+1/x_{k})\right)  ^{2}}\frac{dx_{1}}{x_{1}}\frac{dx_{2}}{x_{2}}%
\ldots\frac{dx_{n}}{x_{n}}%
\]
(also known as McCoy-Tracy-Wu integrals). Clearly $D_{1}=2$ and $D_{2}=1/3$;
we also have
\[
\frac{D_{3}}{8\pi^{2}}=\frac{8+4\pi^{2}/3-27L_{-3}(2)}{8\pi^{2}}%
=0.000814462565...,
\]%
\[
\frac{D_{4}}{16\pi^{3}}=\frac{4\pi^{2}/9-1/6-7\zeta(3)/2}{16\pi^{3}%
}=0.000025448511...,
\]
and the former is sometimes called the ferromagnetic constant
\cite{Ising1-err, Ising2-err}. These integrals are important because
\cite{Ising3-err, Ising4-err}
\[
\pi%
%TCIMACRO{\dsum \limits_{n\equiv1\operatorname*{mod}2}}%
%BeginExpansion
{\displaystyle\sum\limits_{n\equiv1\operatorname*{mod}2}}
%EndExpansion
\frac{D_{n}}{(2\pi)^{n}}=1.0008152604...=2^{3/8}\ln(1+\sqrt{2})^{7/4}%
(0.9625817323...),
\]%
\[
\pi%
%TCIMACRO{\dsum \limits_{n\equiv0\operatorname*{mod}2}}%
%BeginExpansion
{\displaystyle\sum\limits_{n\equiv0\operatorname*{mod}2}}
%EndExpansion
\frac{D_{n}}{(2\pi)^{n}}=\frac{1.0009603287...}{12\pi}=2^{3/8}\ln(1+\sqrt
{2})^{7/4}(0.0255369745...)
\]
and such constants $c_{0}^{+}$, $c_{0}^{-}$ were earlier given in terms of a
solution of the Painlev\'{e} III differential equation.

The number of spanning trees in the $d$-dimensional cubic lattice with
$N=n^{d} $ vertices grows asymptotically as $\exp(h_{d}N)$, where
\begin{align*}
h_{d}  &  =\frac1{(2\pi)^{d}}%
%TCIMACRO{\dint \limits_{-\pi}^{\pi}}%
%BeginExpansion
{\displaystyle\int\limits_{-\pi}^{\pi}}
%EndExpansion%
%TCIMACRO{\dint \limits_{-\pi}^{\pi}}%
%BeginExpansion
{\displaystyle\int\limits_{-\pi}^{\pi}}
%EndExpansion
\cdots%
%TCIMACRO{\dint \limits_{-\pi}^{\pi}}%
%BeginExpansion
{\displaystyle\int\limits_{-\pi}^{\pi}}
%EndExpansion
\ln\left(  2d-2%
%TCIMACRO{\dsum \limits_{k=1}^{d}}%
%BeginExpansion
{\displaystyle\sum\limits_{k=1}^{d}}
%EndExpansion
\cos(\theta_{k})\right)  d\theta_{1}\,d\theta_{2}\cdots\,d\theta_{d}\ \\
\  &  =\ln(2d)+%
%TCIMACRO{\dint \limits_{0}^{\infty}}%
%BeginExpansion
{\displaystyle\int\limits_{0}^{\infty}}
%EndExpansion
\frac{e^{-t}}t\left(  1-I_{0}\left(  \frac td\right)  ^{d}\right)  dt.
\end{align*}
Note the similarity with the formula for $m_{d}$ on p. 323. We have
\cite{FlkLyn-err}
\[%
\begin{array}
[c]{ccc}%
h_{2}=4G/\pi=1.1662436161..., &  & h_{3}=1.6733893029...,
\end{array}
\]
\[%
\begin{array}
[c]{ccccc}%
h_{4}=1.9997076445..., &  & h_{5}=2.2424880598..., &  & h_{6}=2.4366269620....
\end{array}
\]
Other forms of $h_{3}$ have appeared in the literature \cite{Rsngr1-err,
Rsngr2-err, Rsngr3-err}:
\[%
\begin{array}
[c]{ccc}%
h_{3}-\ln(2)=0.9802421224..., &  & h_{3}-\ln(2)-\ln(3)=-0.1183701662....
\end{array}
\]
The corresponding constant for the two-dimensional triangular lattice is
\cite{GlssWu-err}
\[
\hat h=\frac12\ln(3)+\frac6\pi\operatorname*{Ti}\nolimits_{2}\left(
\frac1{\sqrt{3}}\right)  =1.6153297360...
\]
where $\operatorname*{Ti}\nolimits_{2}(x)$ is the inverse tangent integral
(discussed on p. 57). Results for other lattices are known \cite{ChShrk-err,
ChWang-err}; we merely mention a new closed-form evaluation:
\begin{align*}
&  \ \ \ \ \frac{\ln(2)}2+\frac1{16\pi^{2}}%
%TCIMACRO{\dint \limits_{-\pi}^{\pi}}%
%BeginExpansion
{\displaystyle\int\limits_{-\pi}^{\pi}}
%EndExpansion%
%TCIMACRO{\dint \limits_{-\pi}^{\pi}}%
%BeginExpansion
{\displaystyle\int\limits_{-\pi}^{\pi}}
%EndExpansion
\ln\left[  7-3\cos(\theta)-3\cos(\varphi)-\cos(\theta)\cos(\varphi)\right]
\,d\theta\,d\varphi\\
\  &  =\frac G\pi+\frac12\ln(\sqrt{2}-1)+\frac1\pi\operatorname*{Ti}%
\nolimits_{2}(3+2\sqrt{2})=0.7866842753...
\end{align*}
associated with a certain tiling of the plane by squares and octagons.

\textbf{5.23. Monomer-Dimer Constants.} Friedland \&\ Peled \cite{FP-err} and
other authors \cite{FG-err, HLLB-err, Kong2-err, FKLM-err, GK-err, BtrPnc-err}
revisited Baxter's computation of $A$ and confirmed that $\ln
(A)=0.66279897...$. They also examined the three-dimensional analog,
$A^{\prime}$, of $A$, yielding $\ln(A^{\prime})=0.785966...$. \ Butera,
Federbush \& Pernici \cite{BtrFdP-err} estimated $\lambda=0.449...$ which is
inconsistent with some earlier values.

For odd $n$, Tzeng \&\ Wu \cite{Wu1-err, Wu2-err} found the number of dimer
arrangements on the $n\times n$ square lattice with exactly one monomer on the
boundary. If the restriction that the monomer lie on the boundary is removed,
then enumeration is vastly more difficult; Kong \cite{Kong1-err} expressed the
possibility that this problem might be solvable someday. Wu \cite{Wu3-err}
examined dimers on various other two-dimensional lattices.

A \textit{trimer} consists of three adjacent collinear vertices of the square
lattice. The trimer-covering analog of the entropy $\exp(2G/\pi)=1.7916...$ is
$1.60...$, which is variously written as $\exp(0.475...)$ or as $\exp
(3\cdot0.15852...)$ \cite{TRI1-err, TRI2-err, TRI3-err, TRI4-err, TRI5-err,
TRI6-err}.

Ciucu \&\ Wilson \cite{Ciucu-err} discovered a constant $0.9587407138...$ that
arises with regard to the asymptotic decay of monomer-monomer correlation
\textquotedblleft in a sea of dimers\textquotedblright\ on what is called the
critical Fisher lattice.

\textbf{5.24. Lieb's Square Ice Constant. }More on counting Eulerian
orientations is found in \cite{RWhitty-err, FlsnrZk-err}.

\textbf{5.25. Tutte-Beraha Constants.} For any positive integer $r$, there is
a best constant $C(r)\,$ such that, for each graph of maximum degree $\leq r$,
the complex zeros of its chromatic polynomial lie in the disk $|z|\leq C(r)$.
Further, $K=\lim_{r\rightarrow\infty}C(r)/r$ exists and $K=7.963906...$ is the
smallest number for which
\[
\inf_{\alpha>0}\frac1\alpha%
%TCIMACRO{\dsum \limits_{n=2}^{\infty}}%
%BeginExpansion
{\displaystyle\sum\limits_{n=2}^{\infty}}
%EndExpansion
e^{\alpha n}K^{-(n-1)}\frac{n^{n-1}}{n!}\leq1.
\]
Sokal \cite{Sol-err} proved all of the above, answering questions raised in
\cite{BDS-err, BRW-err}. See also \cite{Royle-err}.

\textbf{6.1. Gauss' Lemniscate Constant.} Consider the following game
\cite{Bla-err}. Players $A$ and $B$ simultaneously choose numbers $x$ and $y$
in the unit interval; $B$ then pays $A$ the amount $|x-y|^{1/2}$. The value of
the game (that is, the expected payoff, assuming both players adopt optimal
strategies) is $M/2=0.59907....$ Also, let $\xi_{1}$, $\xi_{2}$, $\ldots$,
$\xi_{n}$, $\eta_{1}$, $\eta_{2}$, $\ldots$, $\eta_{n}$ be distinct points in
the plane and construct, with these points as centers, squares of side $s$ and
of arbitrary orientation that do not overlap. Then
\[
s\leq\frac{L}{\sqrt{2}}\left(  \frac{%
%TCIMACRO{\dprod \limits_{i=1}^{n}}%
%BeginExpansion
{\displaystyle\prod\limits_{i=1}^{n}}
%EndExpansion%
%TCIMACRO{\dprod \limits_{j=1}^{n}}%
%BeginExpansion
{\displaystyle\prod\limits_{j=1}^{n}}
%EndExpansion
|\xi_{i}-\eta_{j}|}{%
%TCIMACRO{\dprod \limits_{i<j}}%
%BeginExpansion
{\displaystyle\prod\limits_{i<j}}
%EndExpansion
|\xi_{i}-\xi_{j}|\cdot%
%TCIMACRO{\dprod \limits_{i<j}}%
%BeginExpansion
{\displaystyle\prod\limits_{i<j}}
%EndExpansion
|\eta_{i}-\eta_{j}|}\right)  ^{1/n}%
\]
and the constant $L/\sqrt{2}=1.85407...=r$ is best possible \cite{Neh-err}.
\ More on $\tilde{r}=1.52995...$ appears in \cite{Nim-err}.

\textbf{6.2. Euler-Gompertz Constant.} We do not yet know whether $C_{2}$ is
transcendental, but it cannot be true that both $\gamma$ and $C_{2}$ are
algebraic \cite{LgrsEC-err, Aptrv-err, TRvoal-err, HPHP2-err}. This result
evidently follows from Mahler \cite{Mah1-err}, who in turn was reporting on
work by Shidlovski \cite{Shdlv-err}. \ Generalizations of $C_{2}$ include
\cite{GldQntc-err, ABKMW-err}%
\[
\frac{1}{(m-1)!}%
%TCIMACRO{\dint \limits_{0}^{\infty}}%
%BeginExpansion
{\displaystyle\int\limits_{0}^{\infty}}
%EndExpansion
t^{m-1}e^{1-e^{t}}dt=\left\{
\begin{array}
[c]{lll}%
0.2659653850... &  & \text{if }m=2,\\
0.0967803251... &  & \text{if }m=3,\\
0.0300938139... &  & \text{if }m=4
\end{array}
\right.
\]
which pertain to statistics governing restricted permutations and set
partitions. \ For actuarial background and history, consult \cite{Ciecka-err}.

The two quantities
\[%
\begin{array}
[c]{ccc}%
I_{0}(2)=%
%TCIMACRO{\dsum \limits_{k=0}^{\infty}}%
%BeginExpansion
{\displaystyle\sum\limits_{k=0}^{\infty}}
%EndExpansion
\dfrac{1}{(k!)^{2}}=2.2795853023..., &  & J_{0}(2)=%
%TCIMACRO{\dsum \limits_{k=0}^{\infty}}%
%BeginExpansion
{\displaystyle\sum\limits_{k=0}^{\infty}}
%EndExpansion
\dfrac{(-1)^{k}}{(k!)^{2}}=0.2238907791...
\end{array}
\]
are similar, but only the first is associated with continued fractions. Here
is an interesting occurrence of the second: letting \cite{AWMv-err}
\[%
\begin{tabular}
[c]{llll}%
$a_{0}=a_{1}=1,$ &  & $a_{n}=n\,a_{n-1}-a_{n-2}$ & for $n\geq2,$%
\end{tabular}
\]
we have $\lim_{n\rightarrow\infty}a_{n}/n!=J_{0}(2)$. The constant
\[
C_{2}=%
%TCIMACRO{\dint \limits_{0}^{\infty}}%
%BeginExpansion
{\displaystyle\int\limits_{0}^{\infty}}
%EndExpansion
\frac{e^{-x}}{1+x}dx=%
%TCIMACRO{\dint \limits_{0}^{1}}%
%BeginExpansion
{\displaystyle\int\limits_{0}^{1}}
%EndExpansion
\frac{1}{1-\ln(y)}dy=0.5963473623...
\]
unexpectedly appears in \cite{MM75-err}, and the constant $2(1-C_{1}%
)=0.6886409151...$ unexpectedly appears in \cite{MBhL-err}. Also, the
divergent alternating series $0!-2!+4!-6!+-\cdots$ has value \cite{Knpp-err}
\[%
%TCIMACRO{\dint \limits_{0}^{\infty}}%
%BeginExpansion
{\displaystyle\int\limits_{0}^{\infty}}
%EndExpansion
\frac{e^{-x}}{1+x^{2}}dx=%
%TCIMACRO{\dint \limits_{0}^{1}}%
%BeginExpansion
{\displaystyle\int\limits_{0}^{1}}
%EndExpansion
\frac{1}{1+\ln(y)^{2}}dy=0.6214496242...
\]
and, similarly, the series $1!-3!+5!-7!+-\cdots$ has value
\[%
%TCIMACRO{\dint \limits_{0}^{\infty}}%
%BeginExpansion
{\displaystyle\int\limits_{0}^{\infty}}
%EndExpansion
\frac{x\,e^{-x}}{1+x^{2}}dx=-%
%TCIMACRO{\dint \limits_{0}^{1}}%
%BeginExpansion
{\displaystyle\int\limits_{0}^{1}}
%EndExpansion
\frac{\ln(y)}{1+\ln(y)^{2}}dy=0.3433779615....
\]

Let $G(z)$ denote the standard normal distribution function and
$g(z)=G^{\prime}(z)$. If $Z$ is distributed according to $G$, then
\cite{KGDD-err}
\[
\operatorname*{E}\left(  Z\;|\;Z>1\right)  =\frac{g(1)}{G(-1)}=\frac1{C_{1}%
}=1.5251352761...,
\]
\[
\operatorname*{E}\left(  \left\{
\begin{array}
[c]{ccc}%
Z &  & \text{if }Z>1,\\
0 &  & \text{otherwise}%
\end{array}
\right.  \right)  =g(1)=\frac1{\sqrt{2\pi e}}=0.2419707245...,
\]
\[
\operatorname*{E}\left(  \max\left\{  Z-1,0\right\}  \right)
=g(1)-G(-1)=0.0833154705...
\]
which contrast interestingly with earlier examples.

\textbf{6.3. Kepler-Bouwkamp Constant.} Additional references include
\cite{Kp1-err, Kp2-err, Kp3-err, ChbrS-err} and another representation is
\cite{Kp4-err}
\[
\rho=\frac{3^{10}\sqrt{3}}{2^{7}5^{2}7\,11\,\pi}\exp\left[  -\sum
_{k=1}^{\infty}\frac{\left(  \zeta(2k)-1-2^{-2k}-3^{-2k}\right)  2^{2k}\left(
\lambda(2k)-1-3^{-2k}\right)  }{k}\right]  ;
\]
the series converges at the same rate as a geometric series with ratio $1/100
$. A relevant inequality is \cite{BBVW-err}
\[%
%TCIMACRO{\dint \limits_{0}^{\infty}}%
%BeginExpansion
{\displaystyle\int\limits_{0}^{\infty}}
%EndExpansion
\cos(2x)%
%TCIMACRO{\dprod \limits_{j=1}^{\infty}}%
%BeginExpansion
{\displaystyle\prod\limits_{j=1}^{\infty}}
%EndExpansion
\cos\left(  \dfrac{x}{j}\right)  dx<\frac{\pi}{8}%
\]
and the difference is less than $10^{-42}$! Powers of two are featured in the
following: \cite{FM96-err, MF06-err}
\[%
%TCIMACRO{\dint \limits_{0}^{\pi}}%
%BeginExpansion
{\displaystyle\int\limits_{0}^{\pi}}
%EndExpansion
\left\vert
%TCIMACRO{\dprod \limits_{m=0}^{n}}%
%BeginExpansion
{\displaystyle\prod\limits_{m=0}^{n}}
%EndExpansion
\sin\left(  2^{m}x\right)  \right\vert dx=\kappa\,\lambda^{n}\left(
1+o(1)\right)
\]
as $n\rightarrow\infty$, where $\kappa>0$ and $0.654336<\lambda<0.663197$. A
prime analog of $\rho$ is \cite{Ktsn2-err, Mathr5-err, Doslic-err}
\[%
%TCIMACRO{\dprod \limits_{p\geq3}}%
%BeginExpansion
{\displaystyle\prod\limits_{p\geq3}}
%EndExpansion
\cos\left(  \frac{\pi}{p}\right)  =0.3128329295...=(3.1965944300...)^{-1}%
\]
and variations abound. Also, the conjecture $%
%TCIMACRO{\tprod \nolimits_{k\geq1}}%
%BeginExpansion
{\textstyle\prod\nolimits_{k\geq1}}
%EndExpansion
\tan(k)=0$ is probably false \cite{DBail-err}.

\textbf{6.4. Grossman's Constant. }Somos \cite{Somos-err} examined the pair of
recurrences%
\[%
\begin{array}
[c]{ccccccc}%
a_{n}=a_{n-1}+b_{n-1}, &  & b_{n}=-a_{n-1}b_{n-1}, &  & a_{0}=-1, &  & b_{0}=x
\end{array}
\]
and conjectured that there exists a unique real number $x=\xi$ for which both
sequences converge (quadratically) to $0$, namely $\xi=0.349587046...$. \ The
resemblance to the AGM\ recursion is striking.

\textbf{6.5. Plouffe's Constant. }This constant is included in a fascinating
mix of ideas by Smith \cite{Smi-err}, who claims that ``angle-doubling'' one
bit at a time was known centuries ago to Archimedes and was implemented
decades ago in binary cordic algorithms (also mentioned in section 5.14).
Another constant of interest is $\arctan(\sqrt{2})=0.9553166181...$, which is
the base angle of a certain isosceles spherical triangle (in fact, the unique
non-Euclidean triangle with rational sides and a single right angle).

Chowdhury \cite{Chwd-err} generalized his earlier work on bitwise XOR sums and
the logistic map: A sample new result is
\[%
%TCIMACRO{\dsum \limits_{n=0}^{\infty}}%
%BeginExpansion
{\displaystyle\sum\limits_{n=0}^{\infty}}
%EndExpansion
\frac{\rho(b_{n}b_{n-1})}{2^{n+1}}=\frac1{4\pi}\oplus\frac1\pi
\]
where $b_{n}=\cos(2^{n})$. The right-hand side is computed merely by shifting
the binary expansion of $1/\pi$ two places (to obtain $1/(4\pi)$) and adding
modulo two without carries (to find the sum).

\textbf{6.6. Lehmer's Constant.} Rivoal \cite{Riv-err} has studied the link
between the rational approximations of a positive real number $x$ coming from
the continued cotangent representation of $x$, and the usual convergents that
proceed from the regular continued fraction expansion of $x$.

\textbf{6.7. Cahen's Constant}. The usual meaning of \textquotedblleft Let $w$
be an infinite sequence\textquotedblright\ (as \textit{fixed} from the start)
became distorted at the bottom of page 435. \ Let $n\geq0$. \ The value
$w_{n}$ isn't actually needed until $q_{n+1}$ is calculated; once this is
done, the values $w_{n+1}$ \&\ $w_{n+2}$ become known; these, in turn, give
rise to $q_{n+2}$ \&\ $q_{n+3}$ and so forth. We miss the author of
\cite{So7-err} and his friendship.

\textbf{6.8. Prouhet-Thue-Morse Constant.} A follow-on to Allouche
\&\ Shallit's survey appears in \cite{AllchT-err}. Simple analogs of the
Woods-Robbins and Flajolet-Martin formulas are \cite{So0-err}
\[%
\begin{array}
[c]{ccc}%
%TCIMACRO{\dprod \limits_{m=1}^{\infty}}%
%BeginExpansion
{\displaystyle\prod\limits_{m=1}^{\infty}}
%EndExpansion
\left(  \dfrac{2m}{2m-1}\right)  ^{(-1)^{m}}=\dfrac{\sqrt{2}\pi^{3/2}}%
{\Gamma(1/4)^{2}}, &  &
%TCIMACRO{\dprod \limits_{m=1}^{\infty}}%
%BeginExpansion
{\displaystyle\prod\limits_{m=1}^{\infty}}
%EndExpansion
\left(  \dfrac{2m}{2m+1}\right)  ^{(-1)^{m}}=\dfrac{\Gamma(1/4)^{2}}%
{2^{5/2}\sqrt{\pi}};
\end{array}
\]
we wonder about the outcome of exponent sequences other than $(-1)^{m}$ or
$(-1)^{t_{m}}$. See also \cite{ChbrS-err, AllSdw1-err, AllSdw2-err,
KVitem-err, Riasat1-err, Riasat2-err}. \ Beware of a shifted version, used in
\cite{AllchP-err}, of our paper folding sequence $(-1)^{s_{m}}$.

Just as the Komornik-Loreti constant $1.7872316501...$ is the unique positive
solution of
\[%
%TCIMACRO{\dsum \limits_{n=1}^{\infty}}%
%BeginExpansion
{\displaystyle\sum\limits_{n=1}^{\infty}}
%EndExpansion
t_{n}q^{-n}=1,
\]
the (transcendental) constants $2.5359480481...$ and $2.9100160556...$ are
unique positive solutions of \cite{VKmrL-err}
\[%
\begin{array}
[c]{ccc}%
%TCIMACRO{\dsum \limits_{n=1}^{\infty}}%
%BeginExpansion
{\displaystyle\sum\limits_{n=1}^{\infty}}
%EndExpansion
(1+t_{n}-t_{n-1})q^{-n}=1, &  &
%TCIMACRO{\dsum \limits_{n=1}^{\infty}}%
%BeginExpansion
{\displaystyle\sum\limits_{n=1}^{\infty}}
%EndExpansion
(1+t_{n})q^{-n}=1.
\end{array}
\]
These correspond to $q$-developments with $0\leq\varepsilon_{n}\leq2$ and
$0\leq\varepsilon_{n}\leq3$ (although our numerical estimates differ from
those in \cite{SBaker-err}). \ Incidently, the smallest $q>\varphi$ possessing
a countably infinite number of $q$-developments with $0\leq\varepsilon_{n}%
\leq1$ is algebraic of degree $5$ \cite{ZWLuB-err}.

\textbf{6.9. Minkowski-Bower Constant. }The question mark satisfies the
functional equation \cite{GbKsTy-err}
\[
?(x)=\left\{
\begin{array}
[c]{lll}%
\dfrac12?\left(  \dfrac x{1-x}\right)  &  & \text{if }0\leq x\leq\dfrac12,\\
1-\dfrac12?\left(  \dfrac{1-x}x\right)  &  & \text{if }\dfrac12<x\leq1.
\end{array}
\right.
\]
See \cite{BG-err, Mrd-err, Pnti-err} for generalizations. Kinney
\cite{Knny-err} examined the constant
\[
\alpha=\frac12\left(
%TCIMACRO{\dint \limits_{0}^{1}}%
%BeginExpansion
{\displaystyle\int\limits_{0}^{1}}
%EndExpansion
\log_{2}(1+x)d?(x)\right)  ^{-1}%
\]
which acts as a threshold for Hausdorff dimension (of sets $\subset[0,1]$).
Lagarias \cite{Lagrs-err} computed that $0.8746<\alpha<0.8749$; the estimate
$0.875$ appears in \cite{TyUtz-err, PVdrB-err, KssStr-err, DshM-err};
Alkauskas \cite{Alkss-err} improved this approximation to $0.8747163051...$.
See also \cite{Fi20-err}.

\textbf{6.10. Quadratic Recurrence Constants}. In our asymptotic expansion for
$g_{n}$, the final coefficient should be $138$, not $137$ \cite{SmsQ1-err,
SmsQ2-err}. The sequence $k_{n+1}=(1/n)k_{n}^{2}$, where $n\geq0 $, is
convergent if and only if
\[
|k_{0}|<%
%TCIMACRO{\dprod \limits_{j=1}^{\infty}}%
%BeginExpansion
{\displaystyle\prod\limits_{j=1}^{\infty}}
%EndExpansion
\left(  1+\frac{1}{j}\right)  ^{2^{-j}}=1.6616879496....
\]
Moreover, the sequence either converges to zero or diverges to infinity
\cite{Jack-err, FeMa-err}. A systematic study of threshold constants like
this, over a broad class of quadratic recurrences, has never been attempted.
The constant $1.2640847353...$ and Sylvester's sequence appear in an
algebraic-geometric setting \cite{Nill-err}. Also, results on Somos' sequences
are found in \cite{Hone1-err, Hone2-err} and on the products
\[%
\begin{array}
[c]{ccc}%
1^{1/2}2^{1/4}3^{1/8}\cdots=1.6616879496..., &  & 1^{1/3}2^{1/9}3^{1/27}%
\cdots=1.1563626843...
\end{array}
\]
in \cite{So0-err, So6-err, SmsQ3-err, SmsQ4-err}. \ Calculating the area of
$M$ continues to attract attention \cite{Mndlbrt1-err, Mndlbrt2-err,
Mndlbrt3-err}.

\textbf{6.11. Iterated Exponential Constants.} Consider the recursion
\[%
\begin{array}
[c]{ccc}%
a_{1}=1, &  & a_{n}=a_{n-1}\exp\left(  \dfrac{1}{e\,a_{n-1}}\right)
\end{array}
\]
for $n\geq2$. It is known that \cite{NiWi-err}
\[%
\begin{array}
[c]{ccc}%
a_{n}=\dfrac{n}{e}+\dfrac{\ln(n)}{2e}+\dfrac{C}{e}+o(1), &  & (n!)^{1/n}%
=\dfrac{n}{e}+\dfrac{\ln(n)}{2e}+\dfrac{\ln(\sqrt{2\pi})}{e}+o(1)
\end{array}
\]
as $n\rightarrow\infty$, where
\[
C=e-1+\frac{\gamma}{2}+\frac{1}{2}%
%TCIMACRO{\dsum \limits_{k=1}^{\infty}}%
%BeginExpansion
{\displaystyle\sum\limits_{k=1}^{\infty}}
%EndExpansion
\frac{k-e\,a_{k}}{e\,k\,a_{k}}+%
%TCIMACRO{\dsum \limits_{k=1}^{\infty}}%
%BeginExpansion
{\displaystyle\sum\limits_{k=1}^{\infty}}
%EndExpansion
\left(  e\,a_{k+1}-e\,a_{k}-1-\frac{1}{2\,e\,a_{k}}\right)  =1.2905502....
\]
Further, $a_{n}-(n!)^{1/n}$ is strictly increasing and
\[
a_{n}-(n!)^{1/n}\leq\left(  C-\ln(\sqrt{2\pi})\right)  /e=0.136708...
\]
for all $n$. The constant is best possible. Putting $b_{n}=1/(e\,a_{n})$
yields the recursion $b_{n}=b_{n-1}\exp(-b_{n-1})$, for which an analogous
asymptotic expansion can be written.

The unique real zero $z_{n}$ of $%
%TCIMACRO{\tsum \nolimits_{k=0}^{n}}%
%BeginExpansion
{\textstyle\sum\nolimits_{k=0}^{n}}
%EndExpansion
z^{k}/k!$, where $n$ is odd, satisfies $\lim_{n\rightarrow\infty}%
z_{n}/n=W(e^{-1})=0.2784645427...=(3.5911214766...)^{-1}$ \cite{SZmyn-err,
HTntr-err}. The latter value appears in number theory \cite{Itrexp1-err,
Itrexp2-err, Itrexp3-err}, random graphs \cite{Itrexp4-err, Itrexp5-err,
Itrexp6-err}, ordered sets \cite{Itrexp7-err}, planetary dynamics
\cite{Itrexp8-err}, search theory \cite{Itrexp9-err, Itrexp10-err},
predator-prey models \cite{Itrexp11-err} and best-constant asymptotics
\cite{Itrexp12-err}.

From the study of minimum edge covers, given a complete bipartite graph, comes
$W(1)^{2}+2W(1)=1.4559380926...=2(0.7279690463...)$ \cite{HsslrW-err}. \ No
analogous formula is yet known for a related constant $0.55872...$%
\cite{Wast1-err}. \ See \cite{Zud6-err} for more on $W(1)$.

Also, $3^{-1}e^{-1/3}=0.2388437701...$ arises in \cite{HBoas-err} as a
consequence of the formula $-W(-3^{-1}e^{-1/3})=1/3$. Note that $-W(-x)$ is
the exponential generating function for rooted labeled trees and hence is
often called the \textit{tree function} \cite{CrHuJf-err}.

The equation $x\,e^{x}=1$ and numerous variations appear in \cite{MM75-err,
MM73-err, Pttl90, Pttl99, KS81-err, AFP98-err, GNS06-err, HB60-err}. For
example, let $S_{n}$ be the set of permutations on $\{1,2,\ldots,n\}$ and
$\sigma_{t}$ be a continuous-time random walk on $S_{n}$ starting from the
identity $I$ with steps chosen as follows: at times of a rate one Poisson
process, we perform a transposition of two elements chosen uniformly at
random, with replacement, from $\{1,2,\ldots,n\}$. Define $d(\sigma_{t})$ to
be the distance from $I$ at time $t$, that is, the minimum number of
transpositions required to return to $I$. For any fixed $c>0$,
\cite{BD06-err}
\[
d(\sigma_{c\,n/2})\sim\left(  1-%
%TCIMACRO{\dsum \limits_{k=1}^{\infty}}%
%BeginExpansion
{\displaystyle\sum\limits_{k=1}^{\infty}}
%EndExpansion
\frac1c\frac{k^{k-2}}{k!}(c\,e^{-c})^{k}\right)  n
\]
in probability as $n\rightarrow\infty$. The coefficient simplifies to $c/2$
for $c<1$ but is $<c/2$ otherwise. It is similar to the expansion
\[
1+\frac1cW(-c\,e^{-c})=1-%
%TCIMACRO{\dsum \limits_{k=1}^{\infty}}%
%BeginExpansion
{\displaystyle\sum\limits_{k=1}^{\infty}}
%EndExpansion
\frac1c\frac{k^{k-1}}{k!}(c\,e^{-c})^{k},
\]
differing only in the numerator exponent.

Consider the spread of a rumor though a population of $n$ individuals. Assume
that the number of ignorants is initially $\alpha\,n$ and that the number of
spreaders is $(1-\alpha)n$, where $0<\alpha<1$. A spreader-ignorant
interaction converts the ignorant into a spreader. When two spreaders
interact, they stop spreading the rumor and become stiflers. A
spreader-stifler interaction results in the spreader becoming a stifler. All
other types of interactions lead to no change. Let $\theta$ denote the
expected proportion of initial ignorants who never hear the rumor, then as
$\alpha$ decreases, $\theta$ increases (which is perhaps surprising!) and
\cite{Rmr1-err, Rmr2-err, Rmr3-err, Rmr4-err, Rmr5-err, Rmr6-err, Rmr7-err}
\[
0.2031878699...=\theta(1^{-})<\theta(\alpha)<\theta(0^{+}%
)=1/e=0.3678794411...
\]
as $n\rightarrow\infty$. The infimum of $\theta$ is the unique solution of the
equation $\ln(\theta)+2(1-\theta)=0$ satisfying $0<\theta<1$, that is,
$\theta=-W(-2e^{-2})/2$.

On the one hand, $\exp(x)=x$ has no real solution and $\sin(x)=x$ has no real
nonzero solution. On the other hand,\ $x=0.7390851332...$ appears in
connection with $\cos(x)=x$ \cite{Dot1-err, Dot2-err, Dot3-err}.

As with the divergent alternating factorial series on p. 425, we can assign
meaning to \cite{Wtsn-err}
\[%
%TCIMACRO{\dsum \limits_{n=0}^{\infty}}%
%BeginExpansion
{\displaystyle\sum\limits_{n=0}^{\infty}}
%EndExpansion
(-1)^{n}\,n^{n}=%
%TCIMACRO{\dsum \limits_{n=0}^{\infty}}%
%BeginExpansion
{\displaystyle\sum\limits_{n=0}^{\infty}}
%EndExpansion
\left(  \frac{(-1)^{n}n^{n}}{n!}%
%TCIMACRO{\dint \limits_{0}^{\infty}}%
%BeginExpansion
{\displaystyle\int\limits_{0}^{\infty}}
%EndExpansion
x^{n}e^{-x}dx\right)  =%
%TCIMACRO{\dint \limits_{0}^{\infty}}%
%BeginExpansion
{\displaystyle\int\limits_{0}^{\infty}}
%EndExpansion
\frac{e^{-x}}{1+W(x)}dx=0.7041699604...
\]
which also appears on p. 263. A variation is \cite{Gspr-err}
\begin{align*}%
%TCIMACRO{\dsum \limits_{n=1}^{\infty}}%
%BeginExpansion
{\displaystyle\sum\limits_{n=1}^{\infty}}
%EndExpansion
(-1)^{n+1}\,(2n)^{2n-1}  &  =%
%TCIMACRO{\dsum \limits_{n=1}^{\infty}}%
%BeginExpansion
{\displaystyle\sum\limits_{n=1}^{\infty}}
%EndExpansion
\left(  \frac{(-1)^{n+1}(2n)^{2n-1}}{(2n)!}%
%TCIMACRO{\dint \limits_{0}^{\infty}}%
%BeginExpansion
{\displaystyle\int\limits_{0}^{\infty}}
%EndExpansion
x^{2n}e^{-x}dx\right) \\
\  &  =%
%TCIMACRO{\dint \limits_{0}^{\infty}}%
%BeginExpansion
{\displaystyle\int\limits_{0}^{\infty}}
%EndExpansion
\ln\left(  \frac{x}{\sqrt{W(-i\,x)W(i\,x)}}\right)  e^{-x}dx=0.3233674316....
\end{align*}
which evidently is the same as \cite{SIAMh1-err, SIAMh2-err, SIAMh3-err}
\[%
%TCIMACRO{\dint \limits_{0}^{\infty}}%
%BeginExpansion
{\displaystyle\int\limits_{0}^{\infty}}
%EndExpansion
\frac{W(x)\cos(x)}{x(1+W(x))}dx=0.3233674316...
\]
although a rigorous proof is not yet known. Another variation is
\cite{Gspr-err}%
\begin{align*}%
%TCIMACRO{\dsum \limits_{n=1}^{\infty}}%
%BeginExpansion
{\displaystyle\sum\limits_{n=1}^{\infty}}
%EndExpansion
(-1)^{n+1}\,(2n-1)^{2n}  &  =\frac{i}{2}%
%TCIMACRO{\dint \limits_{0}^{\infty}}%
%BeginExpansion
{\displaystyle\int\limits_{0}^{\infty}}
%EndExpansion
\left(  \frac{W(-i\,x)}{\left[  1+W(-i\,x)\right]  ^{3}}-\frac{W(i\,x)}%
{\left[  1+W(i\,x)\right]  ^{3}}\right)  e^{-x}dx\\
&  =0.0111203007....
\end{align*}

The only two real solutions of the equation $x^{x-1}=x+1$ are
$0.4758608123...$ and $2.3983843827...$, which appear in \cite{Lim-err}.
Another example of striking coincidences between integrals and sums is
\cite{BPllard-err, BailBB-err}
\[%
%TCIMACRO{\dint \limits_{-\infty}^{\infty}}%
%BeginExpansion
{\displaystyle\int\limits_{-\infty}^{\infty}}
%EndExpansion
\frac{\sin(x)}{x}dx=%
%TCIMACRO{\dint \limits_{-\infty}^{\infty}}%
%BeginExpansion
{\displaystyle\int\limits_{-\infty}^{\infty}}
%EndExpansion
\frac{\sin(x)^{2}}{x^{2}}dx=\pi=%
%TCIMACRO{\dsum \limits_{n=-\infty}^{\infty}}%
%BeginExpansion
{\displaystyle\sum\limits_{n=-\infty}^{\infty}}
%EndExpansion
\frac{\sin(n)}{n}=%
%TCIMACRO{\dsum \limits_{n=-\infty}^{\infty}}%
%BeginExpansion
{\displaystyle\sum\limits_{n=-\infty}^{\infty}}
%EndExpansion
\frac{\sin(n)^{2}}{n^{2}};
\]
more surprises include \cite{OslrTsy-err}%
\[%
%TCIMACRO{\dint \limits_{0}^{1}}%
%BeginExpansion
{\displaystyle\int\limits_{0}^{1}}
%EndExpansion
t^{-x\,t}dt=\frac{1}{x}%
%TCIMACRO{\dsum \limits_{k=1}^{\infty}}%
%BeginExpansion
{\displaystyle\sum\limits_{k=1}^{\infty}}
%EndExpansion
\left(  \frac{x}{k}\right)  ^{k}=-%
%TCIMACRO{\dint \limits_{0}^{1}}%
%BeginExpansion
{\displaystyle\int\limits_{0}^{1}}
%EndExpansion
t^{-x\,t}\ln(t)dt
\]
for all real $x$. The integral \cite{Mathr3-err}
\begin{align*}
\lim_{N\rightarrow\infty}%
%TCIMACRO{\dint \limits_{1}^{2N}}%
%BeginExpansion
{\displaystyle\int\limits_{1}^{2N}}
%EndExpansion
e^{i\pi x}x^{1/x}dx  &  =0.0707760393...-(0.6840003894...)i\\
\  &  =-\frac{2}{\pi}i+\lim_{N\rightarrow\infty}%
%TCIMACRO{\dint \limits_{1}^{2N+1}}%
%BeginExpansion
{\displaystyle\int\limits_{1}^{2N+1}}
%EndExpansion
e^{i\pi x}x^{1/x}dx
\end{align*}
is analogous to the alternating series on p. 450 (since $(-1)^{x}=e^{i\pi x}$).

Exactly two branches of the Lambert function assume real values on $[-1/e,0)$,
one increasing ($W$) and the other decreasing ($\tilde{W}$). \ While
\cite{NIST-err}%
\[
W\left(  -\dfrac{\ln(2)}{2}\right)  =-\ln(2)>\tilde{W}\left(  -\dfrac{\ln
(2)}{2}\right)  =-\ln(4)
\]
and $W$, $\tilde{W}$ coincide at $-\ln(e)/e=-1/e$, we have
\[
\ln(\kappa)=W\left(  -\dfrac{\ln(3)}{3}\right)  =-0.90747304...>\tilde
{W}\left(  -\dfrac{\ln(3)}{3}\right)  =-\ln(3),
\]
i.e., $x^{x}=3^{-1/3}$ has real solutions $x=\kappa=0.40354267...$ and
$x=1/3$. \ What additionally can be said about $\kappa$?

\textbf{6.12. Conway's Constant.} A\ \textquotedblleft
biochemistry\textquotedblright\ based on Conway's \textquotedblleft
chemistry\textquotedblright\ appears in \cite{OMrtn-err}; some
\textquotedblleft stuttering\textquotedblright\ variants are found in
\cite{Stutter1-err, Stutter2-err}.

\textbf{7.1. Bloch-Landau Constants.} In the definitions of the sets $F$ and
$G$, the functions $f$ need only be analytic on the open unit disk $D$ (in
addition to satisfying $f(0)=0$, $f^{\prime}(0)=1$). On the one hand, the
weakened hypothesis doesn't affect the values of $B$, $L$, $K$ or $A$; on the
other hand, the weakening is essential for the existence of $f\in G$ such that
$m(f)=M$. \ We now know that $0.57088586<K\leq0.6563937$ \cite{Xiong-err,
CrlOC-err, Sknnr-err}.

The bounds $0.62\pi<A<0.7728\pi$ were improved by several authors, although
they studied the quantity $\tilde{A}=\pi-A$ instead (the omitted area
constant). Barnard \&\ Lewis \cite{BL-err} demonstrated that $\tilde{A}%
\leq0.31\pi$. Barnard \&\ Pearce \cite{BP-err} established that $\tilde{A}%
\geq0.240005\pi$, but Banjai \&\ Trefethen \cite{BT-err} subsequently computed
that $\tilde{A}=(0.2385813248...)\pi$. It is believed that the earlier
estimate was slightly in error. See \cite{Lew-err, Wan-err, BPC-err, BPS-err}
for related problems.

The spherical analog of Bloch's constant $B$, corresponding to meromorphic
functions $f$ mapping $D$ to the Riemann sphere, was recently determined by
Bonk \&\ Eremenko \cite{BE-err}. This constant turns out to be $\arccos
(1/3)=1.2309594173...$. A proof as such gives us hope that someday the planar
Bloch-Landau constants will also be exactly known \cite{Rttngr1-err,
Rttngr2-err}.

More relevant material is found in \cite{Fi13-err, Fi9-err}.

\textbf{7.2. Masser-Gramain Constant.} It is now known that $1.819776<\delta
<1.819833$, overturning Gramain's conjecture \cite{MelqNZ-err}. Suppose $f(z)
$ is an entire function such that $f^{(k)}(n)$ is an integer for each
nonnegative integer $n$, for each integer $0\leq k\leq s-1$. (We have
discussed only the case $s=1$.) The best constant $\theta_{s}>0$ for which
\[%
\begin{tabular}
[c]{lll}%
$\operatorname*{limsup}\limits_{r\rightarrow\infty}\dfrac{\ln(M_{r})}%
{r}<\theta_{s}$ &  & implies $f$ is a polynomial
\end{tabular}
\]
was proved by Bundschuh \&\ Zudilin \cite{BuZu-err}, building on Gel'fond
\cite{Glfd-err} and Selberg \cite{Slbg-err}, to satisfy
\[
s\cdot\frac{\pi}{3}\geq\theta_{s}>\left\{
\begin{array}
[c]{ccc}%
0.994077... &  & \text{if }s=2,\\
1.339905... &  & \text{if }s=3,\\
1.674474... &  & \text{if }s=4.
\end{array}
\right.
\]
(Actually they proved much more.) Can a Gaussian integer-valued analog of
these integer-valued results be found?

\textbf{7.3. Whittaker-Goncharov Constants}. \ The lower bound $0.73775075<W$,
due to Waldvogel (using Goncharov polynomials), appears only in Varga's
survey; it is not mentioned in \cite{Waldvgl-err}. \ Minimum modulus
zero-finding techniques provide the upper bound $W\leq0.7377507574...$%
(correcting $<$). \ Both bounds are non-rigorous. \ The \textquotedblleft
third constant\textquotedblright\ involves zero-free disks for the
Rogers-Szeg\"{o} polynomials:%
\[%
\begin{array}
[c]{ccc}%
G_{n+1}(z,q)=(1+z)G_{n}(z,q)-\left(  1-q^{n}\right)  G_{n-1}(z,q), &  &
n\geq0,
\end{array}
\]%
\[%
\begin{array}
[c]{ccc}%
G_{-1}(z,q)\equiv0, &  & G_{0}(z,q)\equiv1
\end{array}
\]
where $q\in\mathbb{C}$. \ Let%
\[
r_{n}=\inf\left\{  \left\vert z\right\vert :G_{n}(z,q)=0\text{ and }\left\vert
q\right\vert =1\right\}  ,
\]
then numerical data suggests \cite{Waldvgl-err}%
\[
r_{n}=\left(  3-2\sqrt{2}\right)  +\left(  0.3833...\right)  n^{-2/3}+O\left(
n^{-4/3}\right)
\]
as $n\rightarrow\infty$. \ A\ proof remains open. \ Such asymptotics are
relevant to study of the partial theta function $\sum_{j=0}^{\infty
}q^{j(j-1)/2}z^{j}$ and associated Pad\'{e} approximant convergence properties.

\textbf{7.4. John Constant.} \ Consider analytic functions $f$ defined on the
unit disk $D$ that satisfy $f(0)=0$, $f^{\prime}(0)=1$ and%
\[
\ell\leq\left\vert \frac{z\,f^{\prime}(z)}{f(z)}\right\vert \leq L
\]
at all points $z\in D$. \ The ratio plays the same role as $\left\vert
f^{\prime}(z)\right\vert $ did originally. \ What is the largest number
$\delta$ such that $L/\ell\leq\delta$ implies that $f$ is univalent (on $D $)?
\ Kim \&\ Sugawa \cite{KmSgw1-err, KmSgw2-err} proved that $\exp
(7\pi/25)<\delta<\exp(5\pi/7)$ and stated that tighter bounds are possible.
\ No Gevirtz-like conjecture governing an exact expression for $\delta$ has
yet been proposed.

\textbf{7.5. Hayman Constants.} New bounds \cite{Wng1-err, Wng2-err, Wng3-err,
Shn0-err, Wng4-err, Wng5-err, Wng6-err} for the Hayman-Korenblum constant
$c(2)$ are $0.28185$ and $0.67789$. An update on the Hayman-Wu constant
appears in \cite{Roh-err}.

\textbf{7.6. Littlewood-Clunie-Pommerenke Constants.} The lower limit of
summation in the definition of $S_{2}$ should be $n=0$ rather than $n=1$, that
is, the coefficient $b_{0}$ need not be zero. We have sharp bounds
$|b_{1}|\leq1$, $|b_{2}|\leq2/3$, $|b_{3}|\leq1/2+e^{-6}$ \cite{GaSc-err}. The
bounds on $\gamma_{k}$ due to Clunie \& Pommerenke should be $0.509$ and
$0.83$ \cite{Pmmr-err}; Carleson \&\ Jones' improvement was nonrigorous. While
$0.83=1-0.17$ remains the best established upper bound, the lower bound has
been increased to $0.54=1-0.46$ \cite{MaPm-err, GrPm-err, HdSh-err}. Numerical
evidence for both the Carleson-Jones conjecture and Brennan's conjecture was
found by Kraetzer \cite{Kra-err}. Theoretical evidence supporting the latter
appears in \cite{BVZ-err}, but a complete proof remains undiscovered.
A\ claimed verification that $\alpha=1-\gamma$ \cite{XZW1-err, XZW2-err},
based on \cite{XZW3-err, XZW4-err, XZW5-err}, is unfounded because
\cite{XZW5-err} does not seem to exist.

\textbf{7.7. Riesz-Kolmogorov Constants}. \ The constant $C_{1}$ appears
recently, for example, in \cite{Oskwski-err}.

\textbf{7.8. Gr\"{o}tzsch Ring Constants}. \ The phrase \textquotedblleft
planar ring\textquotedblright\ appearing in the first sentence should be
\textquotedblleft planar region\textquotedblright.

\textbf{8.1. Geometric Probability Constants.} Just as the ratio of a
semicircle to its radius is always $\pi$, the ratio of the latus rectum arc of
any parabola to its semi latus rectum is \cite{Love-err}%
\[
\sqrt{2}+\ln\left(  1+\sqrt{2}\right)  =2.2955871493...=2(1.1477935746...)
\]
Is it mere coincidence that this constant is so closely related to the
quantity $\delta(2)$? Just as the ratio of the area of a circle to its radius
squared is always $\pi$, the ratio of the area of the latus rectum segment of
any equilateral hyperbola to its semi-axis squared is \cite{ReeSndw-err}%
\[
\sqrt{2}-\ln\left(  1+\sqrt{2}\right)  =0.5328399753....
\]
The similarity in formulas is striking:\ length of one conic section
(universal parabolic constant) versus area of another (universal equilateral
hyperbolic constant).

Consider the logarithm $\Lambda$ of the distance between two independent
uniformly distributed points in the unit square. \ The constant%
\[
\exp\left(  \operatorname{E}(\Lambda)\right)  =\exp\left(  \frac{-25+4\pi
+4\ln(2)}{12}\right)  =0.4470491559...=2(0.2235245779...)
\]
appears in calculations of electrical inductance of a long solitary wire with
small rectangular cross section \cite{Induc1-err, Induc2-err, Induc3-err,
Induc4-err}. \ If the wire is fairly short, then more complicated formulas
apply \cite{Induc5-err, Induc6-err, Induc7-err}. \ The constants%
\[%
\begin{array}
[c]{lll}%
e^{-1/4}=0.7788007830..., &  & e^{-3/2}=0.2231301601...
\end{array}
\]
appear instead for cross sections in the form of a disk and an interval, respectively.

The expected distance between two random points on different sides of the unit
square is \cite{BBVW-err}
\[
\frac{2+\sqrt{2}+5\ln\left(  1+\sqrt{2}\right)  }{9}=0.8690090552...
\]
and the expected distance between two random points on different faces of the
unit cube is
\[
\frac{4+17\sqrt{2}-6\sqrt{3}-7\pi+21\ln\left(  1+\sqrt{2}\right)
+21\ln\left(  7+4\sqrt{3}\right)  }{75}=0.9263900551...
\]
See \cite{BoxInt-err, JPhilip-err} for expressions involving $\delta(4)$,
$\Delta(4)$ and $\Delta(5)$. Asymptotics of $\delta_{p}(n)$ and $\Delta
_{p}(n)$ in the $\ell_{p}$ norm as $n\rightarrow\infty$, for fixed
\thinspace$p>0$, are found in \cite{SStnrb-err}. See \cite{ZhngP1-err,
ZhngP2-err, StwrtZ-err, UBsl1-err, UBsl2-err, ZhngP3-err, AhmdP1-err,
AhmdP2-err} for results not in a square, but in an equilateral triangle or
regular hexagon. \ The constant $2\sqrt{\pi}M$ appears in \cite{BiHeTk-err}.
\ Also, the convex hull of random point sets in the unit disk (rather than the
unit square) is mentioned in \cite{Fi10-err}, and properties of random
triangles are extensively covered in \cite{Fi29-err}.

\textbf{8.2. Circular Coverage Constants.} The coefficient of $x^{16}$ in the
minimal polynomial for $r(6)$ should be $-33449976$. \ Fejes T\'{o}th
\cite{GFjsT-err} proved the conjectured formula for $r(N)$ when $8\leq
N\leq10$. \ Here is a variation of the elementary problems at the end. Imagine
two overlapping disks, each of radius 1. If the area $A$ of the intersection
is equal to one-third the area of the union, then clearly $A=\pi/2$. The
distance $w$ between the centers of the two circles is $w=0.8079455065...$,
that is, the unique root of the equation
\[
2\arccos\left(  \frac{w}{2}\right)  -\frac{1}{2}w\sqrt{4-w^{2}}=\frac{\pi}{2}%
\]
in the interval $[0,2]$. If \textquotedblleft one-third\textquotedblright\ is
replaced by \textquotedblleft one-half\textquotedblright, then $\pi/2$ is
replaced by $2\pi/3$ and Mrs. Miniver's constant $0.5298641692...$ emerges
instead. \ Let $z=1.9056957293...$ be the unique solution of $\sin
(z)-z\cos(z)=\pi/2$ satisfying $\left\vert z\right\vert <4$, then
$2\cos(z/2)=1.1587284730...$ is the grazing goat constant \cite{Ullgoat-err,
Hffgoat-err}.

\textbf{8.3. Universal Coverage Constants.} Elekes \cite{Elekes-err} improved
the lower bound for $\mu$ to $0.8271$ and Brass \& Sharifi \cite{BrShrf-err}
improved this further to $0.832$. \ Computer methods were used in the latter
to estimate the smallest possible convex hull of a circle, equilateral
triangle and regular pentagon, each of diameter $1$. \ Hansen evidently made
use of reflections in his convex cover, as did Duff in his nonconvex cover;
Gibbs \cite{PGibbs1-err, PGibbs2-err} claimed a reduced upper bound of
$0.844112$ for the convex case, using reflections. \ It would seem that
Sprague's upper bound remains the best known for displacements, strictly
speaking.\ \ Two additional references for translation covers include
\cite{Rnni-err, Duff-err}.

\textbf{8.4. Moser's Worm Constant.} Coulton \&\ Movshovich \cite{CtMv-err}
proved Besicovitch's conjecture that every worm of unit length can be covered
by an equilateral triangular region of area $7\sqrt{3}/27$. The upper bound
for $\mu$ was decreased \cite{WWg-err} to $0.270912$; the lower bound for
$\mu$ was increased \cite{SrKh-err, KhPaSr-err} to $0.232239$. \ New bounds
$0.096694<\mu^{\prime}<0.112126$ appear in \cite{Wtz0-err}. Relevant progress
is described in \cite{Wtz1-err, JPW-err, Wcha-err, Wtz2-err}. We mention, in
Figure 8.3, that the quantity $x=\sec(\varphi)=1.0435901095...$ is algebraic
of degree six \cite{Zlglr-err, RAlxnd-err}:
\[
3x^{6}+36x^{4}+16x^{2}-64=0
\]
and wonder if this is linked to Figure 8.7 and the Reuleaux triangle of width
$1.5449417003...$ (also algebraic of degree six \cite{Sallee-err}). The latter
is the planar set of maximal constant width that avoids all vertices of the
integer square lattice. \ See \cite{Ghomi-err, GhmW1-err, GhmW2-err, CGLQ-err}
for discussion of Zalgaller's work in $3$-space and \cite{Wa-liaf-err,
KL-liaf-err} for more on \textquotedblleft lost in a forest\textquotedblright\ problems.

\textbf{8.5. Traveling Salesman Constants.} Let $\delta=\left(  \sqrt{2}%
+\ln(1+\sqrt{2})\right)  /6$, the average distance from a random point in the
unit square to its center (page 479). If we identify edges of the unit square
(wrapping around to form a torus), then $\operatorname*{E}(L_{2}(n))/\delta=n$
for $n=2$, $3$ but $\operatorname*{E}(L_{2}(4))/\delta\approx3.609...$. A
closed-form expression for the latter would be good to see \cite{Hateful-err}.
The best upper bound on $\beta_{2}^{\prime}$ is now $0.6321$
\cite{Steinberg-err}; more numerical estimates of $\beta_{2}=0.7124...$ appear
in \cite{Cmprbl-err}.

The random links TSP\ $\beta=2.0415...$ possesses an alternative formulation
\cite{Wast1-err, Wast2-err}: let $y>0$ be defined as an implicit function of
$x$ via the equation%
\[
\left(  1+\frac{x}{2}\right)  e^{-x}+\left(  1+\frac{y}{2}\right)  e^{-y}=1,
\]
then%
\[
\beta=\frac{1}{2}%
%TCIMACRO{\dint \limits_{0}^{\infty}}%
%BeginExpansion
{\displaystyle\int\limits_{0}^{\infty}}
%EndExpansion
y(x)\,dx=2.0415481864...=2(1.0207740932...).
\]
This constant is the same if the lengths are distributed according to
Exponential(1) rather than Uniform[0,1]. \ If instead lengths are equal to the
square roots of exponential variables, the resulting constant is
$1.2851537533...=(1.8174818677...)/\sqrt{2}=(0.7250703609...)\sqrt{\pi}.$

Other proofs that the minimum matching $\beta=\pi^{2}/12$ are known; see
\cite{Wast3-err}. \ If (as in the preceding) lengths are equal to the square
roots of exponential variables, the resulting constant is
$0.5717590495...=(1.1435180991...)/2=(0.3225805000...)\sqrt{\pi}$, recovering
M\'{e}zard \&\ Parisi's calculation \cite{MzPrs-err}. \ An integral
equation-based formula for the latter is \cite{Wast1-err, Wast4-err}%
\[%
\begin{array}
[c]{ccccc}%
\beta=2%
%TCIMACRO{\dint \limits_{0}^{\infty}}%
%BeginExpansion
{\displaystyle\int\limits_{0}^{\infty}}
%EndExpansion%
%TCIMACRO{\dint \limits_{-y}^{\infty}}%
%BeginExpansion
{\displaystyle\int\limits_{-y}^{\infty}}
%EndExpansion
(x+y)f(x)f(y)\,dx\,dy &  & \text{where} &  & f(x)=\exp\left(  -2%
%TCIMACRO{\dint \limits_{0}^{\infty}}%
%BeginExpansion
{\displaystyle\int\limits_{0}^{\infty}}
%EndExpansion
t\,f(t-x)\,dt\right)  .
\end{array}
\]

The cavity method is applied in \cite{ZdMz-err} to matchings on sparse random
graphs. Also, for the cylinder graph $P_{n}\times C_{k}$ on $(n+1)k$ vertices
with independent Uniform $[0,1]$ random edge-lengths, we have
\[
\lim_{n\rightarrow\infty}\frac{1}{n}L_{\text{MST}}(P_{n}\times C_{k}%
)=\gamma(k)
\]
almost surely, where $k$ is fixed and \cite{HutLew-err}
\begin{align*}
\gamma(2)  &  =-%
%TCIMACRO{\dint \limits_{0}^{1}}%
%BeginExpansion
{\displaystyle\int\limits_{0}^{1}}
%EndExpansion
\frac{(x-1)^{2}(2x^{3}-3x^{2}+2)}{x^{4}-2x^{3}+x^{2}-1}dx\\
\  &  =2-\frac{1}{\sqrt{5}}\ln\left(  \frac{\sqrt{5}-1}{\sqrt{5}+1}\right)
-\frac{\pi}{\sqrt{3}}=0.6166095767...,
\end{align*}%
\begin{align*}
\gamma(3)  &  =-%
%TCIMACRO{\dint \limits_{0}^{1}}%
%BeginExpansion
{\displaystyle\int\limits_{0}^{1}}
%EndExpansion
\frac{(x-1)^{3}(3x^{8}-11x^{7}+13x^{6}+x^{5}-18x^{4}+14x^{3}+3x^{2}%
-3x-3)}{x^{10}-5x^{9}+10x^{8}-10x^{7}+x^{6}+11x^{5}-11x^{4}+2x^{3}+x^{2}%
-1}dx\\
\  &  =0.8408530104...
\end{align*}
and $\gamma(4)=1.09178...$.

\textbf{8.6. Steiner Tree Constants.} Doubt has been raised \cite{IKMS-err,
IvTz-err} about the validity of the proof by Du \&\ Hwang of the Gilbert \&
Pollak conjecture.

\textbf{8.7. Hermite's Constants.} A lattice $\Lambda$ in $\mathbb{R}^{n}$
consists of all integer linear combinations of a set of basis vectors
$\{e_{j}\}_{j=1}^{n}$ for $\mathbb{R}^{n}$. If the fundamental parallelepiped
determined by $\{e_{j}\}_{j=1}^{n}$ has Lebesgue measure $1$, then $\Lambda$
is said to be of unit volume. The constants $\gamma_{n}$ can be defined via an
optimization problem
\[
\gamma_{n}=\max_{\substack{\text{unit volume} \\\text{lattices }\Lambda
}}\,\min_{\substack{x\in\Lambda, \\x\neq0 }}\,\left\|  x\right\|  ^{2}%
\]
and are listed in Table 8.10. The precise value of the next constant
$2\leq\gamma_{9}<2.1327$ remains open \cite{CoSl-err, CoEl-err, NeSl-err},
although Cohn \&\ Kumar \cite{CoKu1-err, CoKu2-err} have recently proved that
$\gamma_{24}=4$. A classical theorem \cite{Wats-err, MiHu-err, BeMa-err}
provides that $\gamma_{n}^{n}$ is rational for all $n$. It is not known if the
sequence $\gamma_{1}$, $\gamma_{2}$, $\gamma_{3}$, $\ldots$ is strictly
increasing, or if the ratio $\gamma_{n}/n$ tends to a limit as $n\rightarrow
\infty$. See also \cite{Fi11-err, Bacher-err}.

Table 8.10. \textit{Hermite's constants }$\gamma_{n}$%
\[%
\begin{tabular}
[c]{|l|l|l|l|l|l|}\hline
$n$ & Exact & Decimal & $n$ & Exact & Decimal\\\hline
1 & $1$ & $1$ & 5 & $8^{1/5}$ & $1.5157165665...$\\\hline
2 & $(4/3)^{1/2}$ & $1.1547005383...$ & 6 & $(64/3)^{1/6}$ & $1.6653663553...
$\\\hline
3 & $2^{1/3}$ & $1.2599210498...$ & 7 & $64^{1/7}$ & $1.8114473285...$\\\hline
4 & $4^{1/4}$ & $1.4142135623...$ & 8 & $2$ & $2$\\\hline
\end{tabular}
\]

An arbitrary packing of the plane with disks is called compact if every disk
$D$ is tangent to a sequence of disks $D_{1}$, $D_{2}$, $\ldots$, $D_{n}$ such
that $D_{i}$ is tangent to $D_{i+1}$ for $i=1$, $2$, $\ldots$, $n$ with
$D_{n+1}=D_{1}$. If we pack the plane using disks of radius $1$, then the only
possible compact packing is the hexagonal lattice packing with density
$\pi/\sqrt{12}$. If we pack the plane using disks of radius $1$ and $r<1$
(disks of both sizes must be used), then there are precisely nine values of
$r$ for which a compact packing exists. See Table 8.11. For seven of these
nine values, it is known that the densest packing is a compact packing; the
same is expected to be true for the remaining two values \cite{Kndy1-err,
Kndy2-err, Kndy3-err}. \ Breaking news: the latter assertion is now proved
\cite{BdrFnq-err}.\pagebreak

Table 8.11. \textit{All nine values of }$r<1$\textit{\ which allow compact
packings}

\begin{center}%
\begin{tabular}
[c]{|l|l|}\hline
Exact (expression or minimal polynomial) & Decimal\\\hline
$5-2\sqrt{6}$ & $0.1010205144...$\\\hline
$(2\sqrt{3}-3)/3$ & $0.1547005383...$\\\hline
$(\sqrt{17}-3)/4$ & $0.2807764064...$\\\hline
$x^{4}-28x^{3}-10x^{2}+4x+1$ & $0.3491981862...$\\\hline
$9x^{4}-12x^{3}-26x^{2}-12x+9$ & $0.3861061048...$\\\hline
$\sqrt{2}-1$ & $0.4142135623...$\\\hline
$8x^{3}+3x^{2}-2x-1$ & $0.5332964166...$\\\hline
$x^{8}-8x^{7}-44x^{6}-232x^{5}-482x^{4}-24x^{3}+388x^{2}-120x+9$ &
$0.5451510421...$\\\hline
$x^{4}-10x^{2}-8x+9$ & $0.6375559772...$\\\hline
\end{tabular}

\end{center}

There is space to only mention the circle-packing rigidity constants $s_{n}$
\cite{PHR-err}, their limiting behavior:
\[
\lim_{n\rightarrow\infty}n\,s_{n}=\frac{2^{4/3}}{3}\frac{\Gamma(1/3)^{2}%
}{\Gamma(2/3)}=4.4516506980...
\]
and their connection with conformal mappings. \ Also, the tetrahedral analog
of Kepler's sphere packing density is \textit{possibly}
$4000/4671=0.856347...$ \cite{Tetra1-err, Tetra2-err, Tetra3-err}, but a proof
would likely be exceedingly hard. \ A\ \textquotedblleft second-order
constant\textquotedblright\ in Hales' theorem is now better understood
\cite{NSchrf-err}.

In a recent breakthrough, Viazovska \cite{Viaz1-err, Viaz2-err} determined
that $\Delta_{8}=\pi^{4}/384$ (as expected); follow-on work \cite{Viaz3-err}
gave $\Delta_{24}=\pi^{12}/479001600$.

\textbf{8.8. Tammes' Constants.} Recent conjectures give \cite{BrHdSf-err}%
\[
\lambda=3\left(  \frac{8\pi}{\sqrt{3}}\right)  ^{1/2}\zeta\left(  -\frac{1}%
{2}\right)  \beta\left(  -\frac{1}{2}\right)  =-0.3992556250...
\]
(data fitting earlier predicted $\lambda\approx-0.401$) and%
\[
\mu=\ln(2)+\frac{1}{4}\ln\left(  \frac{2}{3}\right)  +\frac{3}{2}\ln\left(
\frac{\sqrt{\pi}}{\Gamma(1/3)}\right)  =-0.0278026524...=\frac
{-0.0556053049...}{2}.
\]
(improving on $\mu\approx-0.026$). \ Let nonzero $\alpha$ satisfy
$-4<\alpha<2$. The asymptotics for $\alpha=\pm1$ are subsumed by
\[
E(\alpha,N)=\left\{
\begin{array}
[c]{lll}%
\frac{2^{\alpha}}{2+\alpha}N^{2}+3\left(  \frac{8\pi}{\sqrt{3}}\right)
^{\alpha/2}\zeta\left(  -\frac{\alpha}{2}\right)  \beta\left(  -\frac{\alpha
}{2}\right)  N^{1-\alpha/2}+o\left(  N^{1-\alpha/2}\right)  &  & \text{if
}\alpha\neq2,\\
\frac{1}{8}N^{2}\ln(N)+\frac{c}{2}N^{2}+O(1) &  & \text{if }\alpha=2
\end{array}
\right.
\]
as $N\rightarrow\infty$, where
\[
c=\frac{1}{4}\left[  2\gamma+\ln\left(  8/3\right)  +3\ln(\pi)-6\ln\left(
\Gamma\left(  1/3\right)  \right)  \right]  =-0.0857684103...
\]
and the expression for $c$, originally given in terms of generalized Stieltjes
constants, follows from \cite{Mlmstn2-err, Connon-err}.

Consider the problem of covering a sphere by $N$ congruent circles (spherical
caps) so that the angular radius of the circles will be minimal. For
$N=8,9,11$ the conjectured best covering configurations remain unproven
\cite{Tmms1-err, Tmms2-err, Tmms3-err, Tmms4-err, Tmms5-err}.

\textbf{8.9. Hyperbolic Volume Constants}. \ Exponentially improved lower
bounds for $f(n)$ are now known \cite{Glzyrn-err}. \ Let $H(n)=\xi_{n}%
/\eta_{n}$ (due to Smith)\ and $K(n)=(n+1)^{(n-1)/2}$ (due to Glazyrin). \ We
have $f(n)\geq K(n)$ always and%
\[
\lim_{n\rightarrow\infty}\left(  \frac{K(n)}{E(n)}\right)  ^{1/n}=\frac{e}%
{2}=1.3591409142...>1.2615225101...=\sqrt{\frac{e}{2}}c=\lim_{n\rightarrow
\infty}\left(  \frac{H(n)}{E(n)}\right)  ^{1/n}%
\]
where $E(n)=2^{n}(n+1)^{-(n+1)/2}n!$ (simple bound used for comparison).
\ Alternatively,
\[
\lim_{n\rightarrow\infty}\frac{K(n)^{1/n}}{\sqrt{n}}=1>0.9281763921...=\sqrt
{\frac{2}{e}}c=\lim_{n\rightarrow\infty}\frac{H(n)^{1/n}}{\sqrt{n}}.
\]
For $n>2$, a \textit{dissection} of the $n$-cube need not be a triangulation;
the term \textquotedblleft simplexity\textquotedblright\ can be ambiguous in
the literature. \ See also \cite{BlisSu-err}.

\textbf{8.10. Reuleaux Triangle Constants.} In our earlier entry [8.4], we ask
about the connection between two relevant algebraic quantities
\cite{Zlglr-err, Sallee-err}, both zeroes of sextic polynomials.

\textbf{8.11. Beam Detection Constants}. The shortest \textit{opaque set} or
\textit{barrier} for the circle remains unknown; likewise for the square and
equilateral triangle \cite{AGbeam-err, DJPbeam-err, KMOPbeam-err,
DJTbeam-err}. \ See \cite{Ghomi-err, CGLQ-err} for discussion of Zalgaller's
work in $3$-space.

\textbf{8.12. Moving Sofa Constant}. \ The passage of an $\ell\times w$
rectangular piano around a right-angled corner in a hallway of before-width
$u$ and after-width $v$ can be determined by checking the sign of a certain
homogenous sextic polynomial in $\ell,r,u,v$, where $\ell>u\geq v>w$
\cite{YZpnmv-err}. \ Progress toward confirming Gerver's conjecture appears in
\cite{KlsRmk-err}.

\textbf{8.13. Calabi's Triangle Constant.} See \cite{Calabi-err} for details
underlying the main result.

\textbf{8.14. DeVicci's Tesseract Constant}. DeVicci's manuscript
\cite{KPDVS-err} is available online, but unfortunately not\ Ligocki
\&\ Huber's numerical experiments \cite{LgkHbr-err}.

\textbf{8.15. Graham's Hexagon Constant.} Bieri \cite{Bieri-err} partially
anticipated Graham's result. A nice presentation of Reinhardt's isodiametric
theorem is found in \cite{Msngf-err}.

\textbf{8.16. Heilbronn Triangle Constants.} Another vaguely-related problem
involves the maximum $M$ and minimum $m$ of the $\tbinom{n}{2}$ pairwise
distances between $n$ distinct points in $\mathbb{R}^{2}$. What configuration
of $n$ points gives the smallest possible ratio $r_{n}=M/m$? It is known that
\cite{DCntrl1-err, EFrdm1-err}
\[%
\begin{array}
[c]{ccccccccccc}%
r_{3}=1, &  & r_{4}=\sqrt{2}, &  & r_{5}=\varphi, &  & r_{6}=(\varphi\sqrt
{5})^{1/2}, &  & r_{7}=2, &  & r_{8}=\psi
\end{array}
\]
where $\varphi$ is the Golden mean and $\psi=\csc(\pi/14)/2$ has minimal
polynomial $\psi^{3}-2\psi^{2}-\psi+1$. We also have $r_{12}=$ $\sqrt
{5+2\sqrt{3}}$ and an asymptotic result of Thue's \cite{PErdos-err,
CrftFG-err}:
\[
\lim_{n\rightarrow\infty}n^{-1/2}r_{n}=\sqrt{\frac{2\sqrt{3}}{\pi}}.
\]
Erd\H{o}s wrote that the corresponding value of $\lim_{n\rightarrow\infty
}n^{-1/3}r_{n}$ for point sets in $\mathbb{R}^{3}$ is not known. Cantrell
\cite{DCntrl2-err, EFrdm2-err} wrote that it should be $\sqrt[3]{3\sqrt{2}%
/\pi}$, that is, the cube root of the reciprocal of the Kepler packing density
(proved by Hales).

\textbf{8.17. Kakeya-Besicovitch Constants}. \ Reversal of line segments in
higher dimensional regions is the subject of \cite{JJKM-err}.

\textbf{8.18. Rectilinear Crossing Constant.} We now know $\bar{\nu}(K_{n})$
for all $n\leq30$ except $n\in\{28,29\}$ -- see Table 8.12 -- and consequently
$0.379972<\rho<0.380488$ \cite{UVW1-err, UVW2-err, UVW3-err, UVW4-err,
UVW5-err, UVW6-err, UVW7-err, UVW8-err, UVW9-err, UVW10-err, UVW11-err,
UVW12-err}.

Table 8.12. \textit{Values of }$\bar{\nu}(K_{n})$,\textit{\ }$n>12$%

\[%
\begin{tabular}
[c]{|c|c|c|c|c|c|c|c|c|}\hline
$n$ & $13$ & $14$ & $15$ & $16$ & $17$ & $18$ & $19$ & $20$\\\hline
$\bar{\nu}(K_{n})$ & $229$ & $324$ & $447$ & $603$ & $798$ & $1029$ & $1318$ &
$1657$\\\hline
\end{tabular}
\]%
\[%
\begin{tabular}
[c]{|c|c|c|c|c|c|c|c|c|}\hline
$n$ & $21$ & $22$ & $23$ & $24$ & $25$ & $26$ & $27$ & $30$\\\hline
$\bar{\nu}(K_{n})$ & $2055$ & $2528$ & $3077$ & $3699$ & $4430$ & $5250$ &
$6180$ & $9726$\\\hline
\end{tabular}
\]
The validity of Guy's conjectured expression $Z(n)$ (more appropriately named
after Hill \cite{HryHil-err, BnkWs-err}) remains open, although the ratio
$\nu(K_{n})/Z(n)$ is asymptotically $\geq0.8594$ as $n\rightarrow\infty$
\cite{WVU1-err, WVU2-err, WVU3-err, WVU4-err}. \ It is well-known that
$q(R)=25/36\approx0.694$ when $R$ is a rectangle. \ If instead the four points
are bivariate normally distributed, then
\[
q=3\left(  1-2\operatorname{arcsec}(3)/\pi\right)  \approx0.649<2/3.
\]
The proof uses expectation formulas for the number of vertices \cite{Efrn-err,
MCRF-err} and for order statistics \cite{HDvd-err, Fi38-err}.

\textbf{8.19. Circumradius-Inradius Constants.} The phrase \textquotedblleft%
$Z$-admissible\textquotedblright\ in the caption of Figure 8.22 should be
replaced by \textquotedblleft$Z$-allowable\textquotedblright.

\textbf{8.20. Apollonian Packing Constant.} The packing exponent is calculated
to high precision $1.3056867280...$ in \cite{ZQBF2-err}; it appears in a
setting \cite{Fi35-err} which vastly generalizes the circular configurations
portrayed in Figure 8.23.

\textbf{8.21. Rendezvous Constants.} It is now known \cite{Plhmr-err} that
$r(T)\leq R_{2}\leq S_{2}\leq0.678442$; proof that $S_{2}=R_{2}%
=r(T)=0.6675277360...$ remains open.

\textbf{Table of Constants.} The formula corresponding to $0.8427659133...$ is
$(12\ln(2))/\pi^{2}$ and the formula corresponding to $0.8472130848...$ is
$M/\sqrt{2}$.

\section{Second Volume}

\ \ \ \ \textbf{1.5. Multiples and Divisors.} \ With regard to practical
numbers, estimates \cite{Wngrt1-err} for the constant $\kappa$ were improved
\cite{Wngrt2-err} to $1.336073<\kappa<1.336077$.

\textbf{1.7. Unitarism and Infinitarism.} \ The factor $1/N^{2}$ preceding the
summation of $d_{\infty}(n)$ over $1\leq n\leq N$ should be removed.

\textbf{1.22. Central Binomial Coefficients.} A\ new series evaluation
\cite{CLvN-err, SfNim-err}%
\[%
%TCIMACRO{\dsum \limits_{n=0}^{\infty}}%
%BeginExpansion
{\displaystyle\sum\limits_{n=0}^{\infty}}
%EndExpansion
\tbinom{2n}{n}\frac{1}{2^{3n}(2n+1)^{3}}=\frac{1}{\sqrt{2}}\left(  \frac{i}%
{2}\left[  \operatorname*{Li}\nolimits_{3}(1-i)-\operatorname*{Li}%
\nolimits_{3}(1+i)\right]  +\frac{\pi^{3}}{192}\right)
\]
evokes a question: can the expression within square brackets be written in
terms of other mathematical constants? \ Follow-on work appears in
\cite{Dsrdy-err}.

\textbf{1.30. Signum Equations and Extremal Coefficients.} \ The open issue
raised in the final sentence was addressed by Sudler \cite{Sudler-err} (with
numerical work by A. Hurwitz) and Wright \cite{Wrgt-err}:%
\[%
\begin{array}
[c]{ccccc}%
\lambda_{\text{max}}(n)\sim B\dfrac{e^{Kn}}{n}, &  & n\equiv
0\operatorname{mod}4, &  & \text{as }n\rightarrow\infty
\end{array}
\]
where%
\[
K=\ln(2)+G(x_{0})=0.1986176152...=\ln(1.2197154761...),
\]%
\[
B=\frac{2\sqrt{2}e^{K}}{\left(  4-e^{2K}\right)  ^{1/4}}=2.7402229903....
\]
See also \cite{Soae-err, BkUnc-err}.

\textbf{2.13. Prandtl--Blasius Flow.} \ The radius of convergence
$3.1273479155...$ appears in \cite{Shanks-err}: the associated power series is
divergent beyond this point and yet certain techniques exist for inducing convergence.

\textbf{3.12. Electrical Capacitance}. \ Lower bounds for so-called Euclidean
capacities include%
\[
\gamma_{1}(T)\geq\frac{1}{3}\cdot0+\frac{2}{3}\cdot1=\frac{2}{3},
\]%
\[
\gamma_{1}(S)\geq\frac{1}{4}\cdot0+\frac{1}{2}\cdot1+\frac{1}{4}\cdot\sqrt
{2}=\frac{\sqrt{2}}{4}\left(  1+\sqrt{2}\right)  ,
\]%
\[
\gamma_{1}(C)\geq\frac{1}{8}\cdot0+\frac{3}{8}\cdot1+\frac{3}{8}\cdot\sqrt
{2}+\frac{1}{8}\cdot\sqrt{3}=\frac{\sqrt{3}}{8}\left(  1+\sqrt{3}+\sqrt
{6}\right)  .
\]
The constant $2/3$ perhaps cannot be improved for the unit equilateral
triangle (nor can $3/4$ for the regular tetrahedron). \ We are less sure about
the square and cube.

\textbf{4.5. Variants of Brownian Motion.}

\textbf{4.6. Shapes of Binary Trees.} \ Janson \cite{Jans0-err} provided
\textquotedblleft a kind of dictionary\textquotedblright, to facilitate study
of the literature about Brownian excursion area moments and asymptotic
enumeration of connected graphs.

\textbf{4.16. Lyapunov Exponents IV.} \ For random singular matrices of the
form \cite{HlEq1-err, HlEq2-err}%
\[%
\begin{array}
[c]{ccccc}%
\left(
\begin{array}
[c]{cc}%
1 & x_{j}\\
1/x_{j} & 1
\end{array}
\right)  , &  & \text{independent }x_{j}\sim\text{Unif}[0,1]\text{,} &  &
j\geq1
\end{array}
\]
the Lyapunov exponent of their product is $2\ln(2)-1/2$. \ If instead the
underlying distribution is Unif$[-1,1]$, then $\ln(2)-1/2$ emerges.
\ Corresponding central limit variances (expressions in $\ln(2)$ and $\pi$)
are also known.

\textbf{4.23. Electing a Leader.} \ In the expression for $\operatorname*{E}%
(\tilde{V}_{n})$, the fraction $\pi^{2}/16$ is incorrect and should be
$\pi^{2}/8$ instead \cite{FGH5-err}. Thus, in the next sentence, the size
difference is $\pi^{2}/(8\ln(2))=1.7798536656...$, considerably larger than
the height difference. Here is another shooting problem (not to be confused
with riflery described in our entry [5.9]). We start with $n$ assassins in a
room.\ At each integer time $\geq1$, each surviving assassin fatally shoots a
randomly chosen surviving assassin (other than oneself; two people may
conceivably choose each other).\ Eventually we reach a state with either $0$
or $1$ survivors. Let $p_{n}$ denote the probability of $0$ survivors.
Simulation suggests that $p_{n}$ does not converge as $n$ increases, but is
asymptotically periodic on the $\ln(n)$ scale with period $1$ and
\cite{vdBKM-err}
\[%
\begin{array}
[c]{ccc}%
\lim\inf p_{n}\leq0.477449, &  & 0.515428\leq\lim\sup p_{n}%
\end{array}
\]
as $n\rightarrow\infty$.\ Quantitative performance of distributed algorithms
for transmitting data (specifically, in resolving conflicts, maximum finding
and list sorting) is surveyed in \cite{FGH5-err}.

\textbf{4.25.\ Biham--Middleton--Levine Traffic.} Among many models for
traffic, we mention an especially simple one: cars arrive one-by-one at a
stoplight (which alternates between red or green in blocks of $\ell$)
according to a Bernoulli($p$) process. What is the expected maximum line
length over a specified time period $[0,T]$? Such worst-case traffic
congestion analysis, for $\ell\geq2$ and $T\rightarrow\infty$, rests upon what
is known as the \textit{Poisson clumping heuristic}, an unproven yet
remarkably accurate predictive tool \cite{FLLM-err}.

\textbf{5.12. Random Triangles. III.} \ The expected product of two sides of a
random triangle in a unit disk is exactly \cite{Bck01-err, Bck02-err}%
\[
\operatorname*{E}(a\,b)=\frac{3383}{432\pi^{2}}+\frac{35}{96\pi^{2}}\zeta(3)
\]
from which
\[
\operatorname*{E}(\operatorname*{perimeter}\nolimits^{2})=3+\frac{3383}%
{72\pi^{2}}+\frac{35}{16\pi^{2}}\zeta(3)
\]
immediately follows.

\section{Acknowledgements}

There are too many people to thank...

\end{document}